\documentclass[11pt,lettersize,reqno,fullpage]{article}
\usepackage{amsmath,amssymb,amsfonts,mathrsfs,verbatim,enumitem,pstricks,amsthm}
\usepackage[all]{xy}
\usepackage{centernot, xr, accents}
\usepackage{hyperref}

\theoremstyle{plain}
\newtheorem{thm}{Theorem}[section]
\newtheorem{cor}[thm]{Corollary}
\newtheorem{lmm}[thm]{Lemma}
\newtheorem{prp}[thm]{Proposition}

\newtheorem{eg}[thm]{Example}
\newtheorem{defn}[thm]{Definition}

\newtheorem{claim}[thm]{Claim}
 
\newtheorem{rem}[thm]{Remark} 

\newtheorem{notn}[thm]{Notation}
\newtheorem*{mthm}{MAIN THEOREM}

\theoremstyle{definition}
\newtheorem{step}{Step}

\numberwithin{equation}{section}

\def \hf{\hspace*{0.5cm}}                      
\def\bge{\begin{equation}}                
\def\ede{\end{equation}}                
\def\bgd{\begin{displaymath}}         
\def\edd{\end{displaymath}}            
\def\bgee{\begin{equation*}}           
\def\edee{\end{equation*}}           

\def \ni{\noindent}

\def\lra{\longrightarrow}
\def\lrab{\dashrightarrow}
\def\lan{\langle}
\def\ran{\rangle}

\def\BA{\begin{eqnarray}}
\def\EA{\end{eqnarray}}
\def\BAA{\begin{eqnarray*}}
\def\EAA{\end{eqnarray*}}
\def\Bal{\begin{align*}}
\def\Eal{\end{align*}}

\def \C{\mathbb{C}}
\def \PD{\textnormal{PD}}
\def \P{\mathbb{P}}
\def \A{\mathcal{A}}
\def \E{\mathcal{E}}
\def \Num{\mathcal{N}}
\def \PP{\mathcal{P}}
\def \f{\frac}
\def \ff{\tilde{f}}

\def \p{\tilde{p}}

\def \lp{l_{\p}}

\def \w{\tilde{w}}

\def \B{\mathcal{B}}

\def \UL{\mathbb{L}} 
\def \DL{\mathcal{L}}
\def \UV{\mathbb{V}}
\def \DV{\mathcal{V}}

\def \ds{\psi}
\def \us{\Psi}

\def \ct{\mathcal{C}}

\def \y{y}

\def \a{a}
\def \lm{\lambda}
\def \gD{\gamma_{\mathcal{D}}}
\def \gP{\gamma_{_{\mathbb{P}^2}}}
\def \G{\tilde{\gamma}}

\def \25node{A_5} 
\def \62node{A_6}

\def \D{\mathcal{D}}
\def \DD{\mathcal{D}}

\def \X{\mathfrak{X}}

\def \pf{\noindent \textbf{Proof:  }}
\def \ov{\overline}

\def \ts{\zeta} 
\def \nts{\varphi}
\def \td{\PP D}

\def \XC{\mathcal{X}}

\def \AA{\hat{\A}}
\def \DD{\hat{\D}}
 
\def \XX{ \hat{\XC}}

\def \q{\rho}
\def \qq{f}

\def \tt{\tau}
\def \xx{\hat{x}}

\def \J{\mathcal{J}}
\def \U{\mathcal{U}}

\def \hati{\hat}
\def \hatii{\bar} 
\def \ha{\mathring}
\def \ho{\accentset{\star}}

\def \W{\mathbb{W}}
\def \WL{\mathcal{W}}
\def \Q{\mathcal{Q}}
\def \Si{\mathcal{\varphi}}

\def \N{\nabla}

\def \cu{\gamma}

\def \F{\mathcal{F}}

\def \s{s}

\def \GG{\mathrm{G}}
\def \HH{\mathrm{H}}
\def \T{\mathrm{P}}
\def \Z{\mathrm{A}} 
\def \Y{\mathrm{B}}
\def \g{\mathrm{g}}
\def \kap{\mathrm{C}}

\def \sq{s}
\def \rq{r} 
\def \kq{\kappa}
\def \etq{\eta}

\def \xq{\accentset{\star}{x}}
\def \yq{\accentset{\star}{y}} 
\def \Bq{\mathrm{B}}

\def \Hpl{\mathrm{H}} 
\def \Ln{\mathrm{L}}

\def \Bz{\mathrm{B}}
\def \Pz{\mathrm{P}}

\def \TT{\mathrm{T}}

\def \hxt2{\hat{x}_{t_2}} 
\def \hyt2{\hat{y}_{t_2}}

\def \SS{\mathcal{S}}

\title{{\LARGE Enumeration of curves with one singular point}}
\author{Somnath Basu and Ritwik Mukherjee }
\date{}
\usepackage{hyperref}
\hypersetup{
	colorlinks,
	citecolor=blue,
	filecolor=black,
	linkcolor=blue,
	urlcolor=black
}

\begin{document}

\maketitle

\begin{abstract}
In this paper we obtain an explicit formula for the number of curves in $\P^2$, of degree $d$, passing through $(d(d+3)/2 -k)$ generic points 
and having a codimension $k$ singularity, where $k$ is at most $7$. In the past, many of these numbers were computed using techniques from 
algebraic geometry. In this paper we use purely \textit{topological} methods to count curves. Our main tool is a classical fact from 
differential topology: the number of zeros of a generic smooth section of a vector bundle $V$ over $M$, counted with a sign, is 
the Euler class of $V$ evaluated on the fundamental class of $M$. 
\end{abstract}

\tableofcontents

\section{Introduction}
\label{introduction}
\hf\hf Enumerative geometry is a branch of mathematics concerned with the following question: 
\begin{center}
{\it How many geometric objects are there which satisfy prescribed constraints?  }
\end{center}
A well known class of enumerative problems is that of singular curves in $\P^2$ (complex projective space) passing through the appropriate number of points. This question has been studied by algebraic geometers for a long time. 
However, in this paper we use purely topological methods to tackle this problem.  \\    
\hf\hf Let us denote the space of curves of degree $d$ in $\P^2$ by $\D$. It follows that $\D \cong \P^{\delta_d}$, where $\delta_d = d(d+3)/2$. Let $\gP\lra \P^2$ be the tautological line bundle. 
A homogeneous polynomial $f$, of degree $d$ and in $3$ variables, induces a holomorphic section of the line bundle $\gP^{*d} \lra \P^2$. 
If $f$ is non-zero, then we will denote its \textit{equivalence class} in $\D$ by $\ff$. 
Similarly, if $p$ is a non-zero vector in $\C^3$, we will denote its equivalence class in $\P^2$ by $\p$ \footnote{In this paper we will use the symbol $\tilde{A}$ to denote the equivalence class of 
$A$ instead of the standard $[A]$. This will make some of the calculations in section \ref{closure_of_spaces} easier to read.}.     
\begin{defn}
\label{singularity_defn}
Let $\ff \in \D$ and $\p \in \P^2$. A point $\p \in f^{-1}(0)$ \textsf{is of singularity type} $\A_k$,
$\D_k$, $\E_6$, $\E_7$, $\E_8$ or $\XC_8$ if there exists a coordinate system
$(x,y) :(\U,\p) \lra (\C^2,0)$ such that $f^{-1}(0) \cap \U$ is given by 
\begin{align*}
\A_k: y^2 + x^{k+1}   &=0  \qquad k \geq 0, \qquad \D_k: y^2 x + x^{k-1} =0  \qquad k \geq 4, \\
\E_6: y^3+x^4 &=0,  \qquad \E_7: y^3+ y x^3=0, 
\qquad \E_8: y^3 + x^5=0, \\
\XC_8: x^4 + y^4 &=0.   
\end{align*}
\end{defn}
In more common terminology, $\p$ is a {\it smooth} point of $f^{-1}(0)$ if 
it is a singularity of type $\A_0$; a {\it simple node} if its singularity type is $\A_1$; a {\it cusp} if its type is $\A_2$; a {\it tacnode} if its type is $\A_3$; a {\it triple point} if its type is $\D_4$; and a {\it quadruple point} if its type is $\XC_8$. \\
\hf\hf We have several results (cf. Theorem \ref{algopa20}\,-\,\ref{algope7}, section \ref{algorithm_for_numbers})  which can be summarized collectively as our main result. Although \eqref{algopa20}-\eqref{algope7} may appear as equalities, the content of each of these equations is a theorem.
\begin{mthm}
\label{main_result}
Let $\X_k$ be a singularity of type $\A_k$, $\D_k$ or $\E_k$. Denote $\Num(\X_k,n)$ to be the number of degree $d$ curves in $\P^2$ that pass through $\delta_d - (k+n)$ generic points and have a singularity of type $\X_k$ at the intersection of $n$ generic lines. \\
(i) There is a formula for $\Num(\X_k,n)$ if $k \leq 7$, provided  $d \geq \ct_{\X_k}$ where 
\bgd
\ct_{\A_k} = k+1, ~~\ct_{\D_k} = k-1, ~~\ct_{\E_6} = 4, ~~\ct_{\E_7} = 4.
\edd
(ii) There is an algorithm to explicitly compute these numbers.     
\end{mthm}
\begin{rem}
Note that $\Num(\X_k,n)$ is zero if $n>2$, since three or more generic lines do not intersect anywhere.   
Moreover, $\Num(\X_k,2)$ is the the number of degree $d$-curves through $\delta_d-(k+2)$ generic points having 
one singularity of type $\X_k$ lying at a given fixed point (since the intersection of two generic lines 
is a point).
\end{rem}
The numbers $\Num(\X_k,0)$ till 
$k\leq 7$ have also been computed by Maxim Kazarian \cite{Kaz}  and Dimitry Kerner \cite{Ker1} 
using different methods. Our results for $n=0$ 
agree with theirs. The bound $d \geq \ct_{\X_k}$ is imposed to ensure that the 
relevant bundle sections are transverse.$\footnote{However, this bound is not the optimal bound.}$ 
The formulas for $\Num(\A_1,n)$, $\Num(\A_2,n)$ and $\Num(\A_3, n)$ also appear in \cite{Z1}.  
We extend the methods applied by the author to obtain the remaining formulas. 
This method carries over to the case of enumerating curves on any complex surface. 
With some further effort, the method can also be used to enumerate curves with more than one singular point. 
This will be pursued elsewhere. \\[0.2cm]

\hf\hf The second author is grateful to Aleksey Zinger for introducing him to the fascinating subject 
of enumerative geometry using topological methods and explaining $\cite{Z1}$. 
The contents of this paper were part of the second author's PhD thesis. He also thanks Jason Starr and Dusa McDuff for several useful discussions. 
The authors are grateful to Vamsi Pingali for answering several questions in complex geometry. 
And finally, the authors are indebted to Dennis Sullivan for sharing his deep and insightful perspective on this problem and indicating  
how this question is connected to different branches of mathematics. 


\section{Overview} 
\hf\hf Our main tool will be the following well known fact from topology (cf. \cite{BoTu}, Proposition 12.8).
\begin{thm} 
\label{Main_Theorem} 
Let $V\lra X$ be a vector bundle over a manifold $X$. Then the following are true: 
\hspace*{0.5cm}(1) A generic smooth section $s: X\lra V$ is transverse to the zero set. \\
\hspace*{0.5cm}(2) Furthermore, if $V$ and $X$ are oriented with $X$ compact then the zero set of such a section defines an integer homology class  in $X$, 
whose Poincar\'{e} dual is the Euler class of $V$. 
In particular, if the rank of $V$ is same as the dimension of $X$, 
then the signed cardinality of $s^{-1}(0)$ is the Euler class of $V$, evaluated on the fundamental class of 
$X$, i.e., 
\bgd
\pm |s^{-1}(0)| = \lan e(V), [X] \ran. 
\edd
\end{thm}
\begin{rem}
Let $X$ be a compact, complex manifold, $V$ a holomorphic vector bundle and $s$ a holomorphic section that is transverse to the zero set. If the rank of $V$ is same as the dimension of $X$, 
then the signed cardinality of $s^{-1}(0)$ is same as its actual cardinality (provided $X$ and $V$ have their natural orientations). 
\end{rem}
However, for our purposes, the requirement that $X$ is a smooth manifold is too strong. We will typically be dealing with spaces that are smooth but have non-smooth closure. The following result is a stronger  version of Theorem \ref{Main_Theorem}, that applies to singular spaces, provided the set of singular points is of real codimension two or more.   
\begin{thm} 
\label{Main_Theorem_pseudo_cycle} 
Let $M \subset \P^{N}$ be a smooth, compact algebraic variety and $X \subset M$ a smooth subvariety, not necessarily closed. Let $V \lra M$ be an oriented vector bundle, such that the rank of $V$ is same as the dimension of $X$. Then the following are true: \\
\hspace*{0.5cm}(1) The closure of $X$ is an algebraic variety and defines a homology class.\\
\hspace*{0.5cm}(2) The zero set of a generic smooth section $s: M \lra V$ intersects $X$ transversely and does not intersect $\ov{X}-X$ anywhere. \\
\hspace*{0.5cm}(3) The number of zeros of such a section inside $X$, counted with signs, 
is the Euler class of $V$ evaluated on the homology class $[\ov{X}]$, i.e., 
\bgd
\pm |s^{-1}(0) \cap \ov{X}| =  \pm |s^{-1}(0) \cap X| = \lan e(V), ~[\ov{X}] \ran.
\edd
\end{thm}

\begin{rem}
All the subsequent statements we make are true provided $d$ is sufficiently large. The precise bound on $d$ is given in section \ref{bundle_sections}.    
\end{rem}

\hf\hf We will now explain our strategy to compute  $\Num(\X_k, n)$. 
Given a singularity $\X_k$, 
let us also denote by $\X_k$, the \textit{space} of curves of degree $d$ with a marked point $\p$ such that the curve has a singularity of type $\X_k$ at $\p$, i.e.,
\bgd
\X_k := \{ (\ff,\p) \in \D\times \P^2: ~\textnormal{$f$ has a singularity of type $\X_k$ at the point $\p$}\}.
\edd
Let $\p_1, \p_2, \ldots , \p_{\delta_d -(k+n)}$ be 
$\delta_d-(k+n)$ generic points in $\P^2$ 
and $\Ln_1, \Ln_2, \ldots, \Ln_n $ be $n$ generic lines 
in $\P^2$. Define the following sets
\begin{align}
\label{hyperplane} 
\Hpl_i& := \{ \ff \in \D: f(p_i)=0 \}, \qquad \Hpl_i^* := \{ \ff \in \D: f(p_i)=0,  \nabla f|_{p_i} \neq 0 \}  \nonumber \\
\hat{\Hpl}_i& := \Hpl_i \times \P^2, \qquad \hat{\Hpl}_i^* := \Hpl_i^* \times \P^2 \qquad \textnormal{and}  \qquad \hat{\Ln}_i := \D \times \Ln_i.  
\end{align}
By definition, our desired number $\Num(\X_k,n)$ is the cardinality of 
the set
\bge
\label{number_Xk_defn}
\Num(\X_k,n) := |\X_k \cap \hat{\Hpl}_1 \cap \ldots \cap \hat{\Hpl}_{\delta_d-(n+k)} \cap 
\hat{\Ln}_1 \cap \ldots \cap \hat{\Ln}_n|.
\ede
\begin{step} 
If the degree $d$ is sufficiently large then the space $\X_k$ is a smooth algebraic variety and its closure defines a homology class.
\begin{lmm}{\bf (cf. section \ref{bundle_sections})}
The space $\X_k$ is a smooth subvariety of $\D \times \P^2$ of dimension $\delta_d-k$. 
\end{lmm}
\end{step}

\begin{step}
If the points and lines are chosen generically, then the corresponding hyperplanes and lines defined in \eqref{hyperplane} will intersect our space $\X_k$ transversely. 
Moreover, they won't intersect any extra points in the closure. 
\begin{lmm}
\label{gpl2}
Let $\p_1, \p_2, \ldots , \p_{\delta_d -(k+n)}$ be 
$\delta_d-(k+n)$ generic points in $\P^2$ 
and $\Ln_1, \Ln_2, \ldots, \Ln_n $ be $n$ generic lines 
in $\P^2$. Let $\hat{\Hpl}_i$, $\hat{\Hpl}_i^*$  and $\hat{\Ln}_i$ be as defined in 
\eqref{hyperplane}.  
Then 
\bgd
\ov{\X}_k \cap \hat{\Hpl}_1 \cap \ldots \cap \hat{\Hpl}_{\delta_d -(k+n)} 
\cap \hat{\Ln}_1 \cap \ldots \cap \hat{\Ln}_n = 
\X_k \cap \hat{\Hpl}_1^* \cap \ldots \cap \hat{\Hpl}_{\delta_d -(k+n)}^* 
\cap \hat{\Ln}_1 \cap \ldots \cap \hat{\Ln}_n  
\edd
and every intersection is transverse. 
\end{lmm}
Although we omit the details of the proof, this follows from an application of the families transversality theorem and Bertini's theorem. 
The details of this proof  can be found in \cite{BM_Detail}. 
 
\begin{notn}
\label{tau_bundle_defn}
Let $\gD\lra \D$ and $\gP\lra\P^2$ denote the tautological line bundles. If $c_1(V)$ denotes the first Chern class of a vector bundle then we set 
\bgd
\y : = c_1(\gD^*) \in H^{2}(\D; \mathbb{Z}), \qquad 
\a := c_1(\gP^*)  \in H^{2}(\P^2; \mathbb{Z}). 
\edd
\end{notn}
  
\hf\hf As a consequence of Lemma \ref{gpl2}  
we obtain the following fact: 
\begin{lmm}
\label{gpl}
The number $\Num(\X_k,n)$ is given by~  
$ \Num(\X_k, n) = \lan \y^{\delta_d-(n+k)} \a^n, ~[\ov{\X}_k] \ran$. 
\end{lmm}
\end{step}
\pf This follows from  Theorem \ref{Main_Theorem_pseudo_cycle} and Lemma \ref{gpl2}. \qed 

\begin{rem}
\label{an1_again}
Here we are making an abuse of notation 
by referring to $\y, \a \in H^*(\D \times \P^2; \mathbb{Z})$. 
The intended meaning is $\pi_{\D}^* \y$ and 
$\pi_{\P^2}^*\a$, where 
$\pi_{\D}, \pi_{\P^2}: \D \times \P^2 \lra \D, \P^2$ 
are the projection maps. We will make a similar abuse of notation with vector bundles. Our intended meaning should be clear when we say, for instance, $\gD^* \lra \D \times \P^2$.
\end{rem}

\ni The space $\X_k$, unfortunately, is not easy to describe directly. Consequently, computing $\Num(\X_k,n)$ \textit{directly} is not a promising approach. Instead we will look at the space 
$$\PP \X_k \subset \D \times \P T\P^2.$$ 
This is the space of curves $\ff$, of degree $d$, with a marked point $\p \in \P^2$
and a marked direction $\lp \in \P T_{\p}\P^2$, such that the curve $f$ has a singularity of type $\X_k$ at $\p$ and certain directional derivatives \textit{vanish along $\lp$}, and 
certain other derivatives \textit{don't vanish}. To take a simple example, $\PP \A_2$ is the space of curves $\ff$ with a marked point $\p$ and a 
marked direction $\lp$ such that $f$ has an $\A_2$-node at $\p$ and the Hessian is degenerate along $\lp$, 
but the third derivative along $\lp$ is non-zero. It turns out that this space is much easier to describe. 
The precise definition of the space $\PP \X_k$ is given in subsection \ref{definition_of_px}. 

\begin{step}
Since the space $\PP \X_k$ is described locally as the vanishing of certain sections that are transverse to the zero set they are \text{smooth} algebraic varieties.

\begin{lmm}{\bf (cf. section \ref{bundle_sections})}
\label{pr_sp_pseudo}
The space 
$\PP \X_k$ is a smooth subvariety 
of $\D \times \P T\P^2$ of dimension $\delta_d-k$. 
\end{lmm}

\begin{notn}
\label{tau_bundle_pv_defn}
Let $\G \lra \P T\P^2$ be the tautological line bundle. The first Chern class of the dual will be denoted by $\lm = c_1(\G^*)\in H^{2}(\P T\P^2; \mathbb{Z})$.   
\end{notn}
\ni Lemma \ref{pr_sp_pseudo} now motivates the following definition: 
 \begin{defn}
\label{up_number_defn}
 We define the number 
$\Num(\PP \X_k,n,m)$ as 
\bge
\label{num_proj}
\Num(\PP \X_k, n,m) := \lan \y^{\delta_d-(k+n+m)} \a^n \lm^m, 
~[\ov{\PP \X}_k] \ran. 
\ede
\end{defn}
\ni The next Lemma relates the numbers $\Num(\PP\X_k,n,0)$ and 
$\Num(\X_k,n)$. 
\begin{lmm}{\bf (cf. section \ref{appendix_proof_of_Lemmas})}
\label{up_to_down}
The projection map $ \pi: \PP \X_k \lra \X_k $ is one to one if $\X_k = \A_k, \D_k, \E_6, \E_7$ or $\E_8$ except for $\X_k = \D_4$ when it is three to one. In particular, 
\bge
\label{up_down_equation}
\Num(\X_k, n)  = \Num(\PP \X_k, n,0) 
\qquad \textnormal{if} ~~\X_k \neq \D_4   \qquad \textnormal{and} 
\qquad \Num(\D_4, n) = 
\frac{\Num(\PP \D_4,n,0)}{3}.
\ede
\end{lmm}
\end{step}
\ni To summarize, the \textit{definition} of $\Num(\X_k,n)$ is \eqref{number_Xk_defn}. Lemma \ref{gpl} equates this number to a \textit{topological} computation. 
We then introduce another number $\Num(\PP\X_k,n,m)$ in definition \ref{up_number_defn} and relate it to $\Num(\X_k,n) $ in Lemma \ref{up_to_down}. In other words, we do not compute $\Num(\X_k,n)$ \textit{directly}. 
We compute it \textit{indirectly} by first computing $\Num(\PP \X_k, n,m)$ and then using Lemma \ref{up_to_down}. \\
\hf\hf We know give a brief idea of how to compute these numbers. Suppose we want to compute $\Num(\PP \X_k, n, m)$. 
We first find some singularity $\X_l$ for which $\Num(\PP \X_l, n,m)$ has been calculated and which contains $\X_k$ in its closure, i.e., we want $\PP \X_k$ to be a subset of $\ov{\PP \X}_{l}$. 
Usually, $l=k-1$ but it is not necessary. Our next task is to describe the closure of $\PP \X_{l}$ explicitly as
\bge
\label{stratification_general}
\ov{\PP \X}_{l} = \PP \X_{l} \cup \ov{\PP\X}_k \cup \B.
\ede
Equivalently, we want an explicit description of the space $\B$. By definition \ref{up_number_defn} and Theorem \ref{Main_Theorem_pseudo_cycle} 
\bgd
\Num(\PP \X_k, n, m) := \lan e(\W_{n,m,k}), ~[\ov{\PP \X}_k] \ran = \pm|\Q^{-1}(0) \cap \PP \X_k|, 
\edd
where
\bge
\label{generic_Q}  
\Q :\D \times \P T\P^2  \lra \W_{n,m,k}:= \bigg(\bigoplus_{i=1}^{\delta_d -(n+m+k)}\gD^*\bigg)\oplus\bigg(\bigoplus_{i=1}^{n} \gP^*\bigg)\oplus\bigg(  \bigoplus_{i=1}^{m}\G^*  \bigg)
\ede
is a generic smooth section. We now have to construct a section $\us_{\PP \X_k}$ 
of some vector bundle\footnote{The Euler class of this vector bundle is expressible in terms of the Euler classes of three canonical line bundles via the splitting principal.} 
\bgd
\UV_{\PP \X_k}\lra \ov{\PP \X}_{l} = \PP \X_{l} \cup \ov{\PP\X}_k \cup \B
\edd
with the following properties: it should not vanish on $\PP \X_{l}$ and it should vanish \textit{transversely} on $\PP \X_k$. In that case we are led to
\bgd
\lan e(\UV_{\PP \X_{k}} \oplus \W_{n,m,k}), [\ov{\PP \X}_{l}] \ran 
= \Num(\PP \X_k,n,m) + \mathcal{C}_{\B}(\us_{\PP \X_k} \oplus \Q),
\edd
where $\mathcal{C}_{\B}(\us_{\X_k} \oplus \Q)$ is the contribution of the section 
$\us_{\PP \X_k} \oplus \Q$ to the Euler class from the points of $\B$. The left hand side is computable via splitting principle and the fact that $\Num(\PP \X_l, n,m)$ is known. Therefore, once we know $\mathcal{C}_{\B}(\us_{\X_k} \oplus \Q)$, we get a recursive formula for the number $\Num(\PP \X_k,n,m)$ and iterate.   
\begin{eg}
Suppose we wish to compute $\Num(\A_5, n)$. This can be deduced from the knowledge of $\Num(\PP \A_5, n,m)$. The obvious singularities which have $\A_5$-nodes in its closure are $\A_4$-nodes. 
In order to analyze the space $\ov{\PP \A}_4$, we infer that (cf. Lemma \ref{cl}, statement \ref{A4cl})
\bgd
\ov{\PP \A_4} = \PP \A_4 \cup \ov{\PP \A}_5 \cup \ov{\PP \D}_5.
\edd
The corresponding line bundle $\UL_{\PP \A_5} \lra \ov{\PP \A}_4$ with a section $\us_{\PP \A_5}$ that does not vanish on $\PP \A_4$ and vanishes transversely on $\PP \A_5$ is defined in subsection \ref{summary_vector_bundle_definitions}. The verification of these properties of the section is proved in section \ref{bundle_sections} (Proposition \ref{Ak_Condition_prp}).  
Finally, in Corollary \ref{mult_of_pa5_section_around_pd5} we show that if $\Q$ is a generic section of the vector  bundle 
\bgd
\W_{n,m,5} \lra \D \times \P T\P^2 
\edd
then $\us_{\PP \A_5} \oplus \Q $ vanishes on all the points of $\PP \D_5$ with a multiplicity of $2$. 
Hence, we conclude that
\bge
\label{sample_computation}
\lan e(\UL_{\PP \A_{5}} \oplus \W_{n,m,5}), ~[\ov{\PP \A}_{4}] \ran 
= \Num(\PP \A_5, n, m) + 2\Num(\PP \D_5, n, m).
\ede
This gives us a recursive formula for $\Num(\PP \A_5, n, m)$ 
in terms of $\Num(\PP \A_4, n^{\prime}, m^{\prime})$ and $\Num(\PP \D_5, n, m)$ 
which is \eqref{algopa5} in our algorithm. 
\end{eg}
\hf\hf Now we describe the basic organization of our paper. In 
section \ref{algorithm_for_numbers} we state the explicit algorithm to obtain the numbers $\Num(\X_k, n)$ in our MAIN THEOREM in section \ref{introduction}. 
In section \ref{summary_notation_def} we summarize all the spaces, vector bundles 
and sections of vector bundles we will encounter in the  course of our computations. 
In section \ref{condition_for_sing} we describe necessary and sufficient conditions 
for a curve $f^{-1}(0)$ to have a singularity of type $\X_k$ at a point. In section \ref{bundle_sections} we describe the spaces $\PP \X_k$ as the vanishing of certain sections and the non-vanishing of certain other sections. Moreover, we show that these sections are transverse to the zero set. In section \ref{closure_of_spaces} we stratify the space 
$\ov{\PP \X}_k$ as described in \eqref{stratification_general}. Along the way we also compute the \textit{order} to which a certain section vanishes around certain points (i.e., the contribution of the section to the Euler class of a bundle). Finally, using the splitting principal, in section \ref{Euler_class_computation} we compute the Euler class of the relevant bundles and obtain the recursive formula similar to \eqref{sample_computation} above.               


\section{Algorithm}  
\label{algorithm_for_numbers}
\hf\hf We now give an algorithm to compute the numbers $\Num(\X_k,n)$. 
We have implemented this algorithm in a Mathematica program to obtain the final answers. 
The program is available on our web page 
\url{https://www.sites.google.com/site/ritwik371/home}.  We prove these formulas in section \ref{Euler_class_computation}.  \\
\hf\hf The base case for the recursion is:  
\begin{align} 
\label{algoa1}
\Num(\A_1,n) = \begin{cases}
3(d-1)^{2},&\textnormal{if}~n=0;\\
3(d-1),&\textnormal{if}~n=1;\\
1,&\textnormal{if}~n=2;\\
0,&\textnormal{otherwise}.
\end{cases}
\end{align}
\hf\hf Next we will give an algorithm to compute $\Num(\PP \X_k, n,m)$. 
Using Lemma \ref{up_to_down} we get our desired numbers $\Num(\X_k,n)$. We note that 
using the ring structure of $H^*(\D \times \P T\P^2; \mathbb{Z})$, it is easy to see that 
for every singularity  type $\X_k$ we have
\begin{align}
\Num(\PP \X_k,n,m)
& = -3\Num(\PP \X_k,n+1,m-1)
-3\Num(\PP \X_k,n+2,m-2)   \qquad\forall ~~m\ge2. \label{ringp}
\end{align}
This follows from Lemma \ref{total_chern_class_of_tpn} and \ref{cohomology_ring_of_pv}.
Finally, we give recursive formulas for $\Num(\PP \X_k, n, m)$:    
\begin{eqnarray}
\Num(\PP \A_2, n, 0) & = & 2 \Num (\A_1,n) + 2(d-3) \Num(\A_1,n+1) \label{algopa20}\\
\Num(\PP \A_2, n, 1) & = & \Num(\A_1,n) + (2d-9) \Num(\A_1,n+1) + (d^2-9d+18) \Num(\A_1,n+2) \label{algopa21}\\
\Num(\PP \A_3, n, m) & = & \Num(\PP \A_2, n, m ) +  3\Num(\PP \A_2, n, m+1) + d\Num(\PP \A_2, n+1, m) \label{algopa3}\\ 
\Num( \PP \A_4, n, m) & =  & 2\Num( \PP \A_3, n, m) + 2\Num(\PP \A_3, n, m+1) + (2d-6)\Num( \PP \A_3, n+1, m) \label{algopa4} \\  
\Num(\PP \A_5, n, m) & = & 3\Num(\PP \A_4, n, m) + \Num(\PP \A_4, n, m+1) + (3d -12)\Num( \PP \A_4, n+1, m)\nonumber\\ 
&     & - 2 \Num(\PP \D_5, n, m) \label{algopa5} \\  
\Num( \PP \A_6, n, m) & = & 4\Num(\PP \A_5, n, m) +0\Num(\PP \A_5, n, m+1) +  (4d -18)\Num(\PP \A_5, n+1, m) \nonumber \\   
&     & -4\Num(\PP \D_6, n, m) - 3\Num(\PP \E_6, n, m) \label{algopa6}\\
\Num(\PP \A_7, n, 0) & = & 5\Num(\PP \A_6, n, 0) -\Num(\PP \A_6, n, 1) + (5d-24)\Num(\PP \A_6, n+1, 0) \nonumber \\  
&     & -6\Num(\PP \D_7,n,0) -7\Num(\PP \E_7, n, 0) \label{algopa7}  \\
\Num(\PP \D_4, n, m) & = & \Num(\PP \A_3, n, m) -2\Num(\PP \A_3, n, m+1) + (d-6)\Num(\PP \A_3, n+1, m) \label{algopd4} \\ 
\Num(\PP \D_5, n, m) & = & \Num(\PP \D_4, n, m) + \Num(\PP \D_4, n, m+1) + (d-3)\Num(\PP \D_4, n+1, m) \label{algopd5} \\
\Num(\PP \D_6, n, m) & = & \Num(\PP \D_5, n, m) + 4\Num(\PP \D_5, n, m+1) + d\Num(\PP \D_5, n+1, m) \label{algopd6} \\
\Num(\PP \D_7, n, m) & = & 2\Num(\PP \D_6, n, m) + 4\Num(\PP \D_6, n, m+1) + (2d-6)\Num(\PP \D_6, n+1, m) \label{algopd7} \\
\Num(\PP \E_6, n, m) & = & \Num(\PP \D_5, n, m) -\Num(\PP \D_5, n, m+1) + (d-6)\Num(\PP \D_5, n+1, m) \label{algope6} \\
\Num(\PP \E_7, n, m) & = & \Num(\PP \D_6, n, m) -\Num(\PP \D_6, n, m+1) +  (d-6)\Num(\PP \D_6, n+1, m) \label{algope7}
\end{eqnarray}


\section{Summary of definitions and notation}  
\label{summary_notation_def}
 
\subsection{The vector bundles involved}
\label{summary_vector_bundle_definitions}

\hf\hf We now list down all the vector bundles that we will encounter. 
The first three of these, the tautological line bundles, have been defined in notations 
\ref{tau_bundle_defn} and \ref{tau_bundle_pv_defn}. Let $\pi:\D\times \P T\P^2\lra \D\times\P^2$ be the projection map. We have the following bundles over $\D\times\P^2$ :
\begin{eqnarray*}
\DL_{\A_0} &:= & \gD^*\otimes \gP^{*d} \lra \D \times \P^2 \\
\DV_{\A_1} &:= & \gD^*\otimes \gP^{*d} \otimes T^*\P^2 \lra \D \times \P^2 \\
\DL_{\A_2} &:= & (\gD^* \otimes \gP^{*d} \otimes \Lambda^2 T^*\P^2)^{\otimes 2} 
\lra \D \times \P^2 \\ 
\DV_{\D_4} &:= & \gD^*\otimes \gP^{*d} \otimes 
\textnormal{Sym}^2 (T^*\P^2 \otimes T^*\P^2) \lra \D \times \P^2 \\
\DV_{\XC_8} &:= &  \gD^*\otimes \gP^{*d} \otimes 
\textnormal{Sym}^3 (T^*\P^2 \otimes T^*\P^2 \otimes T^*\P^2) 
\lra \D \times \P^2 
\end{eqnarray*}
Associated to the map $\pi$ there are pullback bundles 
\begin{eqnarray*}
\UL_{\AA_0} &:= &\pi^* \DL_{\A_0} \lra \D \times \P T\P^2 \\
\UV_{\AA_1} &:= & \pi^{*} \DV_{\A_1} \lra \D \times \P T\P^2 \\
\UV_{\DD_4} &:= & \pi^{*} \DV_{\D_4} \lra \D \times \P T\P^2 \\
\UV_{\XX_8} &:= & \pi^{*} \DV_{\XC_8} \lra \D \times \P T\P^2 \\
\UV_{\PP \A_2} &:= & \G^*\otimes \gD^*\otimes \gP^{*d} \otimes \pi^* T^*\P^2 
\lra  \D \times \P T\P^2 \\
\UV_{\PP \D_5} &:= & \G^{*2}\otimes \gD^*\otimes \gP^{*d} \otimes \pi^* T^*\P^2 
\lra  \D \times \P T\P^2.
\end{eqnarray*}
Finally, we have 
\begin{align*}
\UL_{\PP \D_4} &:= (T\P^2/\G)^{*2} \otimes \gD^* \otimes \gP^{*d} 
 \lra \D \times \P T\P^2  \\ 
 \UL_{\PP \D_5} &:= \G^{*2} \otimes (T\P^2/\G)^* \otimes 
 \gD^* \otimes \gP^{*d} \lra \D \times \P T\P^2 \\
\UL_{\PP \D_5^{\vee}} &:= \G^{*2} \otimes (T\P^2/\G)^{*4} \otimes 
 \gD^{* 2} \otimes \gP^{*2d} \lra \D \times \P T\P^2 \\
\UL_{\PP \D_6^{\vee}} &:= \G^{*8} \otimes (T\P^2/\G)^{*4} \otimes 
 \gD^{* 5} \otimes \gP^{*5d} \lra \D \times \P T\P^2 \\
 \UL_{\PP \E_6} &:= \G^{*} \otimes (T\P^2/\G)^{*2} \otimes 
 \gD^* \otimes \gP^{*d} \lra \D \times \P T\P^2\\
\UL_{\PP \E_7} &:=  \G^{*4} \otimes \gD^* \otimes \gP^{*d}
\lra \D \times \P T\P^2 \\
\UL_{\PP \E_8} &:=  \G^{*3} \otimes (T\P^2/\G)^* \otimes  \gD^* \otimes \gP^{*d}
\lra \D \times \P T\P^2 \\
\UL_{\PP \XC_8} &:= (T\P^2/\G)^{*3} \otimes \gD^* \otimes \gP^{*d} 
\lra \D \times \P T\P^2 \\
\UL_{\J} &:= \G^{*9}\otimes(T\P^2/\G)^{*3}\otimes \gD^{*3} \otimes \gP^{*3d} 
\lra \D \times \P T\P^2 \\
k \geq 3 \qquad \UL_{\PP \A_k} &:= \G^{*k} \otimes (T\P^2/\G)^{*(2k-6)} \otimes \gD^{*(k-2)} \otimes \gP^{*(d(k+1)-3d)}  
\lra \D \times \P T\P^2 \\  
k \geq 6 \qquad \UL_{\PP \D_k} & := \G^{*(k-2 + \epsilon_k)} \otimes (T\P^2/\G)^{*(2\epsilon_k)} \otimes \gD^{*(1+\epsilon_k)} \otimes \gP^{*(d(1+ \epsilon_k))} \lra \D \times \P T\P^2,
\end{align*}
where $\epsilon_6 =0$, $\epsilon_7=1$ and $\epsilon_8 =3$. 
In general, $\epsilon_k$ is the order of the pole of the section $\D^f_k$ at $f_{12} =0$.
The algorithm to obtain $\D^f_k$ for any $k$ is given in 
Lemma \ref{fstr_prp_Dk}.    
The reason for defining these bundles will become clearer in the subsection \ref{summary_sections_of_vector_bundle_definitions}, when we define 
sections of these bundles.\\

With the abuse of notation as explained in remark \ref{an1_again}, the bundle $T\P^2/\G$ is the quotient of the bundles $V$ and $W$, 
where $V$ is the pullback of the tangent bundle $T\P^2\to \P^2$ via $\D\times\P T\P^2\stackrel{\pi}{\rightarrow}\D\times \P^2\to \P^2$ and $W$ 
is pullback of $\G\to \P T \P^2$ via $\D\times \P T \P^2\to \P T \P^2$.

\subsection{Sections of Vector Bundles}
\label{summary_sections_of_vector_bundle_definitions}

\hf\hf Let us define the notion of \textit{vertical derivatives}. 
\begin{defn}
\label{vertical_derivative_defn}
\ni Let $\pi:V\longrightarrow M$ be a holomorphic vector bundle of rank $k$ and $s:M\lra V$ be a holomorphic section. Suppose $h: V|_{\U} \longrightarrow \U\times \mathbb{C}^{k}$ is a holomorphic trivialization of $V$ and $\pi_{1}, \pi_{2}: \U \times \C^{k} \longrightarrow \U, \C^{k}$ the projection maps.
Let 
\begin{align}
\label{section_local_coordinate}
\hati{s}&:= \pi_{2} \circ h \circ s.
\end{align}
For $q\in \U$, we define the {\it vertical derivative} of $s$ to be the $\C$-linear map  
\begin{align*}
\N s|_{q}: T_{q}M \longrightarrow V_{q}, \qquad  
\N s|_{q} & := (\pi_{2} \circ h)|_{V_{q}}^{-1} \circ d \hati{s}|_q,
\end{align*}
where $V_{q} = \pi^{-1}(q)$, the fibre at $q$. In particular, if $v \in T_q M$  
is given by a 
holomorphic map $\cu: \mathrm{B}_{\epsilon}(0) \lra M$ such that $\cu (0)=q$ and $\frac{\partial \cu}{\partial z}\big|_{z=0} = v$, then
\begin{align*} 
\nabla s|_q (v) & := (\pi_{2} \circ h)|^{-1}_{V_{q}} \circ 
\frac{\partial \hati{s}(\cu(z))}{\partial z}\bigg|_{z=0}
\end{align*}
were $\mathrm{B}_{\epsilon}$ is an open $\epsilon$-ball in $\C$ around the origin.$\footnote{Not every tangent vector is given by a holomorphic map; 
however combined with 
the fact that $\nabla s|_p $ is $\C$-linear, this definition determines $\nabla s|_p$ completely.} $ 
Finally, if $v, w \in T_q M$  are tangent vectors such that  
there exists a family of complex curves  $\cu: \mathrm{B}_\epsilon\times \mathrm{B}_\epsilon \lra M$ such that 
\bgd
\cu (0,0) =q, \qquad 
 \frac{\partial \cu(x,y)}{\partial x}\bigg|_{(0,0)}= v, \qquad 
\qquad \frac{\partial \cu(x,y)}{\partial y}\bigg|_{(0,0)} = w
\edd
then  
\begin{align}
\label{verder}
\nabla^{i+j} s|_q 
(\underbrace{v,\cdots v}_{\textnormal{$i$ times}}, \underbrace{w,\cdots w}_{\textnormal{$j$ times}}) & :=  
(\pi_{2} \circ h)\mid^{-1}_{V_{q}} \circ\left[
\frac{ \partial^{i+j} \hati{s}(\cu(x,y))}{\partial^i x \partial^j y}\right]\bigg|_{(0,0)}.
\end{align}
\end{defn}
\begin{rem}
In general the quantity in \eqref{verder} is not well defined, i.e., it depends on the trivialization and the curve $\cu$. Lemma \ref{wdc} explains on what subspace this quantity is well defined. 
\end{rem}
\begin{rem}
\label{transverse_local}
The section $ s: M \lra V $  is transverse to the zero set if and only if the induced map
\begin{align} 
\label{section_local_coordinate_calculus}
\hatii{s} & := \hati{s} \circ \varphi_{\U}^{-1} : \C^m \lra \C^k 
\end{align}
is transverse to the zero set in the usual calculus sense, where $\varphi_{\U}: \U \lra \C^m $ is a coordinate chart and $\hati{s}$ is as defined in \eqref{section_local_coordinate}.
\end{rem}
\hf\hf Let $f:\P^2 \lra \gP^{*d}$ be a section and $\p \in \P^2$. 
We can think of $p$ as a non-zero vector in $\gP$ and $p^{\otimes d}$ a non-zero vector 
in $\gP^{\otimes d}$ $\footnote{Remember that $p$ is an element of $\C^3-0$ while $\p$ is the corresponding equivalence class in $\P^2$}$. 
The quantity $\nabla f|_{\p}$ acts on a vector in  $\gP^{d}|_{\p}$ and produces an element of $T^*_{\p}\P^2$ . Let us denote this quantity as $\nabla f|_p$, i.e., 
\begin{align}
\nabla f|_p &:= \{\nabla f|_{\p}\}(p^{\otimes d}) \in T^*_{\p}\P^2.   
\end{align}
Notice that  $\nabla f|_{\p}$ is an element of the fibre of $T^*\P^{2} \otimes \gP^{*d}$ at $\p$ while $\nabla f|_{p}$ is an element of $T^*_{\p}\P^{2}$. \\
\hf\hf Now observe that $\pi^{*} T\P^2 \cong \G \oplus \pi^*T\P^2/\G \lra \P T\P^2$, where 
$\pi: \P T\P^2 \lra \P^2~$ is the projection map. Let us denote a vector in $\G$ by $v$ and  a vector 
in $\pi^*T \P^2/\G$ by $\w$. 
Given $\ff \in \D$ and $\p \in \P^2$, let 
\begin{align}
\label{abbreviation}
f_{ij} & := \nabla^{i+j} f|_p 
(\underbrace{v,\cdots v}_{\textnormal{$i$ times}}, \underbrace{w,\cdots w}_{\textnormal{$j$ times}}).
\end{align}
Note that $f_{ij}$ is a \textit{number}. 
In general $f_{ij}$ is not well defined; it depends on the trivialization and  the curve. 
Moreover it is also not well defined on the quotient space. 
Since our sections are not defined on the whole space, 
we will use the notation $ s:M \lrab V$ to indicate that $s$ is defined only on 
a subspace of $M$. With this terminology, we now explicitly define the sections that we will 
encounter in this paper.
\begin{eqnarray*}
\ds_{\A_0}: \D \times \P^2 \lra \DL_{\A_0}, & & \{\ds_{\A_0}(\ff, \p)\} (f \otimes p^{\otimes d}):= f(p) \\
\ds_{\A_1}: \D \times \P^2 \lrab \DV_{\A_1}, & & \{\ds_{\A_1}(\ff, \p)\}(f \otimes p^{\otimes d}):= \nabla f|_p  \\
\ds_{\D_4}: \D \times \P^2 \lrab \DV_{\D_4}, & & \{\ds_{\D_4}(\ff, \p)\}(f\otimes p^{\otimes d}):= \nabla^2 f|_p  \\ 
\ds_{\XC_8}: \D \times \P^2 \lrab \DV_{\XC_8}, & & \{\ds_{\XC_8}(\ff, \p)\}(f \otimes p^{\otimes d}):= \nabla^3 f|_p  \\ 
\ds_{\A_2} : \D \times \P^2 \lrab \DL_{\A_2}, & & \{\ds_{\A_2}(\ff, \p)\}(f \otimes p^{\otimes d}):= \textnormal{det}\, \nabla^2 f|_p \\
\us_{\AA_0}:  \D \times \P T\P^2 \lrab \UL_{\AA_0}, & & \us_{\AA_0}(\ff, \lp):= \ds_{\A_0}(\ff,\p)  \\
\us_{\AA_1}: \D \times \P T\P^2 \lrab \UV_{\AA_1}, & & \us_{\AA_1}(\ff,\lp):= \ds_{\A_1}(\ff, \p) \\
\us_{\DD_4}: \D \times \P T\P^2 \lrab \UV_{\DD_4}, & & \us_{\DD_4}(\ff,\lp):= \ds_{\D_4}(\ff, \p) \\
\us_{\XX_8}: \D \times \P T\P^2 \lrab \UV_{\XX_8}, & & \us_{\XX_8}(\ff,\lp):= \ds_{\XC_8}(\ff, \p).
\end{eqnarray*}
We also have
\begin{eqnarray*}
\us_{\PP \A_2}: \D \times \P T\P^2 \lrab \UV_{\PP \A_2}, & & \{\us_{\PP \A_2}(\ff,\lp)\}(f \otimes p^{\otimes d} \otimes v):= \nabla^2 f|_p (v,\cdot) \\  
\us_{\PP \D_5}^{\UV}: \D \times \P T\P^2 \lrab \UV_{\PP \D_5}, & & \{\us_{\PP \D_5}^{\UV}(\ff,\lp)\}(f \otimes p^{\otimes d} \otimes v^{\otimes 2}):= \nabla^3 f|_p (v,v,\cdot)  \\
\us_{\PP \D_4}: \D \times \P T\P^2 \lrab \UL_{\PP \D_4}, & &  \{\us_{\td_4}(\ff,\lp)\}(f \otimes p^{\otimes d} \otimes w^{\otimes 2}):= f_{02} \\
\us_{\PP \D_5}^{\UL}: \D \times \P T\P^2 \lrab \UL_{\PP \D_5}, & & \{\us_{\PP \D_5}^{\UL}(\ff,\lp)\}(f \otimes p^{\otimes d} \otimes v^{\otimes 2}\otimes w):= f_{21} \\
\us_{\PP \E_6}: \D \times \P T\P^2 \lrab \UL_{\PP \E_6}, & & \{\us_{\PP \E_6}(\ff,\lp)\}(f \otimes p^{\otimes d} \otimes v\otimes w^{\otimes 2}):= f_{12} \\ 
\us_{\PP \E_7}: \D \times \P T\P^2 \lrab \UL_{\PP \E_7}, & & \{\us_{\PP \E_7}(\ff,\lp)\}(f \otimes p^{\otimes d} \otimes v^{\otimes 4}):= f_{40} \\
\us_{\PP \E_8}: \D \times \P T\P^2 \lrab \UL_{\PP \E_8}, & & \{\us_{\PP \E_8}(\ff,\lp)\}(f \otimes p^{\otimes d} \otimes v^{\otimes 3} \otimes w ):= f_{31} \\
\us_{\PP \XC_8}: \D \times \P T\P^2 \lrab \UL_{\PP \XC_8}, & & \{\us_{\PP \XC_8}(f,\lp)\}(f \otimes p^{\otimes d} \otimes w^{\otimes 3}):= f_{03}.
\end{eqnarray*}
We also have sections of the following bundles:  $\us_{\PP \D_5^{\vee}}: \D \times \P T\P^2 \lrab \UL_{\PP \D_5^{\vee}}$ given by
\bge
\label{psi_d5_dual}
\{\us_{\PP \D_5^{\vee}}(\ff,\lp)\}( f^{\otimes 2} \otimes p^{\otimes 2 d} 
\otimes v^{\otimes 2} \otimes w^{\otimes 4}) := 3 f_{12}^2 - 4 f_{21} f_{03}, 
\ede
and  $\us_{\PP \D_6^{\vee}}: \D \times \P T\P^2 \lrab \UL_{\PP \D_6^{\vee}}$ given by
\bge
\label{psi_d6_dual}
\{\us_{\PP \D_6^{\vee}}(\ff,\lp)\}( f^{\otimes 5} \otimes p^{\otimes 5 d} 
\otimes v^{\otimes 8} \otimes w^{\otimes 4}) := \Big(f_{12}^4 f_{40} - 8f_{12}^3 f_{21} f_{31} + 24 f_{12}^2 f_{21}^2 f_{22}-32 f_{12} f_{21}^3 f_{13} + 16 f_{21}^4 f_{04} \Big),
\ede
and $\us_{\J}: \D \times \P T\P^2 \lrab \UL_{\J}$ given by
\bge
\label{J_psi}
\{\us_{\J}(\ff,\lp)\}( f^{\otimes 3} \otimes p^{\otimes 3d} \otimes v^{\otimes 9} \otimes w^{\otimes 3}) := 
\Big(- \frac{f_{31}^3}{8 } + \frac{3 f_{22} f_{31} f_{40}}{16 } - \frac{f_{13} f_{40}^2}{16}  \Big).
\ede

When $k\geq 3$ we have $\us_{\PP \A_k}: \D \times \P T\P^2 \lrab \UL_{\PP \A_k}$ given by
\bgd
\{\us_{\PP \A_k}(\ff,\lp) \}(f^{\otimes (k-2)} \otimes p^{\otimes d(k-2)} \otimes v^{\otimes k} \otimes w^{\otimes (2k-6)}):= f_{02}^{k-3} \A^f_k.
\edd
Similarly, when $k \geq 6$ we have $\us_{\PP \D_k}: \D \times \P T\P^2 \lrab \UL_{\PP \D_k}$ given by
\bgd
\{\us_{\PP \D_k}(\ff,l_p)\}(f^{\otimes (1+ \epsilon_k)} \otimes p^{\otimes d(1+\epsilon_k)} \otimes v^{\otimes (k-2+\epsilon_k)}\otimes w^{\otimes (2\epsilon_k)}):= f_{12}^{\epsilon_k} \D^f_k,
\edd
where, as seen previously, $\epsilon_k$ is the order of the pole of $f_{12} =0$ for the section $\D^f_k$. In particular, $\epsilon_6 = 0$, $\epsilon_7 =1$ and 
$\epsilon_8=3$. The expressions for $\A^f_k$ (resp. $\D^f_k$) are given below explicitly in \eqref{Formula_Ak} (resp. \eqref{Formula_Dk}), till $k=7$ (resp. till $k=8$). 
The algorithm to obtain $\A^f_k$ (resp. $\D^f_k$) for any $k$ is given in 
Lemma \ref{fstr_prp} (resp. Lemma \ref{fstr_prp_Dk}). \\ 
\hf\hf Here is an explicit  formula for  
$\A^f_k$ till $k=7$.  
\begin{align}
\label{Formula_Ak}
\A^f_3&= f_{30},\qquad
\A^f_4 = f_{40}-\frac{3 f_{21}^2}{f_{02}}, \qquad
\A^f_5= f_{50} -\frac{10 f_{21} f_{31}}{f_{02}} + 
\frac{15 f_{12} f_{21}^2}{f_{02}^2} \nonumber \\
\A^f_6 &= f_{60}- \f{ 15 f_{21} f_{41}}{f_{02}}-\f{10 f_{31}^2}{f_{02}} + \f{60 f_{12} f_{21} f_{31}}{f_{02}^2}
   +
   \f{45 f_{21}^2 f_{22}}{f_{02}^2} - \f{15 f_{03} f_{21}^3}{f_{02}^3}
   -\f{90 f_{12}^2 f_{21}^2}{f_{02}^3} \nonumber \\  
\A^f_7 &= f_{70} - \frac{21 f_{21} f_{51}}{f_{02}} 
- \frac{35 f_{31} f_{41}}{f_{02}} + \frac{105 f_{12} f_{21} f_{41}}{f_{02}^2} + \f{105 f_{21}^2 f_{32}}{f_{02}^2} + 
\f{70 f_{12} f_{31}^2}{f_{02}^2}+ \f{210 f_{21}f_{22}f_{31}}{f_{02}^2} \nonumber \\
&
-\f{105 f_{03} f_{21}^2 f_{31}}{ f_{02}^3}
-\f{420 f_{12}^2 f_{21} f_{31}}{f_{02}^3}
-\f{630 f_{12}f_{21}^2 f_{22}}{f_{02}^3}
-\f{105 f_{13} f_{21}^3}{f_{02}^3}
+ \f{315 f_{03} f_{12} f_{21}^3}{f_{02}^4}
+ \f{630 f_{12}^3 f_{21}^2}{f_{02}^4}.
\end{align} 
Here is an explicit formula for 
$\D^f_k$  till $k=8$. 
\begin{align}
\label{Formula_Dk}
\D^f_6 &=  f_{40},\qquad 
\D^f_7 =   f_{50} -\f{5 f_{31}^2}{3 f_{12}},\qquad
\D^f_8 = f_{60} + 
\frac{5 f_{03} f_{31} f_{50}}{3 f_{12}^2} 
-\frac{5 f_{31} f_{41}}{f_{12}} - \frac{10 f_{03} f_{31}^3}{3 f_{12}^3} 
+ \frac{5 f_{22} f_{31}^2}{f_{12}^2}.
\end{align}  

\subsection{The spaces involved.}

\hf\hf We begin by explaining a terminology. If $l_{\p} \in \P T_{\p}\P^2$, then we say that $v \in l_{\p}$ if 
$v$ is a tangent vector in $T_{\p}\P^2$ and lies over the fibre of $l_{\p}$.
We now define the spaces that we will encounter.
\label{definition_of_px}
\begin{align*}
\X_k &:= \{ ( \ff,\p) \in \D \times \P^2~~~~~: \textnormal{$f$ has a singularity of type $\X_k$ at $\p$} \} \\
\hat{\X}_k &:= \{ (\ff,\lp) \in \D \times \P T\P^2: \textnormal{$f$ has a singularity of type $\X_k$ at $\p$} \} ~= \pi^{-1}(\X_k) \\
\textnormal{if $~k>1$ } \quad 
\PP \A_k &:= \{ (\ff,\lp) \in \D \times \P T\P^2: 
\textnormal{$f$ has a singularity of type $\A_k$ at $\p$},\\
& \qquad\qquad\qquad\qquad \qquad  \qquad \nabla^2 f|_p(v, \cdot) =0\,\,\textup{if}\,\,v \in  l_{\p}\}\\ 
\PP \D_4 &:= \{ (\ff,\lp) \in \D \times \P T\P^2: 
\textnormal{$f$ has a singularity of 
type $\D_4$ at $\p$}, \\ 
& \qquad \qquad \qquad \qquad \qquad  \qquad \nabla^3 f|_p(v,v,v) =0\,\,\textup{if}\,\,v  \in  l_{\p}\} \\
\textnormal{if $~k>4$} \qquad 
\PP \D_k &:= \{ (\ff,\lp) \in \D \times \P T\P^2: \textnormal{$f$ has a singularity of 
type $\D_k$ at $\p$} \\
& \qquad \qquad \qquad \qquad \qquad\qquad \nabla^3 f|_p(v, v, \cdot) =0\,\,\textup{if}\,\,v \in l_{\p}\} \\
\textnormal{if $k=6, 7$ or  $8$}\qquad 
\PP \E_k &:= \{ (\ff,\lp) \in \D \times \P T\P^2: \textnormal{$f$ has a singularity of 
type $\E_k$ at $\p$}\\
& \qquad \qquad \qquad \qquad \qquad\qquad \nabla^3 f|_p(v, v, \cdot) =0\,\,\textup{if}\,\,v\in l_{\p}\} \\ 
\textnormal{if $~k>4$} \qquad 
\PP\D_k^{\vee} &:= \{ (\ff,\lp) \in \D \times \P T\P^2: 
\textnormal{$f$ has a singularity of 
type $\D_k$ at $\p$}, \\
& \qquad \qquad \qquad \qquad \qquad  \qquad \nabla^3 f|_p(v,v,v) =0, ~~\nabla^3 f|_p(v,v,w) \neq 0, \\ 
& \qquad \qquad \qquad \qquad \qquad  \qquad  \,\,\textup{if}\,\,v  \in  l_{\p} ~~\textup{and} ~~w \in \frac{T_{\p}\P^2}{l_{\p}}\} 
\end{align*}
\ni We also need the definitions for a few other spaces which will make our computations convenient. 
\begin{align*}
\hat{\A}_1^{\#} := \{ (\ff, l_p) \in \D \times \P T\P^2 &: 
 f(p) =0, \nabla f|_p =0,  \nabla^2 f|_p(v, \cdot) \neq 0, \forall ~v\neq 0 \in l_{\p}  \}  \\
\hat{\D}_4^{\#} := \{ (\ff, \lp) \in \D \times \P T\P^2 &: 
f(p) =0, \nabla f|_p =0, \nabla^2 f|_p \equiv 0, \nabla^3 f|_p (v,v,v) \neq 0, \forall ~v \neq 0 \in l_{\p} \} \\
\hat{\XC}_8^{\#} := \{ (\ff, \lp) \in \D \times \P T\P^2 &: 
f(p) =0, \nabla f|_p =0, \nabla^2 f|_p \equiv 0, \nabla^3 f|_p =0,  \\ 
                & \qquad \nabla^4 f|_p (v,v,v, v) \neq 0 ~\forall ~v \neq 0  \in l_{\p} \}    \\
\XX_{8}^{\# \flat} :=  \{ (\ff, \lp) \in \D \times \P T\P^2 &: (\ff, \lp) \in \hat{\XC}_8^{\#}, \us_{\J}(\ff, \lp)\neq 0, 
\textnormal{where $\us_{\J}$ is defined in \eqref{J_psi} }\}.                        
\end{align*}


\section{Local structure of holomorphic sections}
\label{condition_for_sing}

\hf\hf We give a necessary and sufficient criterion for a curve $f^{-1}(0)$ to have a singularity 
of type $\X_k$ at the point $\p$. Let $\q=\q(x,y)$ be a holomorphic function defined on a neighborhood
of the origin in $\C^2$ and $i,j$ be non-negative integers. We define
\bgd
\q_{ij}:=\frac{\partial^{i+j} \q}{\partial^i x\partial^j y}\bigg|_{(x,y)=(0,0)}\,.
\edd
\begin{lmm}
\label{ift}
Let $\q =\q(x,y)$ be a holomorphic function defined on a neighborhood
of the origin in $\C^2$ such that $\q_{00} =0$ and $\nabla \q |_{(0,0)} \neq 0$. Then there exists a coordinate 
chart $(u,v)$ centered at the origin so that $\q(u,v)= v^2 + u.$ 
\end{lmm}
\pf Follows immediately by considering the Taylor expansion of $\q$. \qed   
\begin{cor}
\label{A0_node_condition_cor}
A curve $\q^{-1}(0)$ has an $\A_0$-node at the origin if and only if it satisfies the hypothesis of Lemma \ref{ift}. 
\end{cor}
\begin{lmm}
\label{ml}
Let $\q = \q(x,y)$ be a holomorphic function defined on a neighbourhood
of the origin in $\C$ such that $\q_{00}, \nabla \q |_{(0,0)} =0 $ and
$\nabla^2 \q|_{(0,0)}$ is non-degenerate.  Then there exists a coordinate 
chart $(u,v)$ centered at the origin so that $\q(u,v) = v^2 + u^2.$
\end{lmm}
\pf This is the Morse Lemma, which again follows by considering the Taylor expansion of $\q$.   \qed 
\begin{cor}
\label{A1_node_condition_cor}
A curve $\q^{-1}(0)$ has an $\A_1$-node if and only if it satisfies the hypothesis of Lemma \ref{ml}. 
\end{cor}
\begin{lmm}\label{fstr_prp}
Let $\q =\q(\rq, \sq)$ be a holomorphic function defined on a neighbourhood of the origin in $\C$ such that 
$\q(0,0), ~\nabla \q|_{(0,0)}=0$ and there exists a non-zero vector $w=(w_1,w_2)$ such that at the origin $ \nabla^2 f (w, \cdot)=0$, i.e., the Hessian is degenerate.
Let $x = w_1 \rq + w_2 \sq, y = -\ov{w}_2 \rq + \ov{w}_1 \sq $ and $\q_{ij}$ be the partial derivatives with respect to the new variables $x$ and $y$. 
If $\q_{02} \neq 0$, there exists a coordinate chart $(u, v)$ centered around the 
origin in $\C^2$ such that 
\bge
\label{ak1}
\q = \left\{\begin{array}{rl}
v^2, & \textup{or}\\
v^2 + u^{k+1}, & \textup{for some $k\geq 2$.}
\end{array}\right.
\ede
\end{lmm}
\begin{rem}
In terms of the new coordinates we have $\q_{00}= \q_{10}= \q_{01}= \q_{20}= \q_{11} =0$ and $\q_{02} \neq 0.$ 
Here $\partial_x + 0 \partial_y = (1,0)$ is the distinguished direction along which the Hessian is degenerate. 
\end{rem}  
\pf Let the Taylor expansion of $\q$ in the new coordinates be given by 
\bgd
\q(x,y) = \Z_0(x) + \Z_1(x)y + \Z_2(x) y^2 + \ldots.
\edd
By our assumption on $\q$, $\Z_2(0) \neq 0.$ We claim that there exists a holomorphic function $\Y(x)$ such that 
after we make a change of coordinates $y = y_1 + \Y(x)$, the function $\q$ is given by 
\bgd
\q = \hat{\Z}_0(x) + \hat{\Z}_2(x) y_1^2 + \hat{\Z}_3(x) y_1^3 + \ldots 
\edd
for some $\hat \Z_k(x)$ (i.e., $\hat{\Z}_1(x) \equiv 0$). To see this, we note that this is possible if $\Y(x)$ satisfies the identity
\begin{align}
\Z_1(x) + 2 \Z_2(x) \Y + 3 \Z_3(x) \Y^2 + \ldots \equiv 0.  \label{psconvgg}
\end{align}
Since ~$\Z_2(0) \neq 0$, $\Y(x)$ exists by the Implicit Function Theorem
\footnote{Moreover it is unique if we require $\Y(0) =0$.}.
Therefore, we can compute $\Y(x)$ as a power series using \eqref{psconvgg} and then
compute $\hat \Z_{0}(x)$. Hence, 
\begin{align}
\label{Ak_node_conditionn}
\q &= v^2 + \frac{\A^{\q}_3}{3!}x^3 + 
\frac{\A^{\q}_4}{4!} x^4 + \ldots,  
~~\textnormal{where} ~~~ v = \sqrt{(\hat{\Z}_2 + \hat{\Z}_3 y_1 + \ldots)} y_1, 
\end{align}
satisfies \eqref{ak1}. \qed \\
\hf\hf Following the above procedure we find $A_i^\rho$ for $i=3,\ldots,7$. In particular,  
\begin{align}
\label{Ak_sections}
\A^{\q}_3&= \q_{30}\,,\qquad
\A^{\q}_4 = \q_{40}-\frac{3 \q_{21}^2}{\q_{02}}\,, \qquad
\A^{\q}_5= \q_{50} -\frac{10 \q_{21} \q_{31}}{\q_{02}} + 
\frac{15 \q_{12} \q_{21}^2}{\q_{02}^2} 
%
\end{align} 
\begin{cor}
\label{Ak_node_condition_cor}
Let the hypothesis be as in Lemma \ref{fstr_prp}. The curve $\q^{-1}(0)$ has an $\A_k$-node (for $k \geq 2$) at the origin if and only if $\q_{02} \neq 0$ and the directional derivatives 
$\A^{\q}_i$ obtained in \eqref{Ak_node_conditionn} are zero for all $i \leq k$ and $\A^{\q}_{k+1} \neq 0$. Furthermore, if $\tt$ is any holomorphic function 
that does not vanish at the origin,
then 
\begin{align}
\label{Ak_trivialization}
\A^{\tt \q}_{k+1} = \tt_{00} \A^{\q}_{k+1} \qquad \textnormal{and} \qquad (\tt \q)_{02}^{k-3}\A^{\tt \q}_{k+1} = \tt_{00}^{k-2} \q_{02}^{k-3} \A^{\q}_{k+1}. 
\end{align}
Finally, if $\A^{\q}_{i}=0$ for $i\leq k$ then the quantity $\A^{\q}_{k+1}$ is invariant under 
\begin{align}
x \lra x+ \TT_1(x,y),\,\,y\lra y + \TT_2(x,y) \label{Ak_independent_curve} \\ 
y \lra y + x, \,\,x\lra x  \label{Ak_wd_quotient}
\end{align}
where $\TT_1$ and $\TT_2$ are holomorphic functions that vanish at the origin and 
whose derivative also vanish at the origin, i.e.,  
\bgd
\TT_i(0,0) =0, ~~\nabla \TT_i(0,0) =0, \qquad i=1,2. 
\edd 
\end{cor}
\pf The first assertion follows immediately from \eqref{Ak_node_conditionn}. 
To prove \eqref{Ak_trivialization}, note that by \eqref{Ak_node_conditionn} 
\begin{align*}
\A^{\q}_{k+1} &= \frac{\partial^{k+1} \q(x,v)}{\partial x^{k+1}}\bigg|_{(0,0)}  \implies \A^{\tt \q}_{k+1} = \frac{\partial^{k+1}\tt \q(x,v)}{\partial x^{k+1}}\bigg|_{(0,0)} = \tt_{00} \A^{\q}_{k+1} 
\end{align*}
which follows from the fact that $\A^{\q}_i =0$ for all $i\leq k$. The second 
equation follows similarly by observing that $ (\tt \q)_{02} = \tt_{00} \q_{02}.$ 
We have omitted here the proofs of \eqref{Ak_independent_curve} and \eqref{Ak_wd_quotient}.  
The details of the proof can be found in \cite{BM_Detail}. \qed 
\begin{rem}
\label{induced_section_Ak}
\ni The quantity $\q_{02}^{k-3} \A^{\q}_{k}$ is defined even when $\q_{02} =0$. 
These quantities induce sections $\us_{\PP \A_k}$ of the line bundles 
$\UL_{\PP \A_k} \lra \D \times \P T\P^2 $ of subsection \ref{summary_vector_bundle_definitions}. 
The induced section is \textit{defined to be} 
\begin{align}
\label{Ak_defining_sections_finally}
\{\us_{\PP \A_k}(\ff, \lp)\}( f^{\otimes (k-2)} \otimes p^{\otimes d} \otimes v^{\otimes k} \otimes w^{\otimes (2k-6)}) := f_{02}^{k-3} \A^{f}_{k},
\end{align}
where $\A^{f}_k$ is the number we get by replacing $\q_{ij}$ by $f_{ij}$ in $\A^{\q}_k$ 
($f_{ij}$ is defined in \eqref{abbreviation}). Note that \eqref{Ak_trivialization} and \eqref{Ak_independent_curve} imply that this section is well 
defined restricted to $\us_{\PP \A_{k-1}}^{-1}(0)$, i.e., it is independent of the trivialization and independent of the curve chosen. 
This is easily seen by unwinding definition \ref{vertical_derivative_defn}. 
The details of this can be found in \cite{BM_Detail}.
Finally, note that \eqref{Ak_wd_quotient} implies that the section $\us_{\PP \A_k}$ is well defined 
on the quotient space $\frac{\pi^*(T \P^2)}{\G^*}$, since the quantity in \eqref{Ak_defining_sections_finally} is invariant under $w \lra w+v$.
\end{rem}
\hf\hf Next we analyze singularities when the Hessian is identically zero. 
\begin{lmm}\label{fstr_prp_D4}
Let $\q = \q(x,y)$ be a holomorphic function defined on a neighbourhood
of the origin in $\C$ such that $\q_{00}, \nabla \q|_{(0,0)}, \nabla^2 \q|_{(0,0)}=0$ and 
there does not exist a non-zero vector $w=(w_1,w_2)$ such that 
at the origin $\nabla^3 \q (w,w,\cdot) =0$.
Then, there exists   
a coordinate chart $(u,v)$ centered at the origin so that $\q(u,v)=u^3+ v^3.$
\end{lmm}
\pf The Taylor expansion of $\q$ is given by 
\begin{align*}
\q &= \frac{\q_{30}}{6}x^3 + \frac{\q_{21}}{2} x^2 y 
+ \frac{\q_{12}}{2} x y^2 + \frac{\q_{03}}{6} y^3 + \ldots 
\end{align*} 
We make an observation from multi-linear algebra 
that if $\nabla ^3 \q$ is non-degenerate, then the 
cubic term in the Taylor expansion has no 
repeated factors (this analogous to the 
similar fact that if $\nabla^2 \q$ is non-degenerate, 
then the quadratic term in the Taylor expansion 
is not a perfect square). Hence, we can make 
a linear change of coordinates so that $\q$ is given by  
\begin{align*}
\q &= x_1^3 + y_1^3 + \eta x_1^2 y_1^2 + x_1^3 h_1 + 
y_1^3 h_2
\end{align*}
where $h_1$ and $h_2$ are holomorphic functions 
of $x_1$ and $y_1$ vanishing at the origin
and $\eta$ is some number. 
We now make a change of coordinate $ x_1 = x_2 + A y_1^2$ to get rid of the coefficient of $x_2^2 y_1^2$. 
Equating coefficients we get $ 3 A + \eta =0.$
Hence,
\bgd
\q = x_2^3 + y_1^3 + x_2^3 h_3 + y_1^3 h_4,
\edd
where $h_3$ and $h_4$ are holomorphic functions of 
$x_2$ and $y_1$ that vanish at the origin. This is equivalent to $u^3+v^3$ after a change of coordinates. \qed 
\begin{cor}
\label{D4_node_condition_cor}
A curve $\q^{-1}(0)$ has a $\D_4$-node if and only if it satisfies the hypothesis of Lemma \ref{fstr_prp_D4}. 
\end{cor}
\begin{lmm}\label{fstr_prp_Dk}
Let $\q=\q(\rq,\sq)$ be a holomorphic function defined on a neighbourhood 
of the origin in $\C$ such that $\q_{00}, \nabla \q|_{(0,0)}, \nabla^2 \q|_{(0,0)}=0$ and 
there exists a non-zero vector $w=(w_1,w_2)$ such that at the origin $\nabla^3 \q(w,w,\cdot) =0$.
Let $x = w_1 \rq + w_2 \sq, ~y = -\ov{w}_2 \rq + \ov{w}_1 \sq$ and $\q_{ij}$ be the partial derivatives with respect to 
the new variables $x$ and $y$. If $\q_{12} \neq 0$, there exists a coordinate chart $(u, v)$
centered around the origin in $\C^2$ such that 
\bgd
\label{Dk_form}
\q(u,v) \equiv \left\{\begin{array}{rl}
v^2 u & \textup{or}\\
v^2 u + u^{k-1} & \textup{for some $k\geq 5$}.
\end{array}\right.
\edd
\end{lmm}
Note that in terms of the new coordinates we have  
\bgd
\q_{00} = \q_{10}= \q_{01}= \q_{20}= \q_{11}= \q_{02}= \q_{30}= \q_{21} =0, \q_{12} \neq 0.
\edd 
\pf In terms of the new coordinates $x$ and $y$, The Taylor expansion of $\q$ is given by 
\begin{align*}
\q &= \frac{\q_{12}}{2} x y^2 + \frac{\q_{03}}{6} y^3 
+ \frac{\q_{40}}{24} x^4 + \ldots 
\end{align*}
We claim that there exists a holomorphic function $\GG(y)$, such that after making a change of coordinate
$x = x_1 + \GG(y)$, the function is given by $\q = x_1 \g (x_1, y)$, i.e, we can kill off all powers of $y$. Assuming such a $\GG(y)$ exists, we can now apply the argument as in Lemma \ref{fstr_prp}. 
Let 
\bgd
\g(x_1, y) = \Z_0(x_1) + \Z_1(x_1) y + \Z_2(x_1) y^2 + \ldots.  
\edd
We can make a change of coordinate $y = y_1 + \Y(x_1)$ so that 
\bgd
\g = \hat{\Z}_0(x_1) + \Z_2(x_1) y_1^2 + \ldots
\edd
i.e., $\hat{\Z}_1(x_1) \equiv 0 $. This is possible since $\Z_2(0) \neq 0$\footnote{Notice that $\Z_2(0) = \frac{\q_{12}}{2}$ is non-zero by hypothesis.}. That gives us 
\begin{align*}
\q &= x_1 ( \hat{\Z}_0 (x_1) + \underbrace{\hat{\Z}_2(x_1) y_1^2 + \hat{\Z}_3(x_1) y_1^3 + \ldots  }_{y_2^2} )
\end{align*} 
If we set 
\bge
\hat{\Z}_0(x_1) = \frac{\D^{\q}_6}{4 !} x_1^3 + \frac{\D^{\q}_7}{5 !} x_1^4 + \ldots \label{Dkdefn}
\ede
then $\q$ is given by 
\begin{align}
\q &= x_1(y_2^2 + \frac{\D^{\q}_6}{4 !} x_1^3 + \frac{\D^{\q}_7}{5 !} x_1^4 + \ldots) \nonumber 
\end{align}
If $\hat{\Z}_0(x_1)\equiv 0$ then $\q=x_1y_2^2$ is of the intended form. Otherwise let $k$ be the smallest integer such that $\D^{\q}_{k+1} \neq 0$. Let 
\begin{align}
\label{Dk_node_conditionn}
x_2  &= \sqrt[k-1]{ \frac{\D^{\q}_{k+1}}{(k-1)!} x_1^{k-1} + \frac{\D^{\q}_{k+2}}{k!} x_1^{k}+\ldots}
\end{align}
with $x_{1} = \kap x_{2} + O(x_2^2)$ and $\kap = \big((k-1)!/ \D^{\q}_{k+1}\big)^{k-1}$. In these new coordinates, $\q$ is given by 
\begin{align*}
\q&= (\kap x_2 + x_2^2 h) y_2^2 + x_2^{k-1}
\end{align*}
for some holomorphic function $h(x_2,y_2)$. Now define $y_3 = y_2\sqrt{\kap  + x_2 h }$. Therefore,
\begin{align}
\label{Dk_obvious}
\q = y_3^2 x_2 + x_2^{k-1}
\end{align}
to get the intended form as in \eqref{Dk_form}. \\
\hf It remains to show that $\GG(y)$ exists. By our assumption on $\q$, we know that  $\q_{30}=\q_{21}=0$ and $\q_{12}\neq0$. 
Hence, the Taylor expansion of $\q$ is given by
\begin{align}
\q(x,y)&=\T_{12} (x,y) xy^2+\T_{03}(y) y^3+\T_{40} (x,y) x^4+ \T_{31}(x,y) x^3 y \label{Dk_existence_of_function}
\end{align}
for some holomorphic functions $\T_{ij}$ with $\T_{12}(0,0)\neq0$. Recall that we want $x = x_1 + \GG(y)$ so that $\q = x_1 g(x_1, y) $, i.e., the coefficients of $y^n$ are killed for all $n$. 
This is equivalent to saying that we want to find a $\GG$ such that $\q(\GG(y), y) =0$\footnote{Note that $\q(\GG(y), y)$ gives us the part of the  
Taylor expansion of $\q(x_1 + \GG(y), y)$ that only depends on $y$. We desire that part to be zero, which is equivalent to killing off all the $y^n$ terms in the expansion of $\q(x_1 + \GG(y), y)$.}. Plugging in $x = x_1 + y \HH(y)$ in \eqref{Dk_existence_of_function} we get that
\bgd
0=\q(y \HH(y), y) =\T_{12}(y \HH(y),y) y^3 \HH(y) +  \T_{03}(y) y^3 + \T_{40}(y \HH(y),y) y^4 \HH(y)^4 + \T_{31}(y \HH(y), y) y^3 \HH(y).
\edd
This implies that
\bgd
\T_{12}(y \HH(y), y) \HH(y) +  \T_{03}(y)  + \T_{40}(y \HH(y),y) y \HH(y)^4 + \T_{31}(y \HH(y), y) \HH(y) =0 
\edd
By the implicit function theorem, $\HH(y)$ exists since $\T_{12}(0,0) \neq 0$, whence $\GG(y)= y \HH(y)$ exists. \qed\\

\hf\hf In practice we first find $\HH(y)$ as a power series and then find $\Z(x_1)$ as a power series. That ultimately gives us $\D^{\q}_i$. Following the above procedure, we prove \eqref{Formula_Dk}.

\begin{cor}
\label{Dk_node_condition_cor}
Let the hypothesis be as in Lemma \ref{fstr_prp_Dk}. Then the curve $\q^{-1}(0)$ has a $\D_k$-node if and only if the directional derivatives 
$\D^{\q}_i$ obtained in \eqref{Dkdefn} are zero for all $i \leq k$ and $\D^{\q}_{k+1} \neq 0$. 
Furthermore, if $\tt$ is a holomorphic function that does not vanish at the origin, then 
\begin{align}
\label{Dk_trivialization}
\D^{\tt \q}_{k+1} = \tt_{00} \D^{\q}_{k+1} \qquad \textnormal{and} \qquad (\tt \q)_{12}^{k-6}\D^{\tt \q}_{k+1} = \tt_{00}^{k-5} \q_{12}^{k-6} \D^{\q}_{k+1}. 
\end{align}
Finally, if $\D_i^{\q}=0$ for $i\leq k$ then the quantity $\D^{\q}_{k+1}$ is invariant under 
\begin{align}
x \lra x+ \TT_1(x,y),\,y\lra y + \TT_2(x,y)\label{Dk_independent_curve} \\ 
x\lra x,\,y \lra y + x,\label{Dk_wd_quotient}
\end{align}
where $\TT_1$ and $\TT_2$ are holomorphic functions that vanish at the origin and 
whose derivative also vanish at the origin, i.e.,  
\bgd
\TT_i(0,0) =0, ~~\nabla \TT_i(0,0) =0, \qquad \textnormal{where} ~~i=1,2. 
\edd 
\end{cor}
\pf The first assertion follows immediately from \eqref{Dk_obvious}. 
To prove the second assertion, note that by \eqref{Dkdefn} 
\begin{align*}
\D^{\q}_{k+1} &= \frac{\partial^{k+1} \q(x_1,y_2)}{\partial x_1^{k+1}}\bigg|_{(0,0)}  \implies \D^{\tt \q}_{k+1} = \frac{\partial^{k+1}\tt \q(x_1,y_2)}{\partial x_1^{k+1}}\bigg|_{(0,0)} = \tt_{00} \D^{\q}_{k+1} 
\end{align*}
which follows from the product rule and the fact that $\D^{\q}_i =0$ for all $i\leq k$. The second 
equation follows similarly by observing that  $(\tt \q)_{12} = \tt_{00} \q_{12}$. 
We have omitted here the proofs of \eqref{Dk_independent_curve} and \eqref{Dk_wd_quotient}.  
The details of the proof can be found in \cite{BM_Detail}. \qed

\begin{rem}
\label{induced_section_Dk}
Similar to remark \ref{induced_section_Ak}, the quantities 
$\q_{12}^{\epsilon_k} \D^{\q}_{k}$    
induce
sections of the 
line bundle  $ \UL_{\PP \D_k} \lra \D\times \P T\P^2$ given by 
\begin{align}
\label{Dk_defining_sections_finally}
\{\us_{\PP \D_k}(\ff, \lp)\}( f^{\otimes (1+\epsilon_k)} \otimes p^{\otimes d(1+\epsilon_k)} \otimes v^{\otimes (k-2+\epsilon_k)} \otimes w^{\otimes (2\epsilon_k)}) := f_{12}^{\epsilon_k} \D^{f}_{k},
\end{align}
where $\epsilon_k$ is the order of the pole of $\D^{\q}_k$ at $\q_{12} =0$. 
Equations \eqref{Dk_trivialization} and \eqref{Dk_independent_curve} imply that   
this section is well defined restricted to $\us_{\PP \D_{k-1}}^{-1}(0)$. Equation \eqref{Dk_wd_quotient} implies that the section $\us_{\PP \D_k}$ is well defined 
on the quotient space.   
\end{rem}

\hf\hf Next we analyze the singularities $\E_6$ and $\E_7$. 
\begin{lmm}\label{fstr_prp_E6}
Let $\q = \q(\rq,\sq)$ be a holomorphic function defined on a neighbourhood
of the origin in $\C$ such that $\q_{00}= \nabla \q|_{(0,0)}= \nabla^2 \q|_{(0,0)}=0$ and 
there exists a non-zero vector $w=(w_1,w_2)$ such that at the origin $\nabla^3 \q(w,w,\cdot) =0$. 
Let $x = w_1 \rq + w_2 \sq,  ~y = -w_2 \rq + w_1 \sq$ and $\q_{ij}$ be partial derivatives with respect 
to the new coordinates, $x$ and $y$. If $\q_{12} =0$ and  $\q_{03}\neq 0, \q_{40} \neq 0$, there exists   
a coordinate chart $(u,v)$ centered at the origin so that
\begin{equation}
\label{E6eqn_e}
\q(u,v)= v^3 + u^4.
\end{equation}
\end{lmm}
Note that in terms of the new coordinates $x$ and $y$, we get that 
\bgd
\q_{00} =\q_{10}= \q_{01}=\q_{20}= \q_{11}= \q_{02}= \q_{30}= \q_{21}= \q_{12} =0, \q_{03}\neq 0, \q_{40} \neq 0.
\edd 
\pf  The Taylor expansion of $\q$ is given by,
\begin{equation}\label{E6fexp_e}
\q(x,y)=\T_{03}(x,y)y^3+\T_{40}(x,y)x^4+\kq_1 x^3y+\kq_2 x^2y^2+\kq_3 x^3y^2
\end{equation}
in terms of the new variables $x$ and $y$, 
for some constants $\kq_1,\kq_2,\kq_3$ and   
holomorphic functions $\T_{03}$ and $\T_{40}$ on a neighborhood of the origin in $\C^2$
such that $\T_{03}(0,0),\T_{40}(0,0)\neq 0$.
We claim that there exists constants $\etq_1$, $\etq_2$ and $\etq_3$, such that if we make the substitution  
$$ x = \xq + \etq_1 \yq, \qquad y = \yq + \etq_2 \xq^2 + \etq_3 \xq^3$$ 
then $\q$ is given by 
\begin{align}
\label{E6_intermediate}
\q &= \hat{\T}_{03}(\xq, \yq) \yq^{3} + \hat{\T}_{40}(\xq, \yq) \xq^4  
\end{align}
such that $\hat{\T}_{03}(0,0), ~\hat{\T}_{40}(0,0) \neq 0$. To see this, we note that this is possible if the following equations are satisfied: 
\begin{align*}
\frac{\q_{40}}{6} \etq_1+\kq_1=0, \qquad
\frac{\q_{03}}{2} \etq_2 + \frac{\q_{40}}{4} \etq_1^2+3\kq_1 \etq_1+\kq_2=0,\\
\frac{\q_{03}}{2} \etq_3 +\frac{\q_{13}}{2} \etq_2 + \frac{\q_{50}}{12} \etq_1^2+ \frac{\q_{41}}{6} \etq_1 +  3 \kq_1 \etq_1^2 \etq_2+4 \kq_2 \etq_1 \etq_2 +\kq_3=0.
\end{align*}
Solutions to $\eta_1,\eta_2,\eta_3$ exist, since $\q_{40} \neq 0$ and $\q_{03}\neq 0$. It is easy to see that \eqref{E6_intermediate} is equivalent 
to \eqref{E6eqn_e} after a change of coordinates, since $\hat{\T}_{03}(0,0), ~\hat{\T}_{40}(0,0) \neq 0$. \qed 
\begin{cor}
\label{E6_node_condition_cor}
A curve $\q^{-1}(0)$ has an $\E_6$-node if and only if 
it satisfies the hypothesis of Lemma \ref{fstr_prp_E6}. 
\end{cor}
\begin{lmm}\label{fstr_prp_E7}
Let $\q =\q(\rq,\sq)$ be a holomorphic function defined on a neighbourhood of the origin in $\C$ such that 
$\q_{00}, \nabla \q|_{(0,0)}, \nabla^2 \q|_{(0,0)}=0$ and there exists a non-zero vector $w=(w_1,w_2)$ such that 
at the origin $\nabla^3 \q(w,w,\cdot) =0$. Let $x = w_1 \rq + w_2 \sq, ~y = -w_2 \rq + w_1 \sq.$ Let $\q_{ij}$ be the partial derivatives with respect 
to the new variables $x$ and $y$. If $\q_{12}= \q_{40}=0$ and $\q_{03}\neq 0, \q_{31}\neq 0$, then there exists   
a coordinate chart $(u,v)$ centered at the origin so that
\begin{equation}
\label{E7eqn_e}
\q(u,v)= v^3 + u^3 v.
\end{equation}
\end{lmm}
In terms of the new coordinates $x$ and $y$, we get that 
\bgd
\q_{00}=\q_{10}= \q_{01}=\q_{20}= \q_{11}= \q_{02} =\q_{30}= \q_{21}= \q_{12}= \q_{40} =0, \q_{03}\neq 0, \q_{31} \neq 0.
\edd
\pf The Taylor expansion of $\q$ is given by,
\begin{equation}
\label{E7fexp_e}
\q(x,y)=\T_{03}(x,y)y^3+\T_{31}(x,y)x^3y+\kq x^2y^2+\T_{50}(x)x^5 
\end{equation}
in terms of the  new coordinates $x$ and $y$, for constant $\kq$ and  holomorphic functions $\T_{03}$, $\T_{31}$, and 
$\T_{50}$ such that  
$\T_{03}(0,0),\T_{31}(0,0)\neq 0$.
We claim that there exists a holomorphic function 
$\Bq(\xq)$ and   
constant $\etq_1$  such that if we make the 
substitution 
$$x = \xq + \etq_1 \yq, \qquad y = \yq + \Bq(\xq) \xq^2 $$ 
then $\q$ is given by 
\begin{align}
\label{E7_intermediate}
\q &= \hat{\T}_{03}(\xq, \yq) \yq^3 + \hat{\T}_{31}(\xq, \yq) \xq^3 \yq  
\end{align}
To see this, note that this is possible if the following equations are satisfied: 
\begin{align*}
\T_{31}(\xq, \Bq(\xq)\xq^2)\Bq(\xq)+
\T_{03}(\xq, \Bq(\xq) \xq^2)\xq \Bq(\xq)^3+k\xq \Bq(\xq)^2+
\T_{50}(\xq)& =0, \\
\frac{\q_{31}}{2} \etq_1 + \frac{\q_{03}}{2} \Bq(0)+\kq & =0.
\end{align*}
A solution to $\Bq(\xq)$ exists since $\q_{31} \neq 0$. We see that
\bgd
\Bq(0) = -\frac{\T_{50}(0)}{\T_{31}(0,0)}  = -\frac{\q_{50}}{20 \q_{31}}
\edd  
which implies the existence of $\eta_1$. It is easy to see that \eqref{E7_intermediate} is equivalent to \eqref{E7eqn_e} after a change of coordinates, since $\hat{\T}_{03}(0,0), ~\hat{\T}_{31}(0,0) \neq 0$.  \qed 
\begin{cor}
\label{E7_node_condition_cor}
A curve $\q^{-1}(0)$ has an $\E_7$-node if and only if it satisfies the hypothesis of Lemma \ref{fstr_prp_E7}. 
\end{cor}
\hf\hf Let us now summarize certain facts about  sections of vector bundles, involving the vertical derivative.   
\begin{lmm}
\label{wdc}
Let $L\lra M$ a complex line bundle over a two dimensional complex manifold $M$, $\s: M\lra L$ a holomorphic section and $q\in M$ a point in $M$. Let $v, w \in T_q M$ be two tangent vectors at the point $q$. Then the following are true: 
\begin{enumerate}
\item \label{wd_nabla_all}  If $\s(q)=0$ and $\nabla^i \s|_{q} =0$ for all $ i<k$ then $\nabla^{k} \s|_{q}$ is well defined. Furthermore, for any tangent vectors $v, w \in T_{q}M$ and non-negative integers $i$ and $j$ such that  $i+j =k$, the quantity 
\begin{align}
\label{pedantic_point}
\ho{\s}_{ij}&:= \nabla^{i+j} \s|_{q} (\underbrace{v,\cdots v}_{\textnormal{$i$ times}}, \underbrace{w,\cdots w}_{\textnormal{$j$ times}}) 
\end{align}
is also well defined. 
\item \label{wd_pa3} If $\s(q)=0$, $\nabla \s|_{q} =0$ and $\nabla^2 \s|_{q}(v, \cdot) =0$
then  $ \nabla^3 \s|_{q}(v,v,v)$ is also well defined.
\item \label{wd_pak} Let $\A^{\ho{\s}}_k$  be the corresponding sections induced from the quantities 
$\A^{\q}_k$ obtained in \eqref{Ak_node_conditionn} 
by replacing $\q_{ij}$ with $\ho{\s}_{ij}$.   
If $\s(q)=0$, $\nabla \s|_{q} =0$, $\nabla^2 \s|_{x}(v, \cdot) =0$
and $\ho{\s}_{02}^{i-3} \A^{\ho{\s}}_i =0$ for all $i <k$,
then  $\ho{\s}_{02}^{k-3} \A^{\ho{\s}}_k$ is well defined. 
Furthermore, the quantity $\ho{\s}_{02}^{k-3} \A^{\ho{\s}}_k$ is invariant under 
$w \lra w+ v$. 
\item \label{wd_pd6}If $\s(q)=0$, $\nabla \s|_{q} =0$, $\nabla^2 \s|_{q} =0$
and  $ \nabla^3 \s|_{q}(v,v,\cdot) =0$, then $\nabla^4 \s|_x(v,v,v,v)$ is well defined. 
\item \label{wd_pdk} Let $\D^{\ho{\s}}_k$  be the corresponding sections induced from the quantities 
$\D^{\q}_k$ obtained in \eqref{Dkdefn} by replacing $\q_{ij}$ with $\ho{\s}_{ij}$. 
If $\s(q)=0$, $\nabla \s|_{q} =0$, $\nabla^2 \s|_{q} =0$, 
$ \nabla^3 \s|_q(v,v,\cdot) =0$, $\nabla^4 \s|_q(v,v,v,v) =0$ and 
$ \ho{\s}_{12}^{\epsilon_i} \D^{\ho{\s}}_{i} =0$ for all $i<k$, then  $\ho{\s}_{12}^{\epsilon_k}\D^{\ho{\s}}_{k}$ 
is well defined. Furthermore, the quantity $\ho{\s}_{12}^{\epsilon_k}\D^{\ho{\s}}_{k}$  
is invariant under $w \lra w+ v$. Here as before, $\epsilon_k$ is the order of the pole of $\D^{\q}_k$ at 
$\q_{12} =0$. 
\item \label{wd_pe6} If $\s(q)=0$, $\nabla \s|_{q} =0$, $\nabla^2 \s|_{q} =0$ and $ \nabla^3 \s|_{q}(v,v,\cdot) =0$, then $\nabla^3 \s|_q(v,w,w)$ is well defined. Furthermore, the quantity $\nabla^3 \s|_q(v,w,w)$ is invariant under $w \lra w+ v$.
\item \label{wd_pe7} If $\s(q)=0$, $\nabla \s|_{q} =0$, $\nabla^2 \s|_{q} =0$, 
$ \nabla^3 \s|_q(v,v,\cdot) =0$, and $\nabla^3 \s|_q(v,w,w) =0$, then $\nabla^4 \s|_q(v,v,v,v)$ is well defined. 
\item \label{wd_pe8} If $\s(q)=0$, $\nabla \s|_{q} =0$, $\nabla^2 \s|_{q} =0$, 
$ \nabla^3 \s|_q(v,v,\cdot) =0$, $\nabla^3 \s|_q(v,w,w) =0$ and $\nabla^4 \s|_q(v,v,v,v)=0$, then 
$\nabla^4 \s|_q(v,v,v,w) $ is well defined. Furthermore, the quantity $\nabla^3 \s|_q(v,v,v,w)$
is invariant under $w \lra w+ v$.
\end{enumerate}
\end{lmm}
\pf We omit the details of the proof here; the details can be found in \cite{BM_Detail}. 
These facts follow in a straightforward way by unwinding definition \ref{vertical_derivative_defn}. As explained in remarks \ref{induced_section_Ak} and \ref{induced_section_Dk}, 
Corollary \ref{Ak_node_condition_cor} and \ref{Dk_node_condition_cor} imply Lemma \ref{wdc}, statement \ref{wd_pak} and \ref{wd_pdk}, respectively. \qed

\begin{rem}
Let us mention a pedantic point about our notation. The $~\ho{}~$ introduced in the notation of \eqref{pedantic_point} might seem strange to the reader.  We have done 
that to be consistent with \eqref{abbreviation}. According to our notation,  if $f:\P^2 \lra \gP^{*d}$ is a section and $\p \in \P^2$ then 
\begin{align*}
\ho{f}_{ij} &:=  \nabla^{i+j} f|_{\p} (\underbrace{v,\cdots v}_{\textnormal{$i$ times}}, \underbrace{w,\cdots w}_{\textnormal{$j$ times}})\in \gP^{*d}\big|_{\p} \\
f_{ij} &:= \{ \nabla^{i+j} f|_{\p} (\underbrace{v,\cdots v}_{\textnormal{$i$ times}}, \underbrace{w,\cdots w}_{\textnormal{$j$ times}}) \}(p^{\otimes d}) =  
\nabla^{i+j} f|_{p} (\underbrace{v,\cdots v}_{\textnormal{$i$ times}}, \underbrace{w,\cdots w}_{\textnormal{$j$ times}})\in\C. 
\end{align*}
Since we encounter the second quantity more in our computations, we have denoted that as $f_{ij}$. Notice that if $\ho{f}_{ij}$ is well defined, then so is $f_{ij}$.  
\end{rem}


\section{Transversality}   
\label{bundle_sections}

\hf\hf In this section we give an explicit description of the spaces $\A_0$,
$\A_1$ and $\PP \X_k$ in terms of vanishing and non-vanishing of bundle sections. Before proceeding, let us state three important Lemmas. We will then show that these spaces satisfy the hypothesis of one of these three Lemmas.

\begin{lmm}
\label{cl1}
Let $M$ be a smooth  manifold and $\{S_i\}_{i=0}^k$ be a family of subspaces defined as: 
\begin{align}
\label{cl1_equation1} 
\SS_{i}:= \{ p \in M: \ts_{0}(p) =0, \ldots, \ts_{i}(p) =0,\ts_{i+1}(p) \neq 0 \}  \qquad \forall ~~0 \leq i\leq k 
\end{align}
where 
$\ts_{i}: M \lrab V_{i}$ are sections of vector bundles only defined on the subspace 
$$ \{ p\in M: \ts_0(p), \ldots, \ts_{i-1}(p) =0 \}.$$ 
If the section 
 $$\ts_{i}: \{ p \in M: \ts_0(p), \ldots, \ts_{i-1}(p) =0 \} 
\lra V_{i} $$ 
is transverse to the zero set for $0\leq i\leq k+1$ then
\begin{align}
\ov{\SS}_{k-1} & = \{ p\in M: \ts_0(p)=0, \ldots, \ts_{k-1}(p)=0 \} = \SS_{k-1} \cup \ov{\SS}_{k}.\label{cl1_equation2}
\end{align}
In particular, $\ov{\SS}_{k-1}$ is a smooth manifold.
\end{lmm}
 
\begin{lmm}
\label{cl2}
Let $M$ be a smooth  manifold and $\{\SS_i\}_{i=-1}^k$ be a family of subspaces defined as follows:
\begin{align*} 
\SS_{-1} &= \{ p \in M: \ts_{0}(p) \neq 0, \nts (p) \neq 0 \}, \\ 
\SS_{i} &= \{ p \in M: \ts_{0}(p) =0, \ldots, \ts_{i}(p) =0, 
\ts_{i+1}(p) \neq 0, \nts (p) \neq 0 \} \qquad 
\forall ~~ 0 \leq i \leq k   
\end{align*}
where $\ts_{i}: M \lrab V_{i}$ are sections of vector bundles only defined on the subspace  
$$ \{ p\in M: \ts_0(p), \ldots, \ts_{i-1}(p) =0 \}$$
and $\nts : M \lra W$ is defined. Suppose that 
$$ \ts_{i}: 
\{ p\in M: \ts_0(p) =0, \ldots, \ts_{i-1}(p) =0, \nts(p) \neq 0 
\} \lra V_i $$
is transverse to the zero set for $0\leq i\leq k+1$. Then 
\begin{align}
\label{cl2_eqn1}
\ov{\SS}_{k-1} &= \SS_{k-1} \cup \ov{\SS}_{k} \cup \B \qquad \textnormal{where} \qquad \B := \{p\in \ov{\SS}_{k-1}: ~~\nts(p) =0 \}.
\end{align}
Furthermore, if $\nts: M \lra W$ is transverse to the zero set, then 
\begin{align}
\label{cl2_extra}
M = \SS_{-1} \cup \ov{\SS}_0 \cup \B^{\prime} \qquad \textnormal{where} \qquad  \B^{\prime} := \{p \in M: \nts(p) =0 \}.
\end{align}
\end{lmm}
 
\begin{lmm}
\label{cl3}
Let $M$ be a smooth manifold and let $\SS_0 \subset M$ be a subspace defined by 
\bgd
\SS_0 = \{ p\in M: \ts_{0}(p) =0, \ts_1(p) \neq 0, \nts(p) \neq 0 \}.
\edd
Here $\ts_{i}: M \lra V_{i}$ 
are sections of vector bundles only defined on the subspace  
$$ \{ p\in M: \ts_0(p), \ldots, \ts_{i-1}(p)=0 \}$$
and $\nts : M \lra W$ is defined. If the sections 
$$\ts_0:M \lra V_0, \qquad \nts: \ts_{0}^{-1}(0) \lra W, \qquad \ts_1: \ts_{0}^{-1}(0)- \nts^{-1}(0) \lra V_1   $$
are transverse to the zero set, then 
\begin{align}
\label{cl3_eqn1}
\ov{\SS}_{0} &= \{ p\in M: \ts_0(p) =0 \}.  
\end{align}
In particular, $\ov{\SS}_0$ is a smooth manifold.
\end{lmm}
\textbf{Proofs of Lemma \ref{cl1}, \ref{cl2} and \ref{cl3}:} See Appendix \ref{appendix_proof_of_Lemmas}. \qed  \\

\subsection{Description of spaces as algebraic varieties}
\ni We will now give an explicit description of our spaces as algebraic varieties. Although we are going to state  
a large number of propositions, the statement of these propositions and their proofs are very similar in nature. 
Moreover it makes it much easier to read the remainder of the document (section \ref{closure_of_spaces} onwards). 
Hence we have decided to organized this section in this particular way.  

\begin{prp}
\label{ift_ml}
The space $\A_0$ can be described as 
\begin{align}
\label{a0_equation}
\A_0 & = \{ (\ff, \p) \in \D\times \P^2: \ds_{\A_0}(\ff, \p) =0, ~\ds_{\A_1}(\ff, \p) \neq 0 \}.
\end{align}
Furthermore, the sections of the vector bundles
\begin{align*}
\ds_{\A_0}:\D \times \P^2   \lra \DL_{\A_0}, ~\ds_{\A_1}: \ds_{\A_0}^{-1}(0) \lra \DV_{\A_1}  
\end{align*}
are transverse to the zero set. In particular, $\A_0$ is a smooth manifold of dimension $\delta_d+1$.
\end{prp}

\begin{cor}
\label{a0_cl_vanish_cor}
The space $\ov{\A}_0$ is a smooth manifold of dimension $\delta_d+1$ that can be described as 
\begin{align*}
\ov{\A}_0 & = \{ (\ff, \p) \in \D\times \P^2: \ds_{\A_0}(\ff, \p) =0\}.  
\end{align*}
\end{cor}

\begin{prp}
\label{mll}
The space $\A_1$ can be described as 
\begin{align}
\label{a1_equation}
\A_1 &= \{ (\ff, \p) \in \ov{\A}_0: \ds_{\A_1} (\ff, \p) =0, ~ \ds_{\A_2} (\ff, \p)  \neq 0, ~\ds_{\D_4}(\ff, \p) \neq 0 \}. 
\end{align}
Furthermore, the sections of the vector bundles,
\begin{align*}
\ds_{\A_1}:\ov{\A}_0 \lra \DV_{\A_1}, ~\ds_{\D_4}:\ds_{\A_1}^{-1}(0) \lra \DV_{\D_4}, ~\ds_{\A_2}:\ds_{\A_1}^{-1}(0)- \ds_{\D_4}^{-1}(0)  \lra \DL_{\A_2} 
\end{align*}
are transverse to the zero set if $d \geq 2$\footnote{The section $\ds_{\D_4}(\ff, \p)$ is non-zero is vacuously true if $\ds_{\A_2} (\ff, \p)  \neq 0$. We have stated the proposition in this way 
so that it is clear that our spaces satisfy the hypothesis of Lemma \ref{cl3}.}. In particular, $\A_1$ is a smooth manifold of dimension $\delta_d-1$ if $d\geq 1$. 
\end{prp}

\begin{cor}
\label{a1_cl_vanish_cor}
The space $\ov{\A}_1$ is a smooth manifold of dimension $\delta_d-1$ that can be described as 
\begin{align*}
\ov{\A}_1 & = \{ (\ff, \p) \in \ov{\A}_0: \ds_{\A_1}(\ff, \p) =0\}, \qquad \textit{provided $d \geq 2$.} 
\end{align*}
\end{cor}

\begin{cor}
\label{a2_and_d4_is_smooth}
The spaces $\A_2$ and $\D_4 $ are smooth manifolds, of dimension $\delta_d-2$  and $\delta_d-4$ respectively,  if $d\geq 2$. 
\end{cor}

\begin{prp}
\label{ift_ml_up}
The space $\AA_0$ can be described as 
\begin{align*}
\AA_0 & = \{ (\ff, \lp) \in \D\times \P T\P^2: \us_{\AA_0}(\ff, \p) =0, ~\us_{\AA_1}(\ff, \lp) \neq 0 \}. 
\end{align*}
Furthermore, the sections of the vector bundles
\begin{align*}
\us_{\AA_0}:\D \times \P^2   \lra \UL_{\AA_0}, ~\us_{\AA_1}: \us_{\AA_0}^{-1}(0) \lra \UV_{\A_1}
\end{align*}
are transverse to the zero set.  
\end{prp}

\begin{cor}
\label{a0_up_cl_vanish_cor}
The space $\ov{\AA}_0$ is a smooth manifold of dimension $\delta_d+2$ that can be described  as 
\begin{align*}
\ov{\AA}_0 = \{ (\ff, \lp) \in \D\times \P T\P^2: \us_{\AA_0}(\ff, \lp) =0\}. 
\end{align*}
\end{cor}

\begin{prp}
\label{mll_up}
The space $\AA_1$ can be described as 
\begin{align*}
\AA_1 &= \{ (\ff, \lp) \in \ov{\AA}_0: \us_{\AA_1} (\ff, \lp) =0, ~ \us_{\AA_2} (\ff, \lp)  \neq 0, ~\us_{\DD_4}(\ff, \lp) \neq 0 \}. 
\end{align*}
Furthermore, the sections of the vector bundles
\begin{align*}
\us_{\AA_1}:\ov{\AA}_0 \lra \UV_{\AA_1}, ~\us_{\DD_4}:\us_{\AA_1}^{-1}(0) \lra \UV_{\DD_4}, ~\us_{\AA_2}:\us_{\AA_1}^{-1}(0)- \us_{\DD_4}^{-1}(0)  \lra \UL_{\AA_2} 
\end{align*}
are transverse to the zero set if $d \geq 2\footnote{Again, $\us_{\DD_4}(\ff, \lp) \neq 0$ is vacuously true if $\us_{\AA_2} (\ff, \lp)  \neq 0$.}$.
\end{prp}

\begin{cor}
\label{a1_up_cl_vanish_cor}
The space $\ov{\AA}_1$ is a smooth manifold of dimension $\delta_d$ that can be described as 
\begin{align*}
\ov{\AA}_1 & = \{ (\ff, \lp) \in \ov{\AA}_0: \us_{\AA_1}(\ff, \lp) =0\}, \textit{ provided $d \geq 2$.} 
\end{align*}
\end{cor}

\begin{prp}
\label{A1_sharp_Condition_prp}
The space $\hat{\A}_1^{\#}$ can be described as 
\begin{align}
\label{a1_sharp_up_equation}
\hat{\A}_1^{\#} = \{ ~(\ff, \lp)\in \D \times \P T\P^2: \us_{\AA_0} (\ff, \lp) =0, ~\us_{\AA_1}(\ff, \lp) =0, ~\us_{\PP \A_2}(\ff, \lp) \neq 0 \}.
\end{align}
Furthermore, the sections of the vector bundles
\begin{gather*}
\us_{\AA_0} : \D \times \P T\P^2 \lra \UL_{\AA_0}, \quad \us_{\AA_1} : \us_{\AA_{0}}^{-1} (0) \lra \UV_{\AA_1},  \quad \us_{\PP \A_2} : \us_{\AA_{1}}^{-1} (0) \lra \UV_{\PP \A_2} 
\end{gather*} 
are transverse to the zero set, provided $d \geq 2$. 
\end{prp}

\begin{cor}
\label{a1_sharp_up_cl_vanish_cor}
The space $\ov{\hat{\A}^{\#}_1}$ is a smooth manifold of dimension $\delta_d$ that can be described as 
\begin{align*}
\ov{\hat{\A}^{\#}_1} & = \{ (\ff, \lp) \in \ov{\AA}_0: \us_{\AA_1}(\ff, \lp) =0\}, \textit{ provided $d \geq 2$.}
\end{align*}
\end{cor}

\begin{prp}
\label{A2_Condition_prp}
The space $\PP\A_2$ can be described as 
\begin{align}
\label{pa2_equation}
\PP \A_2 = \{ ~(\ff, \lp)\in \ov{\hat{\A}_1^{\#}} : ~\us_{\PP \A_2}(\ff, \lp) =0, ~\us_{\PP \A_3}(\ff,\lp) \neq 0, ~\us_{\PP \D_4}(\ff,\lp) \neq 0 \}.
\end{align}
Furthermore, the sections of the vector bundles
\begin{align*}
\us_{\PP \A_2}:  \ov{\hat{\A}_1^{\#}} \lra \UV_{\PP \A_2}, \quad 
\us_{\PP \A_3} : \us_{\PP \A_{2}}^{-1} (0) \lra \UL_{\PP \A_3}, \quad \us_{\PP \D_4} : \us_{\PP \A_{2}}^{-1} (0) \lra \UL_{\PP \D_4} 
\end{align*} 
are transverse to the zero set,   
provided $d \geq 3$.  
\end{prp}

\begin{cor}
\label{pa2_cl_vanish_cor}
The space $\ov{\PP\A}_2$ is a manifold of dimension $\delta_d-2$ and can be described as 
\begin{align*}
\ov{\PP\A}_2 = \{ ~(\ff, \lp)\in \ov{\hat{\A}_1^{\#}} : ~\us_{\PP \A_2}(\ff, \lp) =0\}, \textit{ provided $d \geq 3$.}
\end{align*}
\end{cor}

\begin{prp}
\label{D4_sharp_Condition_prp}
The space $\hat{\D}_4^{\#}$ can be described as 
\begin{align}
\label{d4_sharp_up_equation}
\hat{\D}_4^{\#} = \{(\ff, \lp)\in \D \times \P T\P^2: &  
\us_{\AA_0} (\ff, \lp) =0, \us_{\AA_1}(\ff, \lp) =0, \us_{\DD_4}(\ff,\lp) =0, \us_{\PP \A_3}(\ff,\lp) \neq 0 \}. 
\end{align}
Furthermore, the sections of the vector bundles
\begin{align*}
\us_{\AA_0}: \D \times \P T\P^2 \lra \UL_{\AA_0}, ~\us_{\AA_1} : \us_{\AA_{0}}^{-1} (0) \lra \UV_{\AA_1}, ~\us_{\DD_4} : \us_{\AA_{1}}^{-1} (0) \lra \UV_{\DD_4},  
~\us_{\PP \A_3} : \us_{\DD_4}^{-1} (0) \lra \UL_{\PP \A_3}
\end{align*} 
are transverse to the zero set,   
provided $d \geq 3$.  
\end{prp}

\begin{cor}
\label{d4_up_sharp_cl_vanish_cor}
The space $\ov{\hat{\D}^{\#}_4}$ is a manifold of dimension $\delta_d-4 $ and can be described as 
\begin{align*}
\ov{\hat{\D}_4^{\#}} = \{ (\ff, \lp)\in \D \times \P T\P^2: &  
~\us_{\AA_0} (\ff, \lp) =0, ~\us_{\AA_1}(\ff, \lp) =0, ~\us_{\DD_4}(\ff,\lp) =0 \}, \textit{ provided $d \geq 3$.}
\end{align*}
\end{cor}

\begin{prp}
\label{D4_Condition_prp}
The space $\PP\D_4$ can be described as 
\begin{align}
\label{pd4_equation}
\PP \D_4 = \{ (\ff, \lp)\in \ov{\PP \A}_2 :&  
~\us_{\PP \A_3}(\ff, \lp) =0, 
~\us_{\PP \D_4} (\ff,\lp)=0, ~\us_{\PP \D_5}^{\UL}(\ff,\lp) \neq 0, ~\us_{\PP \D_5^{\vee}}(\ff,\lp) \neq 0 \}.
\end{align}
Furthermore, the sections of the vector bundles 
\begin{align*}
\us_{\PP \A_3}&: \ov{\PP \A}_2 \lra \UL_{\PP \A_3}, \qquad \us_{\PP \D_4}: \us_{\PP \A_3}^{-1} (0) \lra \UL_{\PP \D_4},  \\
\us_{\PP \D_5}^{\UL}&: \us_{\PP \D_4}^{-1} (0) \lra \UL_{\PP \D_5}, \qquad \us_{\PP \D_5^{\vee}}: \us_{\PP \D_4}^{-1} (0)- \us_{\PP \D_5^{\vee}}^{-1}(0) \lra \UL_{\PP \D_5} 
\end{align*} 
are transverse to the zero set, provided $d \geq 3$.
\end{prp}

\begin{cor}
\label{pd4_cl_vanish_cor}
The space $\ov{\PP\D}_4$ is a manifold of dimension $\delta_d-3 $ and can be described as 
\begin{align*}
\ov{\PP\D}_4 = \{(\ff, \lp)\in \ov{\PP \A}_2 :&  
~\us_{\PP \A_3}(\ff, \lp) =0, 
~\us_{\PP \D_4} (\ff,\lp)=0 \}, \textit{~provided $d \geq 3$.}
\end{align*}
\end{cor}

\begin{prp}
\label{A3_Condition_prp}
The space $\PP\A_3$ can be described as 
\begin{align}
\label{pa3_equation}
\PP \A_3 = \{(\ff, \lp)\in \ov{\PP \A}_2:~\us_{\PP \A_3}(\ff, \lp) =0, ~\us_{\PP \A_{4}}(\ff,\lp) \neq 0, ~\us_{\PP \D_4} (\ff, \lp) \neq 0 \}.
\end{align}
Furthermore, the sections of the vector bundles 
\begin{align*}
\us_{\PP \A_3}:  \ov{\PP \A}_2 \lra \UL_{\PP \A_3}, \quad \us_{\PP \D_4} : \us_{\PP \A_{3}}^{-1} (0) \lra \UL_{\PP \D_4}, \quad 
\us_{\PP \A_4} : \us_{\PP \A_{3}}^{-1} (0) - \us_{\PP \D_4}^{-1}(0) \lra \UL_{\PP \A_4}
\end{align*} 
are transverse to the zero set provided $d\geq 4$.
\end{prp}

\begin{cor}
\label{pa3_cl_vanish_cor}
The space $\ov{\PP\A}_3$ is a manifold of dimension $\delta_d-3$ and can be described  as 
\begin{align*}
\ov{\PP \A}_3 = \{(\ff, \lp)\in \ov{\PP \A}_2:~\us_{\PP \A_3}(\ff, \lp) =0\}, \textit{~provided $d\geq 4$.}
\end{align*}
\end{cor}

\begin{cor}
\label{a3_is_smooth}
The space $\A_3$ is a smooth manifold of dimension $\delta_d-3$, if $d\geq 3$. 
\end{cor}

\begin{prp}
\label{Ak_Condition_prp}
If $k > 3$, the space $\PP\A_k$ can be described as 
\begin{align}
\label{pak_equation}
\PP \A_k = \{(\ff, \lp)\in \ov{\PP \A}_3: \us_{\PP \A_j}(\ff, \lp) =0\,\,\textup{for}\,\,4\leq j\leq k, \us_{\PP \A_{k+1}}(\ff,\lp) \neq 0, \us_{\PP \D_4} (\ff, \lp) \neq 0   \}.
\end{align}
Furthermore, the sections of the vector bundles 
$\us_{\PP \A_i}: \us_{\PP\A_{i-1}}^{-1} (0) - \us_{\PP \D_4}^{-1}(0) \lra \UL_{\PP \A_i}$
are transverse to the zero set for all $4\leq i \leq k+1$, provided $d \geq k+1$.
\end{prp}

\begin{cor}
\label{ak_is_smooth}
The space $\A_k$ is a smooth manifold of dimension $\delta_d-k$, if $d \geq k$. 
\end{cor}

\begin{prp}
\label{D5_Condition_prp}
The space $\PP\D_5$ can be described as 
\begin{align}
\label{pd5_equation}
\PP \D_5 = \{(\ff, \lp)\in \ov{\PP \D}_4:  
~\us_{\PP \D_5}^{\UL} (\ff,\lp)=0 , ~ \us_{\PP \D_6}(\ff,\lp) \neq 0, ~\us_{\PP \E_6}(\ff,\lp) \neq 0 \}
\end{align}
Furthermore, the sections of the vector bundles 
\begin{align*}
\us_{\PP \D_5^{\UL}}: \ov{\PP \D}_4: \lra \UL_{\PP \D_5},  
~~\us_{\PP \D_6}: \us_{\PP \D_5^{\UL}}^{-1} (0) \lra \UL_{\PP \D_6}, 
~~\us_{\PP \E_6} : \us_{\PP \D_5^{\UL}}^{-1} (0) \lra \UL_{\PP \E_6}
\end{align*} 
are transverse to the zero set, provided $d \geq 4$.
\end{prp}

\begin{cor}
\label{pd5_cl_vanish_cor}
The space $\ov{\PP\D}_5$ is a manifold of dimension $\delta_d-5$  and can be described as 
\begin{align*}
\ov{\PP \D}_5 = \{(\ff, \lp)\in \ov{\PP \D}_4:~\us_{\PP \D_5}^{\UL} (\ff,\lp)=0 \}, \textit{~provided $d \geq 4$.}
\end{align*}
\end{cor}

\begin{prp}
\label{D5_Condition_prp_dual}
The space $\PP\D_5^{\vee}$ can be described as 
\begin{align}
\label{pd5_equation_dual}
\PP \D_5^{\vee} = \{(\ff, \lp)\in \ov{\PP \D}_4:  
~\us_{\PP \D_5^{\vee}} (\ff,\lp)=0 , ~ \us_{\PP \D_6^{\vee}}(\ff,\lp) \neq 0, ~\us_{\PP \D_5}(\ff,\lp) \neq 0 \}.
\end{align}
Furthermore, the sections of the vector bundles 
\begin{align*}
\us_{\PP \D_5}^{\UL}: \ov{\PP \D}_4: \lra \UL_{\PP \D_5},  
~~\us_{\PP \D_5^{\vee}}: \ov{\PP \D}_4 - \us_{\PP \D_5}^{\UL^{-1}}(0) \lra \UL_{\PP \D_6^{\vee}}, 
~~\us_{\PP \D_6^{\vee}} : \us_{\PP \D_5^{\vee}}^{-1} (0) \lra \UL_{\PP \D_6^{\vee}}
\end{align*} 
are transverse to the zero set, provided $d \geq 4$.
\end{prp}

\begin{cor}
\label{pd5_cl_vanish_cor_dual}
The space $\ov{\PP\D^{\vee}_5}$ is a manifold of dimension $\delta_d-5$  and can be described as 
\begin{align*}
\ov{\PP \D^{\vee}_5} = \{(\ff, \lp)\in \ov{\PP \D}_4:~\us_{\PP \D_5^{\vee}} (\ff,\lp)=0 \}, \textit{~provided $d \geq 4$.}
\end{align*}
\end{cor}

\begin{cor}
\label{d5_is_smooth}
The space $\D_5$ is a smooth manifold of dimension $\delta_d-5$, if $d \geq 3$. 
\end{cor}

\begin{prp}
\label{E6_Condition_prp}
The space $\PP \E_6 $ can be described as
\begin{align}
\label{pe6_equation}
\PP \E_6 = \{ (\ff,\lp) \in \ov{\PP \D}_5 &: \us_{\PP \E_6}(\ff,\lp) =0, 
~\us_{\PP \E_7}(\ff, \lp) \neq 0, ~\us_{\PP \XC_8}(\ff, \lp) \neq 0  \}
\end{align}
Furthermore, the sections of the vector bundles 
\begin{align*}
\us_{\PP \E_6} : \ov{\PP \D}_5 \lra \UL_{\PP \E_6}, \qquad \us_{\PP \E_7}: \us_{\PP \E_6}^{-1} (0) \lra \UL_{\PP \E_7}, 
\qquad \us_{\PP \XC_8}: \us_{\PP \E_6}^{-1} (0) \lra \UL_{\PP \XC_8}
\end{align*} 

%
are transverse to the zero set, provided $d \geq 4$.
\end{prp}

\begin{cor}
\label{pe6_cl_vanish_cor}
The space $\ov{\PP \E}_6$ is a  manifold of dimension $\delta_d-6$ and can be described as
\begin{align*}
\ov{\PP \E}_6 = \{ (\ff,\lp) \in \ov{\PP \D}_5 &: \us_{\PP \E_6}(\ff,\lp) =0\}, \textit{~provided $d \geq 4$.}
\end{align*}
\end{cor}

\begin{cor}
\label{e6_is_smooth}
The space $\E_6$ is a  manifold  of dimension $\delta_d-6$, if $d \geq 3$. 
\end{cor}

\begin{prp}
\label{D6_Condition_prp}
The space $\PP\D_6$ can be described as 
\begin{align}
\label{pd6_equation}
\PP \D_6 = \{(\ff, \lp)\in \ov{\PP \D}_5: \us_{\PP \D_6}(\ff, \lp) =0, ~\us_{\PP \D_7}(\ff,\l) \neq 0, ~\us_{\PP \E_6} (\ff,\lp) \neq 0 \}.
\end{align}
Furthermore, the sections of the vector bundles 
\begin{align}
\us_{\PP \D_6}:  \ov{\PP \D}_5 \lra \UL_{\PP \D_6}, \quad \us_{\PP \E_6} : \us_{\PP \D_6}^{-1} (0) \lra \UL_{\PP \E_6}, \quad 
\us_{\PP \D_7} : \us_{\PP \D_6}^{-1} (0) - \us_{\PP \E_6}^{-1}(0) \lra \UL_{\PP \D_7}
\end{align} 
are transverse to the zero set, provided $ d \geq 5$. 
\end{prp}

\begin{cor}
\label{pd6_cl_vanish_cor}
The space $\ov{\PP\D}_6$  is  a manifold of dimension $\delta_d-6$ and can be described as  
\begin{align*}
\ov{\PP \D}_6 = \{ ~(\ff, \lp)\in \ov{\PP \D}_5: \us_{\PP \D_6}(\ff, \lp) =0 \}, \textit{~provided $d \geq 5$.}
\end{align*}
\end{cor}

\begin{cor}
\label{d6_is_smooth}
The space $\D_6$ is a manifold of dimension $\delta_d-6$, if $d \geq 4$. 
\end{cor}

\begin{prp}
\label{Dk_Condition_prp}
If $k > 6$, then the space $\PP\D_k$ can be described as 
\begin{align}
\label{pdk_equation}
\PP \D_k = \{(\ff, \lp)\in \ov{\PP \D}_6: \us_{\PP \D_j}(\ff, \lp) =0\,\,\textup{if}\,\, 7\leq j\leq k, 
 ~\us_{\PP \D_{k+1}}(\ff,\lp) \neq 0, ~\us_{\PP \E_6} (\ff,\lp) \neq 0   \}.
\end{align}
Furthermore, the section of the vector bundle $ ~\us_{\PP \D_i}: \us_{\PP \D_{i-1}}^{-1} (0) - \us_{\PP \E_6}^{-1}(0) \lra \UL_{\PP \D_i} ~$
is transverse to the zero set for all $ 7 \leq i \leq k+1$ 
provided $d >k-2$.
\end{prp}

\begin{cor}
\label{dk_is_smooth}
The space $\D_k$ is a manifold of dimension $\delta_d-k$, if $d \geq k-2$. 
\end{cor}

\begin{prp}
\label{E7_Condition_prp}
The space $\PP \E_7 $ can be described as 
\begin{align}
\label{pe7_equation}
\PP \E_7 = \{ (\ff,\lp) \in \ov{\PP \E}_6: ~~\us_{\PP \E_7}(\ff,\lp) =0, 
 ~~\us_{\PP \E_8}(\ff, \lp) \neq 0, ~~\us_{\PP \XC_8}(\ff, \lp) \neq 0  \}.
\end{align}
Furthermore, the sections of the vector bundles 
\begin{align*}
\us_{\PP \E_7}:  \ov{\PP \E}_6 \lra \UL_{\PP \E_7}, \quad \us_{\PP \E_8} : \us_{\PP \E_7}^{-1} (0) \lra \UL_{\PP \E_8}, \quad 
\us_{\PP \XC_8} : \us_{\PP \E_7}^{-1} (0) \lra \UL_{\PP \XC_8}
\end{align*} 
are transverse to the zero set, provided $ d \geq 4$. 
\end{prp}

\begin{cor}
\label{pe7_cl_vanish_cor}
The space $\ov{\PP \E}_7$ is a manifold of dimension $\delta_d-6$ and can be described as 
\begin{align*}
\ov{\PP \E}_7 = \{ (\ff,\lp) \in \ov{\PP \E}_6: ~~\us_{\PP \E_7}(\ff,\lp) =0 \}, \textit{~provided $d \geq 4$.}
\end{align*}
\end{cor}

\begin{cor}
\label{e7_is_smooth}
The space $\E_7$ is a manifold of dimension $\delta_d-7$, if $d \geq 4$. 
\end{cor}

\begin{prp}
\label{E7_Condition_prp_using_D6}
The space $\PP \E_7 $ can also be described as 
\begin{align}
\label{pe7_equation_also}
\PP \E_7 = \{ (\ff,\lp) \in \ov{\PP \D}_6: ~~\us_{\PP \E_6}(\ff,\lp) =0, 
 ~~\us_{\PP \E_8}(\ff, \lp) \neq 0, ~~\us_{\PP \XC_8}(\ff, \lp) \neq 0  \}.
\end{align}
Furthermore, the sections of the vector bundles 
\begin{align*}
\us_{\PP \E_6}:  \ov{\PP \D}_6 \lra \UL_{\PP \E_6}, \quad \us_{\PP \E_8} : \us_{\PP \E_7}^{-1} (0) \lra \UL_{\PP \E_8}, \quad 
\us_{\PP \XC_8} : \us_{\PP \E_7}^{-1} (0) \lra \UL_{\PP \XC_8}
\end{align*} 
are transverse to the zero set, provided $ d \geq 4$. 
\end{prp}

\begin{cor}
\label{pe7_cl_vanish_cor_also}
The space $\ov{\PP \E}_7$ is a manifold of dimension $\delta_d-6$ and can also be described as 
\begin{align*}
\ov{\PP \E}_7 = \{ (\ff,\lp) \in \ov{\PP \D}_6: \us_{\PP \E_6}(\ff,\lp) =0 \}, \textit{provided $ d \geq 4$.}
\end{align*}
\end{cor}

\begin{prp}
\label{X8_sharp_Condition_prp}
The spaces $\hat{\XC}_8^{\#}$ and  $\hat{\XC}_8^{\# \flat}$ can be described as 
\begin{align}
\label{x8_sharp_and_flat_equation}
\hat{\XC}_8^{\#} & = \{(\ff, \lp)\in \D \times \P T\P^2: \us_{\AA_0} (\ff, \lp) =0,\us_{\AA_1}(\ff, \lp) =0, \us_{ \DD_4}(\ff,\lp) =0, \nonumber \\ 
& \qquad \qquad \qquad \qquad  \qquad \,\,\,\us_{\XX_8}(\ff,\lp) =0, \us_{\PP \E_7} (\ff,\lp)\neq 0\} \nonumber \\
\hat{\XC}_8^{\# \flat} & = \{(\ff, \lp)\in \hat{\XC}_8^{\#}:  \us_{\J}(\ff,\lp) \neq 0\}.  
\end{align}
Furthermore, the sections of the vector bundle
\begin{align*}
\us_{\AA_0}& : \D \times \P T\P^2 \lra \UL_{\AA_0}, \qquad \us_{\AA_1} : \us_{\AA_{0}}^{-1} (0) \lra \UV_{\AA_1}, \qquad \us_{\DD_4}: \us_{\AA_{1}}^{-1} (0) \lra \UV_{\DD_4} \\
\us_{\XX_8}&: \us_{\DD_4}^{-1} (0) \lra \UV_{\XX_8}, \qquad \us_{\PP\E_7}: \us_{\XX_8}^{-1}(0) \lra \UL_{\PP \E_7}, \qquad \us_{\J}: \us_{\XX_8}^{-1}(0)- \us_{\PP \E_7}^{-1}(0)  \lra \UL_{\J} 
\end{align*} 
are transverse to the zero set,   
provided $d \geq 4$.  
\end{prp}

\begin{cor}
\label{X8_and_X8_sharp_up_cl_vanish_cor}
The spaces $\ov{\hat{\XC}_8^{\#}}$ and $\ov{\hat{\XC}_8^{\# \flat}}$ are equal to each other. They are manifolds of dimension $\delta_d-7$ and can be described as 
\begin{align*}
\ov{\hat{\XC}_8^{\#}} = \ov{\hat{\XC}_8^{\# \flat}} = \{ (\ff, \lp)\in \D \times \P T\P^2:    
\us_{\AA_0} (\ff, \lp) =0,\us_{\AA_1}(\ff, \lp) =0, \us_{ \DD_4}(\ff,\lp) =0, \us_{\XX_8}(\ff,\lp)=0  \}
\end{align*}
provided $ d \geq 4$.
\end{cor}

\subsection{Proofs of Propositions} 

\hf\hf Let $\F \cong \C^{\delta_d+1}$ denote the space of homogeneous polynomials of degree $d$. Let $\F^*$ be the subspace of non-zero polynomials. 
This can also be thought of as the space of polynomials in two variables of degree at most $d$. 
If $V \lra M$ is any vector bundle then a section 
$$\psi: \D \times M \lra \pi_{\D}^*\gD^* \otimes \pi_M^* V  $$
induces a section 
\begin{align}
\label{polynomial_induced_section}
\ha{\psi}: \F^* \times M \lra \pi_{M}^*V     \qquad \textnormal{given by}  \qquad \ha{\psi}(f, x):= \{\psi(\ff, x)\}(f).
\end{align} 
We note that $\psi$ is transverse to zero at $(\ff, x)$ if and only if $\ha{\psi}$ is transverse to  zero  at $(f,x)$.\\

\ni \textbf{Proof of Proposition \ref{ift_ml}:}   \ni Equation \eqref{a0_equation} follows from  Corollary \ref{A0_node_condition_cor}. We will now show transversality. We will use the setup of definition \ref{vertical_derivative_defn}, 
remark \ref{transverse_local} and \eqref{polynomial_induced_section}. Suppose $(\ff,\p)\in\ds_{\A_0}^{-1}(0)$. Choose homogeneous coordinates $[X:Y:Z]$ on $\P^2$ so that $\p =[0:0:1]$  and let 
\begin{align*}
\U & :=\big\{[X:Y:Z] \in \P^2: Z\neq0\big\}, ~~\varphi_{\U} : \U \lra \C^2,\,\,~~\varphi_{\U}([X:Y:Z]) = (X/Z,Y/Z).  
\end{align*}
Let us  also denote $x:=X/Z$   and $y:=Y/Z$. The section
\bgd
\ds_{\A_0} : \D \times \P^2 \lra \DL_{\A_0}:=\pi_{\D}^\ast\gD \otimes \pi_{\P^2}^\ast \gamma_{\P^2}^\ast
\edd
induces a section $\ha{\ds}_{\A_0} : \F^* \times \P^2 \lra \pi_{\P^2}^\ast \gamma_{\P^2}^\ast$ as given in \eqref{polynomial_induced_section}. We shall denote the bundle by $\gamma_{\P^2}^\ast$. We now observe that with respect to the \textit{standard} trivialization of $\gamma_{\P^2}^\ast \lra \F^* \times \P^2 $, the 
induced map $ \hatii{\ha{\ds}}_{\A_0}:\F^* \times \C^2 \lra\C$ of the section $\ha{\ds}_{\A_0}$ (cf. \eqref{section_local_coordinate_calculus}, remark \ref{transverse_local}) is 
given by 
$$\hatii{\ha{\ds}}_{\A_0}(\qq,x,y) = \qq(x,y):= \qq_{00} + \qq_{10} x + \qq_{01} y + \frac{f_{20}}{2} x^2 + f_{11} x y + \frac{f_{02}}{2} y^2 +\cdots.$$
By remark \ref{transverse_local}, 
it suffices to show that this induced map $ \hatii{\ha{\ds}}_{\A_0}$ 
is transverse to zero  at~$(\qq,0,0)$. Since the Jacobian matrix of this map at $(\qq,0,0)$ is  
\bgd
d \hatii{\ha{\ds}}_{\A_0}\big|_{(\qq,0,0)} = (1 \,\,\,\, 0 \,\,\,\, 0 \,\,\,\,\cdots ),
\edd
where the first column is partial derivative with respect to $\qq_{00}$, transversality follows. Next we will prove that $\ds_{\A_1}$ restricted to $\ds_{\A_0}^{-1}(0)$ is transverse to zero. 
By Lemma \ref{wdc}, statement \ref{wd_nabla_all} we conclude that if $\ds_{\A_0}(\ff, \p) =0$, 
then $\ds_{\A_1}(\ff, \p)$ is well defined. Let   
\bgd
(\ff,\p)\in \ds_{\A_1}^{-1}(0) \subset  \ds_{\A_0}^{-1}(0).
\edd 
With respect to the \textit{standard} trivialization of 
$ \pi_{D}^* \gD \otimes \DV_{\A_1} \lra \ha{\ds}_{\A_0}^{-1}(0)$, the 
induced map of the section $\hatii{\ha{\ds}}_{\A_1}$ (cf. definition \ref{vertical_derivative_defn}, \eqref{section_local_coordinate}) is given by
\bgd
\hatii{\ha{\ds}}_{\A_1}: (\F^* \times \C^2) \cap \hatii{\ha{\ds}}_{\A_0}^{-1}(0)  \lra\C^2, \qquad \hatii{\ha{\ds}}_{\A_1}( \qq, x,y) = \big(\qq_x(x,y), ~\qq_y(x,y)\big).
\edd
Since the function $\hatii{\ha{\ds}}_{\A_0}$ is transverse to the zero set at $(\qq, 0,0)$, showing that  $\hatii{\ha{\ds}}_{\A_1}$ is transverse to zero at $(\qq,0,0)$ is equivalent to showing that the map
\bgd
\hatii{\ha{\ds}}_{\A_0} \oplus \hatii{\ha{\ds}}_{\A_1}: \F^* \times \C^2 \lra\C^3, \qquad \hatii{\ha{\ds}}_{\A_0} \oplus \hatii{\ha{\ds}}_{\A_1}(\qq, x,y) = \big(\qq(x,y), \qq_x(x,y), \qq_y(x,y)\big)
\edd
is transverse to zero. Since $\qq(x,y) = \qq_{00} + \qq_{10} x + \qq_{01} y + \cdots$, the Jacobian at $(\qq,0,0)$ is 
\[ d (\hatii{\ha{\ds}}_{\A_0} \oplus \hatii{\ha{\ds}}_{\A_1}) |_{(\qq,0,0)} = \left( \begin{array}{cccc}
1 & 0 & 0 & \ldots  \\
0 & 1 & 0 & \ldots  \\
0 & 0 & 1 & \ldots  \end{array} \right),\]
where the first three columns are partial derivatives
with respect to $\qq_{00}$, $\qq_{10}$ and $\qq_{01}$.\qed \\

\ni \textbf{Proof of Corollary \ref{a0_cl_vanish_cor}:} 
\ni This follows immediately from Lemma \ref{cl1} and Proposition \ref{ift_ml}. \qed \\ 

\ni \textbf{Proof of Proposition \ref{mll}:} 
\ni Equation \eqref{a1_equation} follows from Corollary \ref{A1_node_condition_cor} and Corollary \ref{a0_cl_vanish_cor}. We have already shown the transversality of $\ds_{\A_1}$ in the proof of Proposition \ref{ift_ml}. We will now show the transversality of $\ds_{\D_4}$ and $\ds_{\A_2}$. Let us start with $\ds_{\D_4}$. By Lemma \ref{wdc}, statement \ref{wd_nabla_all} we conclude that $\ds_{\D_4}$ is well defined restricted to $\ds_{\A_1}^{-1}(0)$. As in the proof of Proposition \ref{ift_ml}, proving transversality is equivalent to showing that the map
\begin{align*}
\hatii{\ha{\ds}}_{\A_0}\oplus \hatii{\ha{\ds}}_{\A_1} \oplus \hatii{\ha{\ds}}_{\D_4}:\F^* \times \C^2 \lra\C^6, \qquad 
(\qq,x,y)\mapsto \big(\qq(x,y), \qq_x, \qq_y, \qq_{xx}, \qq_{xy}, \qq_{yy}\big)
\end{align*}
is transverse to zero at the point $(\qq,0,0)$. The Jacobian matrix of this map at $(\qq,0,0)$ is a $6 \times (\delta_d+3)$ matrix which has full rank if $d\geq 2$; this follows, for instance, if the first six columns of the matrix are partial derivatives with respect to $\qq_{00}$, $\qq_{10}$, $\qq_{01}$, $\qq_{20}$, $\qq_{11}$ and $\qq_{02}$.  \\
\hf Next let us show transversality of the section $\ds_{\A_2}$. The section is well defined on $\ds_{\A_1}^{-1}(0)$ by Lemma \ref{wdc}, statement \ref{wd_nabla_all}.    
As before, proving transversality  is equivalent to showing that the map
\begin{align*}
\hatii{\ha{\ds}}_{\A_0}\oplus \hatii{\ha{\ds}}_{\A_1} \oplus \hatii{\ha{\ds}}_{\A_2}:\F^*\times \C^2 \lra\C^4, \qquad 
(\qq,x,y)\mapsto \big(\qq(x,y), ~\qq_x, ~\qq_y, ~\qq_{xx} \qq_{yy} - \qq_{xy}^2\big) 
\end{align*}
is transverse to zero at the point $(\qq, 0,0)$.
The Jacobian matrix of this map at $(\qq,0,0)$ is a $4 \times (\delta_d+3)$ matrix which has full rank if $d\geq 2$ and $\ds_{\D_4}(\ff, \p) \neq 0$ : take the first four columns of the matrix as the partial derivatives with respect to $\qq_{00}$, $\qq_{10}$, $\qq_{01}$, and $\qq_{20}$, $\qq_{02}$ or $\qq_{11}$, depending on whether $\qq_{02}$, $\qq_{20}$ or $\qq_{11}$ is non-zero. One of them is guaranteed to be non-zero if $\ds_{\D_4}(\ff, \p) \neq 0$. \qed \\

\ni \textbf{Proof of Corollary \ref{a1_cl_vanish_cor}:} This follows immediately from Lemma \ref{cl3} and Proposition \ref{mll}. \qed \\

\textbf{Proof of Corollary \ref{a2_and_d4_is_smooth}:} This follows from Proposition \ref{mll}, Corollary \ref{Ak_node_condition_cor} and Corollary \ref{D4_node_condition_cor}.  \qed \\ 

\ni \textbf{Proof of Proposition \ref{ift_ml_up}:}  This is identical to the proof of Proposition \ref{ift_ml}.  \qed \\

\ni \textbf{Proof of Corollary \ref{a0_up_cl_vanish_cor}:} This follows immediately from Lemma \ref{cl1} and Proposition \ref{ift_ml_up}. \qed \\

\ni \textbf{Proof of Proposition \ref{mll_up}:}  This is identical to the proof of Proposition \ref{mll}.  \qed \\

\ni \textbf{Proof of Corollary \ref{a1_up_cl_vanish_cor}:} This follows immediately from Lemma \ref{cl3} and Proposition \ref{mll_up}. \qed \\

\ni \textbf{Proof of Proposition \ref{A1_sharp_Condition_prp}:} Equation \eqref{a1_sharp_up_equation} is the \textit{definition} of $\AA_1^{\#}$, so there is nothing to prove. To show transversality we continue with the setup of the proof of Proposition \ref{ift_ml}, but choose coordinate chart so that 
\begin{align*}
\tilde{\U} & := \{ [a \partial_x, ~b \partial_y] \in \P T\P^2|_{\U}: a \neq 0 \}, ~~\varphi_{\tilde{\U}} : \tilde{\U} \lra \C^3,\,\,\varphi_{\tilde{\U}} ([a \partial_x, ~b\partial_y]) = (x,y, \eta),  
\end{align*}
where $\eta:= b/a$. By Lemma \ref{wdc},  statement \ref{wd_nabla_all} we conclude that $\us_{\PP \A_2}(\ff, \lp)$ is well defined restricted to $\us_{\hat{\A}_1}^{-1}(0)$. 
With respect to the \textit{standard} trivialization, the induced map $ \hatii{\ha{\us}}_{\PP \A_2}$ restricted to $\hatii{\ha{\us}}_{\hat{\A}_1}^{-1}(0) $  is given by
\begin{align*} 
\hatii{\ha{\us}}_{\PP \A_2}: (\F^* \times \C^3) \cap \hatii{\ha{\us}}_{\hat{\A}_1}^{-1}(0) \lra\C^2, \qquad  \qq, x,y, \eta\mapsto \big(\qq_{xx} + \eta \qq_{xy}, ~\qq_{xy} + \eta \qq_{yy}\big).
\end{align*}
As before, this is equivalent to showing that the map 
\begin{align*}
\hatii{\ha{\us}}_{\A_0} \oplus \hatii{\ha{\us}}_{\hat{\A}_1} \oplus \hatii{\ha{\us}}_{\PP \A_2}: \F^* \times \C^3 \lra\C^5, \quad 
(\qq, x,y, \eta)\mapsto \big(\qq(x,y), \qq_x, \qq_y, \qq_{xx} + \eta \qq_{xy}, \qq_{xy} + \eta \qq_{yy}\big)
\end{align*}
is transverse to zero at $(\qq, 0,0,0)$. The Jacobian matrix of this map at $(\qq,0,0,0)$ is a $5 \times (\delta_d+4)$ matrix which has full rank if $d\geq 2$. This is easy to see if the first five columns of the matrix are partial derivatives with respect to   $\qq_{00}$, $\qq_{10}$, $\qq_{01}$, $\qq_{20}$ and $\qq_{11}$.  \qed \\

\ni \textbf{Proof of Corollary \ref{a1_sharp_up_cl_vanish_cor}:} This follows immediately from Lemma \ref{cl1} and Proposition \ref{A1_sharp_Condition_prp}.\qed \\ 

\ni \textbf{Proof of Proposition \ref{A2_Condition_prp}:} Equation \eqref{pa2_equation} follows from Lemma \ref{fstr_prp} and Corollary \ref{Ak_node_condition_cor}. To see this, we observe that Corollary \ref{Ak_node_condition_cor} (for $k=2$) gives us a necessary and sufficient condition for a curve $\q^{-1}(0)$ to have an $\A_2$-node. In terms of bundle sections, this is equivalent to the statement that the space $\PP\A_2$ can be described as 
\begin{align*}
\PP \A_2 = \{ ~(\ff, \lp)\in \D \times \P T\P^2 :&  
 ~\us_{\AA_0} (\ff, \lp) =0, ~\us_{\AA_1}(\ff, \lp) =0, 
~\us_{\PP \A_2}(\ff, \lp) =0, \\ 
& ~\us_{\PP \A_3}(\ff,\lp) \neq 0, ~\us_{\PP \D_4}(\ff,\lp) \neq 0 \}.
\end{align*}
The desired equation \eqref{pa2_equation} now follows from Corollary \ref{a0_up_cl_vanish_cor} and \ref{a1_sharp_up_cl_vanish_cor}. \\
\hf We will now show transversality, having already proved it for the section $\us_{\PP \A_2}$ in the proof of Proposition \ref{A1_sharp_Condition_prp}. Let us start with $\us_{\PP \A_3}$. 
By Lemma \ref{wdc}, statement \ref{wd_pa3} we conclude that restricted to $\us_{\PP \A_2}^{-1}(0)$, the section $\us_{\PP \A_3} $ is well defined. As before, showing that the section $\us_{\PP \A_3} $ is transverse to the zero set is equivalent to showing that the map 
\bgd
\hatii{\ha{\us}}_{\AA_0} \oplus \hatii{\ha{\us}}_{\AA_1} \oplus \hatii{\ha{\us}}_{\PP \A_2} \oplus \hatii{\ha{\us}}_{\PP \A_3} : \F^* \times \C^3 \lra\C^6,
\edd
\bgd
(\qq, x,y, \eta)\mapsto \big(\qq(x,y), ~\qq_x, ~\qq_y, ~\qq_{xx} + \eta \qq_{xy}, ~\qq_{xy} + \eta \qq_{yy}, ~\qq_{\xx \xx \xx}\big)
\edd
is transverse to zero, where 
\bgd
\xx := x + \eta y \qquad \textnormal{and} \qquad \qq_{\xx \xx \xx} := (\partial_{x} + \eta \partial_{y})^3 f(x,y).
\edd
The Jacobian matrix of this map at $(\qq,0,0,0)$ is a $6 \times (\delta_d+4)$ matrix which has full rank if $d\geq 3$. This follows by looking at the first six columns of the matrix which are partial derivatives with respect to $\qq_{00}$, $\qq_{10}$,  $\qq_{01}$, $\qq_{20}$, 
$\qq_{11}$ and $\qq_{30}$.   Next let us show the transversality of $\us_{\PP \D_4}$. 
Note that $\us_{\PP \D_4}$ is well defined restricted to $\us_{\PP \A_2}^{-1}(0)$ since 
$\nabla^2 f|_p$ is well defined (cf. Lemma \ref{wdc}, statement \ref{wd_nabla_all}). 
Showing that the section $\us_{\PP \D_4}$ is transverse to the zero is equivalent to showing that  the map 
\bgd
\hatii{\ha{\us}}_{\AA_0} \oplus \hatii{\ha{\us}}_{\AA_1} \oplus  \hatii{\ha{\us}}_{\PP \A_2} \oplus \hatii{\ha{\us}}_{\PP \D_4} : \F^* \times \C^3 \lra\C^6,
\edd
\bgd
(\qq, x,y, \eta)\mapsto \big(\qq(x,y), ~\qq_x, ~\qq_y, ~\qq_{xx} + \eta \qq_{xy}, ~\qq_{xy} + \eta \qq_{yy}, ~\qq_{yy}\big)
\edd
is transverse to zero at $(\qq, 0,0,0)$. The Jacobian matrix of this map at $(\qq,0,0,0)$ is a $6 \times (\delta_d+4)$ matrix which has full rank if $d\geq 2$; look at the first six columns which are partial derivatives with respect to $\qq_{00}$, $\qq_{10}$,  $\qq_{01}$, $\qq_{20}$, 
$\qq_{11}$ and $\qq_{02}$.  \qed \\

\ni \textbf{Proof of Corollary \ref{pa2_cl_vanish_cor}:} This follows immediately from Lemma \ref{cl3} and Proposition \ref{A2_Condition_prp}. \qed \\

\ni \textbf{Proof of Proposition \ref{D4_sharp_Condition_prp}:} Equation \eqref{d4_sharp_up_equation} is the \textit{definition} of $\DD_4^{\#}$. Towards showing transversality, note that the transversality of the sections $\us_{\AA_0}$ and $\us_{\AA_1}$ is taken care of in the proof of Proposition \ref{A1_sharp_Condition_prp}. Moreover, proving the transversality of $\us_{\DD_4}$ is identical  to the proof of transversality of $\ds_{\D_4}$ in Proposition \ref{mll}. We will now show the transversality of the section 
$\us_{\PP \A_3}$ restricted to $\us_{\DD_4}^{-1}(0)$. Let  $(\ff, \lp) \in \us_{\DD_4}^{-1}(0)$. Then it is easy to see that $\us_{\PP \A_2}(\ff, \lp) =0$. Hence, $\us_{\PP \A_3}$ is well defined (cf. Lemma \ref{wdc}, statement \ref{wd_pa3}). As before, showing that the section $\us_{\PP \A_3}$ is transverse to the zero set  is equivalent to showing that the map 
\bgd
\hatii{\ha{\us}}_{\AA_0} \oplus \hatii{\ha{\us}}_{\AA_1} \oplus \hatii{\ha{\us}}_{\DD_4} \oplus \hatii{\ha{\us}}_{\PP \A_3}  : \F^*\times \C^{3}  \lra\C^7, \,\,\, (\qq, x,y, \eta)\mapsto \big(\qq(x,y), ~\qq_x, ~\qq_y, ~\qq_{xx}, ~\qq_{xy}, ~\qq_{yy}, ~\qq_{\xx \xx \xx}\big),
\edd
is transverse to zero and $\xx = x + \eta y$. The Jacobian matrix of this map at $(\qq,0,0,0)$ is a $7 \times (\delta_d+4)$ matrix which has full rank if $d\geq 3$; choose the first seven columns as the partial derivatives with respect to $\qq_{00}$, $\qq_{10}$,  $\qq_{01}$, $\qq_{20}$, 
$\qq_{11}$, $\qq_{02}$ and $\qq_{30}$. \qed \\

\ni \textbf{Proof of Corollary \ref{d4_up_sharp_cl_vanish_cor}:} This follows immediately from Lemma \ref{cl1} and Proposition \ref{D4_sharp_Condition_prp}. \qed \\       
 
\ni \textbf{Proof of Proposition \ref{D4_Condition_prp}:} Proving \eqref{pd4_equation} requires some care. 
First let us define 
\[ \beta := f_{30}^2 f_{03}^2-6 f_{03} f_{12} f_{21} f_{30} + 4 f_{12}^3 f_{30} + 4 f_{03} f_{21}^3 - 3 f_{12}^2 f_{21}^2. \]
This is the discriminant of the cubic term in the Taylor expansion of $f$. 
Since $(\ff, \lp) \in \PP \D_4$, we conclude  
using Corollary \ref{D4_node_condition_cor}
that 
\[ \beta \neq 0, ~~f_{30} =0 ~~\implies f_{21} \neq 0 ~~\textnormal{and} ~~3f_{12}^2 - 4 f_{21} f_{03} \neq 0.   \]
Hence, 
\begin{align}
\PP \D_4 & = \{ (\ff,\lp) \in \D \times \P T\P^2: \us_{\AA_0}(\ff, \lp) =0, \us_{\AA_1}(\ff, \lp) =0, \us_{\DD_4}(\ff, \lp) =0, \nonumber\\
& \qquad \us_{\PP \A_3}(\ff, \lp) =0, ~~\us_{\PP \D_5}^{\UL}(\ff, \lp) \neq 0, ~~\us_{\PP \D_5^{\vee}}^{\UL}(\ff, \lp) \neq 0  \}. \label{pd4_def1} 
\end{align}
It is evident via linear algebra that
\begin{align}
\us_{\DD_4}(\ff, \lp) & =0,  
~~\us_{\PP \D_5}^{\UL}(\ff, \lp) \neq 0, ~~\us_{\PP \D_5^{\vee}}(\ff, \lp) \neq 0, 
~~\us_{\PP \A_3}(\ff, \lp) =0 \nonumber \\
\iff ~~\us_{\PP \A_2}(\ff, \lp) & =0, ~~\us_{\PP \A_3}(\ff, \lp) =0, 
~~\us_{\PP \D_4} (\ff,\lp) =0, ~~\us_{\PP \D_5}^{\UL}(\ff, \lp) \neq 0, ~~\us_{\PP \D_5^{\vee}}(\ff, \lp) \neq 0. \label{pd4_equiv_def}
\end{align}
Using equations \eqref{pd4_def1} and \eqref{pd4_equiv_def}, we get that  
\begin{align}
\PP \D_4 = \{ ~(\ff, \lp)\in \D \times \P T\P^2 :&  
~\us_{\AA_0} (\ff, \lp) =0, ~\us_{\AA_1}(\ff, \lp) =0, 
~\us_{\PP \A_2}(\ff, \lp) =0, ~\us_{\PP \A_3}(\ff, \lp) =0,   \nonumber \\
& ~\us_{\PP \D_4} (\ff, \lp)=0, ~\us_{\PP \D_5}^{\UL}(\ff, \lp) \neq 0, 
~~\us_{\PP \D_5^{\vee}}^{\UL}(\ff, \lp) \neq 0\}.  \label{pd4_def3} 
\end{align}
Hence, using \eqref{pd4_def3}, Corollary \ref{pa2_cl_vanish_cor}, \ref{a1_sharp_up_cl_vanish_cor} and \ref{a0_up_cl_vanish_cor} 
we obtain \eqref{pd4_equation}. Next we prove transversality of the bundle sections. 
We have already proved the transversality of the sections $\us_{\PP \A_3} $ in Proposition \ref{A2_Condition_prp}.  
We will now show the transversality of the sections 
$\us_{\PP \D_4}$, $\us_{\PP \D_5}^{\UL}$ and 
$\us_{\PP \D_5^{\vee}}^{\UL}$. 
Since $\nabla f|_{p} =0$ restricted to $\us_{\PP \A_3}^{-1}(0)$, 
we conclude that $\us_{\PP \D_4}$ is well defined by Lemma \ref{wdc}, statement \ref{wd_nabla_all}. 
As before, showing that the section $\us_{\PP \D_4} $ is transverse to the zero set is equivalent to showing that the map 
\bgd
\hatii{\ha{\us}}_{\AA_0} \oplus \hatii{\ha{\us}}_{\AA_1} \oplus \hatii{\ha{\us}}_{\PP \A_2} \oplus \hatii{\ha{\us}}_{\PP \A_3} \oplus \hatii{\ha{\us}}_{\PP \D_4}  : \F^* \times \C^3   \lra\C^7, 
\edd
\bgd
(\qq, x,y, \eta)\mapsto  \big(\qq(x,y), ~\qq_x, ~\qq_y, ~\qq_{xx} + \eta \qq_{xy}, ~\qq_{xy} + \eta \qq_{yy},~\qq_{\xx \xx \xx}, ~\qq_{yy}\big)
\edd
is transverse to zero at $(\qq,0,0,0)$, where $\xx = x + \eta y$. The Jacobian matrix of this map at $(\qq,0,0,0)$ is a $7 \times (\delta_d+4)$ matrix which has full rank if $d\geq 3$; take the first seven columns as the partial derivatives with respect to $f_{00}$, $f_{10}$,  $f_{01}$, $f_{20}$, $f_{11}$, $f_{30}$ and $f_{02}$. \\
\hf Next let us show that $\us_{\PP \D_5}^{\UL}$ is transverse. Note that restricted to $\us_{\PP \D_4}^{-1}(0)$, the section $\us_{\PP \D_5}^{\UL}$ is well defined; when restricted 
to $\us_{\PP \D_4}^{-1}(0)$, we infer that $ f(p) =0$, $\nabla f|_{p} =0 $ and $\nabla^2 f|_{p} =0$. Hence, by Lemma \ref{wdc}, statement \ref{wd_nabla_all}, 
the quantity  $\nabla^3 f|_p $ is well defined. Consequently, the section $\us_{\PP \D_5}^{\UL}$ is also well defined. 
Showing that the section $\us_{\PP \D_5}^{\UL}$ is transverse to the zero set is equivalent to showing that the map 
\bgd
\hatii{\ha{\us}}_{\AA_0} \oplus \hatii{\ha{\us}}_{\AA_1} \oplus \hatii{\ha{\us}}_{\PP \A_2} \oplus \hatii{\ha{\us}}_{\PP \A_3} \oplus \hatii{\ha{\us}}_{\PP \D_4} \oplus \hatii{\ha{\us}}_{\PP \D_5}^{\UL}: 
\F^* \times \C^3  \lra\C^8, 
\edd
\bgd
(\qq,x,y,\eta)\mapsto \big(\qq(x,y), ~\qq_x, ~\qq_y, ~\qq_{xx} + \eta \qq_{xy}, ~\qq_{xy} + \eta \qq_{yy}, ~\qq_{\xx \xx \xx}, ~\qq_{yy}, ~\qq_{\xx \xx y}\big)
\edd
is transverse to zero at $(\qq, 0,0,0)$, where $\xx = x + \eta y$. The Jacobian matrix of this map at $(\qq,0,0,0)$ is an $8 \times (\delta_d+4)$ matrix 
which has full rank if $d\geq 3$; take the first eight columns as the partial derivatives with 
respect to $\qq_{00}$, $\qq_{10}$,  $\qq_{01}$, $\qq_{20}$, $\qq_{11}$, $\qq_{30}$ ,$\qq_{02}$ and $\qq_{21}$. \\ 
\hf Finally, let us show that $\us_{\PP \D_5^{\vee}}$ is transverse. 
Since  $\nabla^3 f|_p $ is well defined restricted to $\us_{\PP \D_4}^{-1}(0)$, 
we conclude the section $\us_{\PP \D_5^{\vee}}$ is well defined 
when restricted 
to $\us_{\PP \D_4}^{-1}(0)$.  
Showing that the section $\us_{\PP \D_5^{\vee}}$ is transverse to the zero set is equivalent to showing that the map 
\bgd
\hatii{\ha{\us}}_{\AA_0} \oplus \hatii{\ha{\us}}_{\AA_1} \oplus \hatii{\ha{\us}}_{\PP \A_2} \oplus 
\hatii{\ha{\us}}_{\PP \A_3} \oplus \hatii{\ha{\us}}_{\PP \D_4} \oplus \hatii{\ha{\us}}_{\PP \D_5^{\vee}}: 
\F^* \times \C^3  \lra\C^8, 
\edd
\bgd
(\qq,x,y,\eta)\mapsto \big(\qq(x,y), ~\qq_x, ~\qq_y, ~\qq_{xx} + \eta \qq_{xy}, ~\qq_{xy} + \eta \qq_{yy}, ~\qq_{\xx \xx \xx}, ~\qq_{yy}, 
~3 \qq_{\xx y y}^2 - 4 \qq_{\xx \xx y} \qq_{y y y}\big)
\edd
is transverse to zero at $(\qq, 0,0,0)$, when $\qq_{21} \neq 0$ where $\xx = x + \eta y$. The Jacobian matrix of this map at $(\qq,0,0,0)$ is an $8 \times (\delta_d+4)$ matrix 
which has full rank if $d\geq 3$; take the first eight columns as the partial derivatives with respect to 
$\qq_{00}$, $\qq_{10}$,  $\qq_{01}$, $\qq_{20}$, $\qq_{11}$, $\qq_{30}$ ,$\qq_{02}$ and $\qq_{03}$.  \qed \\

\ni \textbf{Proof of Corollary \ref{pd4_cl_vanish_cor}:} This follows immediately from Lemma \ref{cl3} 
and Proposition \ref{D4_Condition_prp}.\footnote{Take $M$ to be $\us_{\PP \A_3}^{-1}(0)$ in Lemma \ref{cl3}.} \qed \\

\ni \textbf{Proof of Proposition \ref{A3_Condition_prp}:} By Corollary \ref{Ak_node_condition_cor} (for $k=3$),  we get that 
\begin{align}
\label{pa3_explicit}
\PP \A_3 = \{ ~(\ff, \lp)\in \D \times \P T\P^2 :& ~\us_{\AA_0} (\ff, \lp) =0, ~\us_{\AA_1}(\ff, \lp) =0, ~\us_{\PP\A_2}(\ff, \lp) =0, \nonumber \\ 
                                                 & ~\us_{\PP \A_3}(\ff, \lp) =0, ~\us_{\PP \A_4}(\ff, \lp) \neq 0, ~\us_{\PP \D_4} (\ff, \lp) \neq 0 \}. 
\end{align}
Equation \eqref{pa3_equation} now follows from \eqref{pa3_explicit} and Corollary \ref{pa2_cl_vanish_cor}, \ref{a1_sharp_up_cl_vanish_cor} and \ref{a0_up_cl_vanish_cor}. \\
\hf Towards showing transversality of the bundle sections, note that the sections $\us_{\PP \A_3}$ and $\us_{\PP \D_4}$ are already transverse (cf. Proposition \ref{A2_Condition_prp} and \ref{D4_Condition_prp}). We now prove transversality of $\us_{\PP \A_4}$. By Lemma \ref{wdc}, statement \ref{wd_pak} we conclude that  $\us_{\PP \A_4}$ is well defined. 
Showing that this section is transverse to the zero set is equivalent to showing that the map 
\bgd
\hatii{\ha{\us}}_{\AA_0} \oplus \hatii{\ha{\us}}_{\AA_1} \oplus \hatii{\ha{\us}}_{\PP \A_2} \oplus \hatii{\ha{\us}}_{\PP \A_3} \oplus \hatii{\ha{\us}}_{\PP \A_4}: \F^* \times \C^3  \lra\C^7,
\edd
\bge
(\qq,x,y,\eta)\mapsto \Big(\qq(x,y), ~\qq_x, ~\qq_y, ~\qq_{xx} + \eta \qq_{xy}, ~\qq_{xy} + \eta \qq_{yy}, ~\qq_{\xx \xx \xx}, ~\qq_{yy}\A^{\qq(\xx,y)}_4\Big)  \label{pa3_transverse_local}
\ede
is transverse to zero at $(\qq, 0,0,0)$, where $\xx = x + \eta y$. From \eqref{Ak_sections} it follows that
\bgd
\qq_{yy} \A^{\qq(\xx,y)}_4 = \qq_{yy} \qq_{\xx \xx \xx \xx} - 3 \qq_{\xx \xx y}^2. 
\edd
Hence, the Jacobian matrix of this map at $(\qq,0,0,0)$ is a $7 \times (\delta_d+4)$ matrix which has full rank if $d\geq 4$ and $\qq_{02} \neq 0$; take the first seven columns as the partial derivatives with respect to $\qq_{00}$, $\qq_{10}$,  $\qq_{01}$, $\qq_{20}$, $\qq_{11}$, $\qq_{30}$ and $\qq_{40}$. Notice that the condition $\us_{\PP \D_4}(\ff, \lp) \neq 0$ (which is equivalent to $\qq_{02} \neq 0$) is necessary to conclude that the Jacobian matrix has full rank.  \qed \\

\ni \textbf{Proof of Corollary \ref{pa3_cl_vanish_cor}:} This follows immediately from Lemma \ref{cl3} and Proposition \ref{A3_Condition_prp}. \qed \\

\ni \textbf{Proof of Corollary \ref{a3_is_smooth}:} We observe that $\A_3 = \pi(\PP \A_3)$, where $\pi: \D \times \P T\P^2 \lra \D \times \P^2$ 
is the projection map. Let us continue with the setup of Proposition \ref{A3_Condition_prp} and consider \eqref{pa3_transverse_local}. It suffices to show that the zero set of the map 
\bgd
\F^* \times \C^2 \lra \C^5, \qquad (\qq, x,y)\mapsto \bigg(\qq(x,y), ~\qq_x, ~\qq_y, ~\qq_{xx} \qq_{yy} - \qq_{xy}^2, ~\Big ( \partial_x -\frac{\qq_{xy}}{\qq_{yy}}  \partial_y \Big )^3 \qq(x,y)\bigg)
\edd
is smooth submanifold of $\F^* \times \C^2$ at $(\qq,0,0)$. The Jacobian of this map at $(\qq,0,0)$ is a $5 \times (\delta_{d}+3)$ matrix which has full rank if $d\geq 3$; take the first five columns as the partial derivatives with respect to $\qq_{00}$, $\qq_{10}$,  $\qq_{01}$, $\qq_{20}$ and $\qq_{30}$. \qed \\

\ni \textbf{Proof of Proposition \ref{Ak_Condition_prp}:} By Corollary \ref{Ak_node_condition_cor},  we get that 
\begin{align}
\label{pak_explicit}
\PP \A_k = \{ ~(\ff, \lp)\in \D \times \P T\P^2 :&  
~\us_{\AA_0} (\ff, \lp) =0, ~\us_{\AA_1}(\ff, \lp) =0, 
~\us_{\PP\A_2}(\ff, \lp) =0, \ldots, \us_{\PP \A_k}(\ff, \lp) =0, \nonumber \\ 
& ~\us_{\PP \A_{k+1}}(\ff,\lp) \neq 0, ~~\us_{\PP \D_4} (\ff,\lp) \neq 0\}. 
\end{align}
Equation \eqref{pak_equation} now follows from \eqref{pak_explicit} and Corollary \ref{pa2_cl_vanish_cor}, \ref{a1_sharp_up_cl_vanish_cor} and \ref{a0_up_cl_vanish_cor}. \\
\hf Towards proving transversality of the bundle section $\us_{\PP \A_i}$, note that by Lemma \ref{wdc} and statement \ref{wd_pak} this section is well defined. Showing that the section $\us_{\PP \A_i}$ is transverse to the zero set is equivalent to showing that the map  
\bgd
\hatii{\ha{\us}}_{\AA_0} \oplus \hatii{\ha{\us}}_{\AA_1} \oplus \hatii{\ha{\us}}_{\PP \A_2} \oplus \hatii{\ha{\us}}_{\PP \A_3} \oplus \ldots  \hatii{\ha{\us}}_{\PP \A_i}: \F^* \times \C^3 \lra\C^{i+3}, 
\edd
\bge
(\qq,x,y,\eta)\mapsto \Big(\qq, \qq_x, \qq_y, \qq_{xx} + \eta \qq_{xy}, \qq_{xy} + \eta \qq_{yy}, \qq_{\xx \xx \xx}, \qq_{yy}\A^{\qq(\xx,y)}_4, \qq_{yy}^2 \A^{\qq(\xx,y)}_5, \ldots, \qq_{yy}^{i-3}\A^{\qq(\xx,y)}_i\Big) \label{pak_transverse_local}
\ede
is transverse to zero at $(\qq, 0,0,0)$, where $\xx = x + \eta y$. Recall that $\A^{\qq}_i$ are defined in \eqref{Ak_node_conditionn}  (implicitly) and in \eqref{Ak_sections} (explicitly till $i=7$). Hence, the Jacobian matrix of this map at $(\qq,0,0,0)$ is an $ (i+3) \times (\delta_d+4)$ matrix which has full rank if $d\geq i$ and $\qq_{02} \neq 0$; take the first $(i+3)$ columns of the matrix as partial derivatives with respect to $\qq_{00}$, $\qq_{10}$,  $\qq_{01}$, $\qq_{20}$, $\qq_{11}$, $\qq_{30}, \qq_{40} \ldots, \qq_{i0}$. Notice that the condition $\us_{\PP \D_4}(\ff, \lp) \neq 0$ (which is equivalent to $\qq_{02} \neq 0$) is necessary 
to conclude that the Jacobian matrix has full rank. \qed \\

\ni \textbf{Proof of Corollary \ref{ak_is_smooth}:} This is similar to the proof of Corollary \ref{a3_is_smooth}. We observe that $\A_k = \pi(\PP \A_k)$, where $\pi: \D \times \P T\P^2 \lra \D \times \P^2$ is the projection map. It suffices to show that the zero set of the map $
\F^* \times \C^2 \lra \C^{k+2}$ given by
\bgd
(\qq, x,y)\mapsto \Big(\qq(x,y), ~\qq_x, ~\qq_y, ~\qq_{xx} \qq_{yy} - \qq_{xy}^2, ~\qq_{\xx \xx \xx}, ~~\qq_{yy}\A^{\qq(\xx,y)}_4, ~\qq_{yy}^2 \A^{\qq(\xx,y)}_5, \ldots, ~\qq_{yy}^{k-3}\A^{\qq(\xx,y)}_k\Big),
\edd
where $\xx = x -\frac{\qq_{xy}}{\qq_{yy}} y$, is smooth submanifold of $\F^* \times \C^{k+2}$ at $(\qq,0,0)$. The Jacobian of this map at $(\qq,0,0)$ is a $(k+2) \times (\delta_{d}+3)$ matrix which has full rank if $d\geq k$; take the first $(k+2)$ 
columns as the partial derivatives with respect to $\qq_{00}$, $\qq_{10}$,  $\qq_{01}$, $\qq_{20}$, $\qq_{30}, \ldots, \qq_{k0}$. Notice that we require $\qq_{02} \neq 0 $ 
to conclude that the matrix has full rank. \qed \\

\ni \textbf{Proof of Proposition \ref{D5_Condition_prp}:} Proving \eqref{pd5_equation} requires some care. 
Unravelling the definition of $\PP \D_5$ and using Corollary \ref{Dk_node_condition_cor} for $k=5$ we see that
\begin{align}
\PP \D_5 &:= \{ (\ff,\lp) \in \D \times \P T\P^2: \us_{\AA_0}(\ff, \lp) =0, ~~\us_{\AA_1}(\ff, \lp) =0, ~~\us_{\DD_4}(\ff, \lp) =0, \nonumber \\ 
& \qquad \qquad \qquad\qquad\qquad \us_{\PP \D_5}^{\UV}(\ff, \lp) = 0, ~~\us_{\PP \D_6}(\ff, \lp) \neq 0, ~~\us_{\PP \E_6}(\ff, \lp) \neq 0 \}.   \label{pd5_def1}
\end{align}
Standard linear algebra implies that
 \begin{align}
\us_{\DD_4}(\ff, \lp) & =0, ~~\us_{\PP \D_5}^{\UV}(\ff, \lp) = 0   \nonumber \\
\iff ~\us_{\PP \A_2}(\ff, \lp) &=0, ~~\us_{\PP \A_3}(\ff, \lp) =0,  ~~\us_{\PP \D_4} (\ff, \lp)=0, ~~\us_{\PP \D_5^{\UL}} (\ff, \lp) =0.  \label{pd6_equiv_def}
\end{align}
Hence, using \eqref{pd5_def1} and \eqref{pd6_equiv_def}, we get that  
\begin{align}
\PP \D_5 &= \{(\ff, \lp)\in \D \times \P T\P^2:  
\us_{\AA_0} (\ff, \lp) =0, \us_{\AA_1}(\ff, \lp) =0, 
~\us_{\PP \A_2}(\ff, \lp) =0, ~\us_{\PP \A_3}(\ff, \lp) =0,  \nonumber \\
& \qquad \qquad \qquad \us_{\PP \D_4} (\ff,\lp)=0, \us_{\PP \D_5}^{\UL}(\ff,\lp) = 0, 
\us_{\PP \D_6}(\ff,\lp) \neq 0, \us_{\PP \E_6}(\ff,\lp) \neq 0 \}.   \label{pd5_def3}
\end{align}
Equation \eqref{pd5_equation} now follows from \eqref{pd5_def3}, Corollary \ref{pd4_cl_vanish_cor}, \ref{pa2_cl_vanish_cor}, \ref{a1_sharp_up_cl_vanish_cor} and \ref{a0_up_cl_vanish_cor}. \\ 
\hf We have already proved the transversality of the section $\us_{\PP \D_5}^{\UL}$  in Proposition \ref{D4_Condition_prp}. 
We will now show the transverslaity of the sections 
$\us_{\PP \D_6}$ and  $\us_{\PP \E_6}$. By Lemma \ref{wdc}, statement \ref{wd_pd6} and \ref{wd_pe6}, these two sections are well defined. Showing that the section $\us_{\PP \D_6} $ is transverse to the zero set is equivalent to showing that the map 
\bgd
\hatii{\ha{\us}}_{\AA_0} \oplus \hatii{\ha{\us}}_{\AA_1} \oplus \hatii{\ha{\us}}_{\PP \A_2} \oplus \hatii{\ha{\us}}_{\PP \A_3} \oplus \hatii{\ha{\us}}_{\PP \D_4} \oplus \hatii{\ha{\us}}_{\PP \D_5}^{\UL} \oplus \hatii{\ha{\us}}_{\PP \D_6}: \F^*\times \C^3 \lra\C^9, 
\edd
\bgd
(\qq, x, y, \eta)\mapsto \big(\qq(x,y), ~\qq_x, ~\qq_y, ~\qq_{xx} + \eta \qq_{xy}, ~\qq_{xy} + \eta \qq_{yy}, ~\qq_{\xx \xx \xx}, ~f_{yy}, ~f_{\xx \xx y}, ~f_{\xx \xx \xx \xx} \big)
\edd
is transverse to zero at $(\qq, 0,0,0)$, where $\xx = x + \eta y$. The Jacobian matrix of this map at $(\qq,0,0,0)$ is a $ 9 \times (\delta_d+4)$ matrix which has full rank if $d\geq 4$; take the first $9$ columns as the partial derivatives with respect to $\qq_{00}$, $\qq_{10}$,  $\qq_{01}$, $\qq_{20}$, $\qq_{11}$, $\qq_{30}$, $\qq_{02}$, $\qq_{21}$ and $\qq_{40}$.\\
\hf Showing that the section $\us_{\PP \E_6}$ is transverse to the zero set is equivalent to showing that the map 
\bgd
\hatii{\ha{\us}}_{\AA_0} \oplus \hatii{\ha{\us}}_{\AA_1} \oplus \hatii{\ha{\us}}_{\PP \A_2} \oplus \hatii{\ha{\us}}_{\PP \A_3} \oplus \hatii{\ha{\us}}_{\PP \D_4} \oplus \hatii{\ha{\us}}_{\PP \D_5}^{\UL} 
\oplus \hatii{\ha{\us}}_{\PP \E_6}: \F^* \times \C^3 \lra\C^9,
\edd
\bgd
(\qq, x, y, \eta)\mapsto \big(\qq(x,y), ~\qq_x, ~\qq_y, ~\qq_{xx} + \eta \qq_{xy}, ~\qq_{xy} + \eta \qq_{yy}, ~\qq_{\xx \xx \xx}, ~\qq_{yy}, ~\qq_{\xx \xx y}, ~\qq_{\xx yy}\big)
\edd
is transverse to zero at $(\qq, 0,0,0)$, where $\xx = x + \eta y$. The Jacobian matrix of this map at $(\qq,0,0,0)$ is a $ 9 \times (\delta_d+4)$ matrix which has full rank if $d\geq 4$; take the first $9$ columns as the partial derivatives with respect to 
$\qq_{00}$, $\qq_{10}$,  $\qq_{01}$, $\qq_{20}$, $\qq_{11}$, $\qq_{30}$, $\qq_{02}$, $\qq_{21}$ and $\qq_{12}$. \qed \\

\ni \textbf{Proof of Corollary \ref{pd5_cl_vanish_cor}:} This follows immediately from Lemma \ref{cl3} and Proposition \ref{D5_Condition_prp}. \qed \\

\ni \textbf{Proof of Proposition \ref{D5_Condition_prp_dual}:} 
Proving \eqref{pd5_equation_dual} requires some care. 
Since $(\ff, \p, \lp) \in \Delta \PP \D_5^{\vee}$, we claim 
\begin{align}
f_{21} & \neq 0, \qquad \beta_1 := \frac{(f_{12} \partial_x - 2 f_{21} \partial_y)^3 f}{2 f_{21}^2} = 3 f_{12}^2 -4f_{21} f_{03} =0  
\qquad \textnormal{and} \label{pd5_dual_nv1} \\
\beta_2&:= (f_{12} \partial_x - 2 f_{21} \partial_y)^4 f = 
f_{12}^4 f_{40} - 8f_{12}^3 f_{21} f_{31} + 24 f_{12}^2 f_{21}^2 f_{22}-32 f_{12} f_{21}^3 f_{13} + 16 f_{21}^4 f_{04} \neq 0. 
\label{pd5_dual_nv_2}
\end{align}
It is easy to see that \eqref{pd5_dual_nv1} and \eqref{pd5_dual_nv_2} imply \eqref{pd5_equation_dual}. 
Let us now justify \eqref{pd5_dual_nv1} and \eqref{pd5_dual_nv_2}.
Since $(\ff, \p) \in \D_5$ there exists a non zero vector $u = m_1v + m_2w$ such that 
\begin{align}
\nabla^3 f|_{\p} (u,u, v) & = m_1^2 f_{30} + 2 m_1 m_2 f_{21} + m_2^2 f_{12} =0,  \label{nabla_cube_f_1} \\
\nabla^3 f|_{\p} (u,u, w) &= m_1^2 f_{21} + 2 m_1 m_2 f_{12} + m_2^2 f_{03} =0, \label{nabla_cube_f_2} \\
\nabla^4 f|_{\p} (u,u,u,u) & \neq 0. \label{nabla_fourth_nv}
\end{align}
Since $(\ff, \lp) \in \PP \D_5^{\vee}$, we conclude \textit{by definition} that 
\begin{align}
f_{30} &=0, \qquad f_{21} \neq 0, \qquad m_2 \neq 0. \label{m_pd5_dual} 
\end{align}
(If $m_2 =0$ then $f_{21}$ would be zero). 
Equations \eqref{m_pd5_dual} and \eqref{nabla_cube_f_1} now imply that 
\begin{align}
\frac{m_1}{m_2} &= -\frac{f_{12}}{2 f_{21}}. \label{m1_m2_value}
\end{align}
Equation \eqref{m1_m2_value} and \eqref{nabla_cube_f_2} implies \eqref{pd5_dual_nv1}. 
Finally, \eqref{nabla_fourth_nv} implies \eqref{pd5_dual_nv_2}. \\ 
\hf We have already shown $\us_{\PP \D_5}^{\UL}$ and 
$\us_{\PP \D^{\vee}_5}$ 
are  transverse to the zero set in Proposition \ref{D4_Condition_prp}. 
We will now show the transverslaity of the section 
$\us_{\PP \D_6^{\vee}}$. 
The reason why $\us_{\PP \D_6^{\vee}}$ is well defined is similar to why 
$\us_{\PP \D_6}$ is well defined. 
Showing that the section $\us_{\PP \D_6^{\vee}}$ 
is transverse to the zero set is equivalent to showing that the map 
\bgd
\hatii{\ha{\us}}_{\AA_0} \oplus \hatii{\ha{\us}}_{\AA_1} \oplus \hatii{\ha{\us}}_{\PP \A_2} \oplus \hatii{\ha{\us}}_{\PP \A_3} \oplus 
\hatii{\ha{\us}}_{\PP \D_4} \oplus \hatii{\ha{\us}}_{\PP \D_5^{\vee}} \oplus \hatii{\ha{\us}}_{\PP \D_6^{\vee}}: \F^*\times \C^3 \lra\C^9, 
\edd
\bgd
(\qq, x, y, \eta)\mapsto \big(\qq(x,y), ~\qq_x, ~\qq_y, ~\qq_{xx} + \eta \qq_{xy}, ~\qq_{xy} + \eta \qq_{yy}, ~\qq_{\xx \xx \xx}, ~f_{yy}, 
~\beta_1(\xx), ~\beta_2(\xx) \big)
\edd
is transverse to zero at $(\qq, 0,0,0)$ when $f_{21} \neq 0$, 
where $\xx = x + \eta y$ and $\beta_i(\xx)$ is $\beta_i$, but with partial with respect to $x$ replaced by  
partial with respect to $\xx$. 
The Jacobian matrix of this map at $(\qq,0,0,0)$ is 
a $ 9 \times (\delta_d+4)$ matrix which has full rank if $d\geq 4$; take the first $9$ columns as the partial derivatives with respect to 
$\qq_{00}$, $\qq_{10}$,  $\qq_{01}$, $\qq_{20}$, $\qq_{11}$, $\qq_{30}$, $\qq_{02}$, $\qq_{03}$ and $\qq_{04}$. Note that $f_{21} \neq 0$ 
is crucial here. \qed \\

\ni \textbf{Proof of Corollary \ref{pd5_cl_vanish_cor_dual}:} This follows immediately from Lemma \ref{cl3} and Proposition \ref{D5_Condition_prp_dual}. \qed \\

\ni \textbf{Proof of Corollary \ref{d5_is_smooth}:} This basically follows from the setup of Proposition \ref{D5_Condition_prp}. 
The proof is similar to the proof of Corollary \ref{a3_is_smooth} and \ref{ak_is_smooth}.    \qed \\

\ni \textbf{Proof of Proposition \ref{E6_Condition_prp}:} Proving \eqref{pe6_equation} requires some care. Unravelling the definition of $\PP \E_6$ and using Corollary \ref{E6_node_condition_cor} we conclude that
\begin{align}
\PP \E_6 &= \{ (\ff,\lp) \in \D \times \P T\P^2: \us_{\AA_0}(\ff, \lp) =0, ~~\us_{\AA_1}(\ff, \lp) =0, ~~\us_{\DD_4}(\ff, \lp) =0, ~~\us_{\PP \D_5}^{\UV}(\ff, \lp) = 0, \nonumber \\ 
& \qquad \qquad \qquad \us_{\PP \E_6}(\ff, \lp) = 0, ~~\us_{\PP \E_7}(\ff, \lp) \neq 0, ~~\us_{\PP \XC_8}(\ff, \lp) \neq 0 \}.  \label{pe6_def1}
\end{align}
Using \eqref{pe6_def1} and \eqref{pd6_equiv_def} we get that  
\begin{align}
\PP \E_6 &= \{ ~(\ff, \lp)\in \D \times \P T\P^2:  
~~\us_{\AA_0} (\ff, \lp) =0, ~\us_{\AA_1}(\ff, \lp) =0, 
~\us_{\PP \A_2}(\ff, \lp) =0, ~\us_{\PP \A_3}(\ff, \lp) =0,  \nonumber \\
& ~\us_{\PP \D_4} (\ff, \lp)=0, ~\us_{\PP \D_5}^{\UL}(\ff, \lp) = 0, 
~\us_{\PP \E_6}(\ff, \lp) = 0, ~\us_{\PP \E_7}(\ff, \lp) \neq 0, ~\us_{\PP \XC_8}(\ff, \lp) \neq 0 \}.   \label{pe6_def3}
\end{align}
\ni Equation \eqref{pe6_equation} now follows from \eqref{pe6_def3} and Corollary \ref{pd5_cl_vanish_cor}, \ref{pd4_cl_vanish_cor}, \ref{pa2_cl_vanish_cor}, \ref{a1_sharp_up_cl_vanish_cor} and \ref{a0_up_cl_vanish_cor}. \\
\hf We have already proved the transversality of the sections $\us_{\PP \E_6}$ in Proposition \ref{D5_Condition_prp}. We will now show the transverslaity of the sections $\us_{\PP \E_7}$ and $\us_{\PP \XC_8}$. Let us start with $\us_{\PP \E_7}$. By Lemma \ref{wdc}, statement \ref{wd_pe7}, the  section is well defined.   
Showing that the section $\us_{\PP \E_7}$ is transverse to the zero set is equivalent to showing that the map 
\bgd
\hatii{\ha{\us}}_{\AA_0} \oplus \hatii{\ha{\us}}_{\AA_1} \oplus \hatii{\ha{\us}}_{\PP \A_2} \oplus \hatii{\ha{\us}}_{\PP \A_3} \oplus \hatii{\ha{\us}}_{\PP \D_4} \oplus \hatii{\ha{\us}}_{\PP \D_5}^{\UL} 
\oplus \hatii{\ha{\us}}_{\PP \E_6} \oplus \hatii{\ha{\us}}_{\PP \E_7}: \F^* \times \C^3 \lra\C^{10}, 
\edd
\bgd
(\qq, x, y, \eta)\mapsto \big(\qq(x,y), ~\qq_x, ~\qq_y, ~\qq_{xx} + \eta \qq_{xy}, ~\qq_{xy} + \eta \qq_{yy}, ~\qq_{\xx \xx \xx}, ~\qq_{yy}, ~\qq_{\xx \xx y}, ~\qq_{\xx yy}, ~\qq_{\xx \xx \xx \xx}\big)
\edd
is transverse to zero at $(\qq, 0,0,0)$, where $\xx = x + \eta y$. The Jacobian matrix of this map at $(f,0,0,0)$ is a $ 10 \times (\delta_d+3)$ matrix which has full rank if $d\geq 4$; take the first $10$ columns as the partial derivatives with respect to $f_{00}$, $f_{10}$,  $f_{01}$, $f_{20}$, 
$f_{11}$, $f_{30}$, $f_{02}$, $f_{21}$, $f_{12}$ and $f_{40}$.\\
\hf Observe that restricted to $\us_{\PP\E_6}^{-1}(0)$,~$\nabla^2  f|_p =0$, whence $\us_{\PP \XC_8} $ is well defined by Lemma \ref{wdc}, statement \ref{wd_nabla_all}. Showing that this section is transverse to the zero set is equivalent to showing that the map 
\bgd
\hatii{\ha{\us}}_{\AA_0} \oplus \hatii{\ha{\us}}_{\AA_1} \oplus \hatii{\ha{\us}}_{\PP \A_2} \oplus \hatii{\ha{\us}}_{\PP \A_3} \oplus \hatii{\ha{\us}}_{\PP \D_4} \oplus \hatii{\ha{\us}}_{\PP \D_5}^{\UL} 
\oplus \hatii{\ha{\us}}_{\PP \E_6} \oplus \hatii{\ha{\us}}_{\PP \XC_8}: \F^*\times \C^3 \lra\C^{10}, 
\edd
\bgd
(\qq, x, y, \eta)\mapsto \big(\qq(x,y), ~\qq_x, ~\qq_y, ~\qq_{xx} + \eta \qq_{xy}, ~\qq_{xy} + \eta \qq_{yy}, ~\qq_{\xx \xx \xx}, ~\qq_{yy}, ~\qq_{\xx \xx y}, ~\qq_{\xx yy}, ~\qq_{yyy}\big)
\edd
is transverse to zero at $(\qq, 0,0,0)$, where $\xx = x + \eta y$. The Jacobian matrix of this map at $(\qq,0,0,0)$ is a $ 10 \times (\delta_d+4)$ matrix which has full rank if $d\geq 3$; take the first $10$ columns as the partial derivatives with respect to 
$\qq_{00}$, $\qq_{10}$,  $\qq_{01}$, $\qq_{20}$, $\qq_{11}$, $\qq_{30}$, $\qq_{02}$, $\qq_{21}$, $\qq_{12}$ and $\qq_{03}$. \qed \\

\ni \textbf{Proof of Corollary \ref{pe6_cl_vanish_cor}:} This follows immediately from Lemma \ref{cl3} and Proposition \ref{E6_Condition_prp}. \qed \\

\ni \textbf{Proof of Corollary \ref{e6_is_smooth}:} This is identical to the setup of Proposition \ref{E6_Condition_prp}. The proof is similar to the proof of Corollary \ref{a3_is_smooth} and \ref{ak_is_smooth}.    \qed \\

\ni \textbf{Proof of Proposition \ref{D6_Condition_prp}:}  
Proving \eqref{pd6_equation} requires some care. We begin by observing that Corollary \ref{Dk_node_condition_cor} (for $k=6$) implies that
\begin{align}
\label{pd6_explicit}
\PP \D_6 = \{ ~(\ff, \lp)\in \D \times \P T\P^2 :& ~\us_{\AA_0} (\ff, \lp) =0, ~\us_{\AA_1}(\ff, \lp) =0, ~\us_{\DD_4}(\ff, \lp) =0, ~\us_{\PP \D_5}^{\UV}(\ff, \lp) =0, \nonumber \\ 
                                                 & ~\us_{\PP \D_6}(\ff, \lp) = 0, ~\us_{\PP \D_7} (\ff, \lp) \neq 0, ~\us_{\PP \E_6} (\ff, \lp) \neq 0 \}. 
\end{align}
Using \eqref{pd6_explicit} and \eqref{pd6_equiv_def} we conclude that 
\begin{align}
\label{pd6_explicit_also}
\PP \D_6 & = \{ ~(\ff, \lp)\in \D \times \P T\P^2 : ~\us_{\AA_0} (\ff, \lp) =0, ~\us_{\AA_1}(\ff, \lp) =0, ~\us_{\PP \A_2}(\ff, \lp) =0, ~\us_{\PP \A_3}(\ff, \lp) =0, \nonumber \\ 
                                                 & ~\us_{\PP \D_4}(\ff, \lp) = 0, ~\us_{\PP \D_5}^{\UL} (\ff, \lp) =0, ~\us_{\PP \D_6} (\ff, \lp) = 0,  
                                                  ~\us_{\PP \D_7} (\ff, \lp) \neq 0, ~\us_{\PP \E_6} (\ff, \lp) \neq 0 \}. 
\end{align}
Equation \eqref{pd6_equation} now follows from \eqref{pd6_explicit_also} and Corollary \ref{pd5_cl_vanish_cor}, \ref{pd4_cl_vanish_cor}, \ref{pa2_cl_vanish_cor}, \ref{a1_sharp_up_cl_vanish_cor} and \ref{a0_up_cl_vanish_cor}. \\
\hf Towards proving transversality, note that we have already shown the transversality of $\us_{\PP \D_6}$ in Proposition \ref{D5_Condition_prp}. We now prove transversality of $\us_{\PP \E_6}$ and $\us_{\PP  \D_7}$. Let us start with $\us_{\PP \E_6}$. Note that restricted to $\us_{\PP \D_6}^{-1}(0)$, the quantity $\nabla ^2 f|_{p}$ vanishes, whence the section $\us_{\PP \E_6}$ is well defined (cf. Lemma \ref{wdc}, statement \ref{wd_nabla_all}). Showing that this section is transverse to the zero set is equivalent to showing that the map 
\bgd
\hatii{\ha{\us}}_{\AA_0} \oplus \hatii{\ha{\us}}_{\AA_1} \oplus \hatii{\ha{\us}}_{\PP \A_2} \oplus \hatii{\ha{\us}}_{\PP \A_3} \oplus \hatii{\ha{\us}}_{\PP \D_4} 
\oplus \hatii{\ha{\us}}_{\PP \D_5}^{\UL} 
\oplus \hatii{\ha{\us}}_{\PP \D_6} \oplus \hatii{\ha{\us}}_{\PP \E_6} : \F^* \times \C^3  \lra\C^{10},
\edd
\bgd
(\qq,x,y,\eta)\mapsto  \big(\qq(x,y), ~\qq_x, ~\qq_y, ~\qq_{xx} + \eta \qq_{xy} , ~\qq_{xy} + \eta \qq_{yy}, ~\qq_{ \xx \xx \xx}, ~\qq_{yy}, ~\qq_{\xx \xx y},  ~\qq_{\xx \xx \xx \xx}, ~\qq_{\xx y y}\big)
\edd
is transverse to zero at $(\qq, 0,0,0)$, where $\xx = x + \eta y$. The Jacobian matrix of this map at $(\qq,0,0,0)$ is a $10 \times (\delta_d+4)$ matrix which has full rank if $d\geq 4$ and 
$\qq_{02} \neq 0$. This follows by taking the first ten columns to be partial derivatives with respect to $\qq_{00}$, $\qq_{10}$,  $\qq_{01}$, $\qq_{20}$, 
$\qq_{11}$, $\qq_{30}$, $\qq_{02}$, $\qq_{21}$, $\qq_{40}$ and $\qq_{12}$. \\
\hf Observe that $\us_{\PP \D_7}$ is well defined by Lemma \ref{wdc}, statement \ref{wd_pdk}. 
Showing that this section is transverse to the zero set  is equivalent to showing that the map 
\bgd
\hatii{\ha{\us}}_{\AA_0} \oplus \hatii{\ha{\us}}_{\AA_1} \oplus \hatii{\ha{\us}}_{\PP \A_2} \oplus \hatii{\ha{\us}}_{\PP \A_3} 
\oplus \hatii{\ha{\us}}_{\PP \D_4} \oplus \hatii{\ha{\us}}_{\PP \D_5}^{\UL} 
\oplus \hatii{\ha{\us}}_{\PP \D_6} \oplus \hatii{\ha{\us}}_{\PP \D_7} : \F^* \times \C^3  \lra\C^{10},
\edd
\bgd
(\qq,x,y,\eta) \mapsto \Big(\qq(x,y), ~\qq_x, ~\qq_y, ~\qq_{xx} + \eta \qq_{xy} , ~\qq_{xy} + \eta \qq_{yy}, ~\qq_{\xx \xx \xx}, ~\qq_{yy}, ~\qq_{\xx \xx y},  
~\qq_{\xx \xx \xx \xx}, ~\qq_{\xx yy} \D^{f(\xx,y)}_7\Big)  
\edd
is transverse to zero at $(\qq, 0,0,0)$, where $\xx = x + \eta y$. From \eqref{Formula_Dk} we see that 
\bgd
\qq_{\xx yy} \D^{\qq(\xx,y)}_7 = \qq_{\xx yy}\qq_{\xx \xx \xx \xx \xx} - \frac{5 \qq_{\xx \xx \xx y}^2}{3}.
\edd
Hence, the Jacobian matrix of this map at $(\qq,0,0,0)$ is a $10 \times (\delta_d+4)$ matrix which has full rank if $d\geq 5$ and $\qq_{12} \neq 0$. 
This is evident by taking the first ten columns to be partial derivatives with respect to 
$\qq_{00}$, $\qq_{10}$,  $\qq_{01}$, $\qq_{20}$, 
$\qq_{11}$, $\qq_{30}$, $\qq_{02}$, $\qq_{21}$, $\qq_{40}$ and $\qq_{50}$. 
Notice that the condition $\qq_{12}\neq 0$ is necessary to conclude that the matrix has full rank. \qed \\

\ni \textbf{Proof of Corollary \ref{pd6_cl_vanish_cor}:} This follows immediately from Lemma \ref{cl3} and Proposition \ref{D6_Condition_prp}. \qed \\

\ni \textbf{Proof of Corollary \ref{d6_is_smooth}:} This basically follows from the setup of Proposition \ref{D6_Condition_prp}. 
The proof is similar to the proof of Corollary \ref{a3_is_smooth} and \ref{ak_is_smooth}.   \qed  \\

\ni \textbf{Proof of Proposition \ref{Dk_Condition_prp}:} Proving \eqref{pdk_equation} requires some care. Notice that Corollary \ref{Dk_node_condition_cor} implies that
\begin{align}
\label{pdk_explicit}
\PP \D_k & = \{ ~(\ff, \lp)\in \D \times \P T\P^2 : ~\us_{\AA_0} (\ff, \lp) =0, ~\us_{\AA_1}(\ff, \lp) =0, ~\us_{\DD_4}(\ff, \lp) =0, ~\us_{\PP \D_5}^{\UV}(\ff, \lp) =0, \nonumber \\ 
                                                 &\us_{\PP \D_6}(\ff, \lp) = 0, \ldots, ~\us_{\PP \D_k} (\ff, \lp) = 0, ~\us_{\PP \D_{k+1}} (\ff, \lp) \neq 0,  ~\us_{\PP \E_6} (\ff, \lp) \neq 0 \}. 
\end{align}
Equation \eqref{pdk_equation} now follows from \eqref{pdk_explicit}, \eqref{pd6_equiv_def} and Corollary \ref{pd6_cl_vanish_cor}, \ref{pd5_cl_vanish_cor}, \ref{pd4_cl_vanish_cor}, \ref{pa2_cl_vanish_cor}, \ref{a1_sharp_up_cl_vanish_cor} and \ref{a0_up_cl_vanish_cor}. \\
\hf We now prove transversality of the bundle sections $\us_{\PP \D_i}$. By Lemma \ref{wdc}, statement \ref{wd_pdk} these sections are well defined. Showing that these sections are transverse to the zero set is equivalent to showing that the map 
\bgd
\hatii{\ha{\us}}_{\AA_0} \oplus \hatii{\ha{\us}}_{\AA_1} \oplus \hatii{\ha{\us}}_{\PP \A_2} \oplus \hatii{\ha{\us}}_{\PP \A_3} \oplus \hatii{\ha{\us}}_{\PP \D_4} \oplus \hatii{\ha{\us}}_{\PP \D_5^{\UL}} 
\oplus \hatii{\ha{\us}}_{\PP \D_6}  \oplus \ldots  \hatii{\ha{\us}}_{\PP \D_i}   : \F^* \times \C^3 \lra\C^{i+3},
\edd
\bgd
(\qq,x,y,\eta) \mapsto \bigg(\qq(x,y) , \qq_x, \qq_y, \qq_{xx} + \eta \qq_{xy}, \qq_{xy} + \eta \qq_{yy}, \qq_{\xx \xx \xx}, \qq_{yy}, \qq_{ \xx \xx y}, \D^{\qq(\xx,y)}_6,  \qq_{12} \D^{\qq(\xx,y)}_7,  \ldots, \qq_{12}^{i-6} \D^{\qq(\xx,y)}_i\bigg)
\edd
is transverse to zero at $(\qq, 0,0,0)$, where $\xx = x + \eta y$. Recall that $\D^f_i$ is defined in \eqref{Dk_node_conditionn} implicitly and in \eqref{Formula_Dk} explicitly till $i=8$. The Jacobian matrix of this map at $(\qq,0,0,0)$ is an $ (i+3) \times (\delta_d+4)$ matrix which has full rank if $d\geq i-2$ and $\qq_{12} \neq 0$. This follows by choosing the first $(i+3)$ columns to be the partial derivatives with respect to $\qq_{00}$, $\qq_{10}$,  $\qq_{01}$, $\qq_{20}$, $\qq_{11}$, $\qq_{30}$, $\qq_{02}$, $\qq_{21}$, $\qq_{40}$, $\qq_{50}, \ldots$ and $\qq_{i-2, 0}$. Notice that the condition $\us_{\PP \E_6}(\ff, \lp) \neq 0$ (which is equivalent to $\qq_{12} \neq 0$) is necessary to conclude that the Jacobian matrix has full rank.  \qed \\

\ni \textbf{Proof of Corollary \ref{dk_is_smooth}:} This basically follows from the setup of Proposition \ref{Dk_Condition_prp}. The proof is similar to the proof of Corollary \ref{a3_is_smooth} and \ref{ak_is_smooth}.  \qed   \\

\ni \textbf{Proof of Proposition \ref{E7_Condition_prp}:} Proving \eqref{pe7_equation} requires some care. Unravelling the definition of $\PP \E_7$ in conjunction with Corollary \ref{E7_node_condition_cor} we gather
\begin{align}
\PP \E_7 &= \{ (\ff,\lp) \in \D \times \P T\P^2: \us_{\AA_0}(\ff, \lp) =0, ~~\us_{\AA_1}(\ff, \lp) =0, ~~\us_{\DD_4}(\ff, \lp) =0, ~~\us_{\PP \D_5}^{\UV}(\ff, \lp) = 0, \nonumber \\ 
    & \us_{\PP \E_6}(\ff, \lp)=0,  ~~\us_{\PP \E_7}(\ff, \lp) = 0, ~~\us_{\PP \E_8}(\ff, \lp) \neq 0, ~~\us_{\PP \XC_8}(\ff, \lp) \neq 0 \}.  \label{pe7_def1}
\end{align}
Hence, using \eqref{pe7_def1} and \eqref{pd6_equiv_def}, we get that  
\begin{align}
\label{pe7_explicit}
\PP \E_7 &= \{ ~(\ff, l_p)\in \D \times \P T\P^2:  \us_{\AA_0} (\ff, l_p) =0, ~\us_{\AA_1}(\ff, l_p) =0, ~\us_{\PP \A_2}(\ff, l_p) =0, ~\us_{\PP \A_3}(\ff, l_p) =0,  \nonumber \\
& \qquad\qquad\qquad\qquad\qquad ~~\us_{\PP \D_4} (\ff,l_p)=0, ~\us_{\PP \D_5}^{\UL}(\ff,l_p) = 0, ~\us_{\PP \E_6}(\ff,l_p) = 0,\nonumber \\
& \qquad\qquad\qquad\qquad\qquad ~~\us_{\PP \E_7}(\ff,l_p) = 0 ~~\us_{\PP \E_8}(\ff,l_p) \neq 0, ~~\us_{\PP \XC_8}(\ff,l_p) \neq 0 \}.   
\end{align}
Equation \eqref{pe7_equation} now follows from equation \eqref{pe7_explicit} and Corollary \ref{pe6_cl_vanish_cor}, \ref{pd5_cl_vanish_cor}, \ref{pd4_cl_vanish_cor}, \ref{pa2_cl_vanish_cor}, \ref{a1_sharp_up_cl_vanish_cor} and \ref{a0_up_cl_vanish_cor}. \\ 
\hf We have already proved transversality of $\us_{\PP \E_7}$ in Proposition \ref{E6_Condition_prp}. 
We now prove transversality of $\us_{\PP \E_8}$ and $\us_{\PP  \XC_8}$. \ni Lets us start with $\us_{\PP \E_8}$. By Lemma \ref{wdc}, statement \ref{wd_pe8} the section is well defined. 
Showing that this section is transverse to the zero set, is equivalent to showing that the map 
\bgd
\hatii{\ha{\us}}_{\AA_0} \oplus \hatii{\ha{\us}}_{\AA_1} \oplus \hatii{\ha{\us}}_{\PP \A_2} \oplus \hatii{\ha{\us}}_{\PP \A_3} 
\oplus \hatii{\ha{\us}}_{\PP \D_4} \oplus \hatii{\ha{\us}}_{\PP \D_5}^{\UL} 
\oplus \hatii{\ha{\us}}_{\PP \E_6} \oplus \hatii{\ha{\us}}_{\PP \E_7} \oplus \hatii{\ha{\us}}_{\PP \E_8}  : \F^* \times \C^3  \lra\C^{11},
\edd
\bgd(\qq,x,y,\eta)\mapsto\big(\qq(x,y), ~\qq_x, ~\qq_y, ~\qq_{xx} + \eta \qq_{xy} , ~\qq_{xy} + \eta \qq_{yy}, ~\qq_{ \xx \xx \xx}, ~\qq_{yy}, ~\qq_{\xx \xx y},  ~\qq_{\xx yy}, ~\qq_{\xx \xx \xx \xx}, ~\qq_{\xx \xx \xx y}\big).  
\edd
is transverse to zero at $(\qq, 0,0,0)$, where $\xx = x + \eta y$. The Jacobian matrix of this map at $(\qq,0,0,0)$ is an $11 \times (\delta_d+4)$ matrix which has full rank if $d\geq 4$ and 
$\qq_{12} \neq 0$. This follows by taking the first eleven columns as the partial derivatives with respect to $\qq_{00}$, $\qq_{10}$,  $\qq_{01}$, $\qq_{20}$, 
$\qq_{11}$, $\qq_{30}$, $\qq_{02}$, $\qq_{21}$, $\qq_{12}$, $\qq_{40}$ and $\qq_{31}$. \\
\hf To show that $\us_{\PP \XC_8}$ is transverse to the zero set, we observe that restricted to $\us_{\PP \E_7}^{-1}(0)$, the quantity $\nabla^3 f|_{p}$ is identically zero (via linear algebra). Hence, the section is well defined by Lemma \ref{wdc}, statement \ref{wd_nabla_all}. Showing that this section is transverse to the zero set is equivalent to showing that the map  
\bgd
\hatii{\ha{\us}}_{\AA_0} \oplus \hatii{\ha{\us}}_{\AA_1} \oplus \hatii{\ha{\us}}_{\PP \A_2} \oplus \hatii{\ha{\us}}_{\PP \A_3} \oplus 
\hatii{\ha{\us}}_{\PP \D_4} \oplus \hatii{\ha{\us}}_{\PP \D_5}^{\UL} 
\oplus \hatii{\ha{\us}}_{\PP \E_6} \oplus \hatii{\ha{\us}}_{\PP \E_7} \oplus \hatii{\ha{\us}}_{\PP \XC_8}  : \F^* \times \C^3  \lra\C^{11},
\edd
\bgd
(\qq,x,y,\eta)\mapsto \big(\qq(x,y), ~\qq_x, ~\qq_y, ~\qq_{xx} + \eta \qq_{xy} , ~\qq_{xy} + \eta \qq_{yy}, ~\qq_{\xx \xx \xx}, ~\qq_{yy}, ~\qq_{\xx \xx y},  ~\qq_{\xx yy}, ~\qq_{\xx \xx \xx \xx}, ~\qq_{yyy}\big)  
\edd
is transverse to zero at $(\qq, 0,0,0)$, where $\xx = x + \eta y$. The Jacobian matrix of this map at $(\qq,0,0,0)$ is an $11 \times (\delta_d+4)$ matrix which has full rank if $d\geq 4$ and $\qq_{12} \neq 0$; take the first eleven columns as the partial derivatives with respect to 
$\qq_{00}$, $\qq_{10}$,  $\qq_{01}$, $\qq_{20}$, $\qq_{11}$, $\qq_{30}$, $\qq_{02}$, $\qq_{21}$, $\qq_{12}$, $\qq_{40}$ and $\qq_{03}$.  \qed \\

\ni \textbf{Proof of Corollary \ref{pe7_cl_vanish_cor}:} This follows immediately from Lemma \ref{cl3} and Proposition \ref{E7_Condition_prp}. \qed \\

\ni \textbf{Proof of Corollary \ref{e7_is_smooth}:} This basically follows from the setup of Proposition \ref{E7_Condition_prp}. 
The proof is similar to the proof of Corollary \ref{a3_is_smooth} and \ref{ak_is_smooth}.  \qed  \\

\ni \textbf{Proof of Proposition \ref{E7_Condition_prp_using_D6}:}  
Observe that $\us_{\PP \D_6} = \us_{\PP \E_7}$. Equation \eqref{pe7_equation_also} follows from Proposition \ref{E7_Condition_prp}, Corollary \ref{pd6_cl_vanish_cor} and \ref{pe6_cl_vanish_cor}. Transversality of $\us_{\PP \E_7}$ has been proven in Proposition \ref{D6_Condition_prp} (there it was denoted as $\us_{\PP \D_6}$). Proving that the sections $\us_{\PP \E_8}$ and $\us_{\PP \XC_8}$ are transverse to the zero set is almost identical to the proof in Proposition \ref{E7_Condition_prp}.  \qed \\ 

\ni \textbf{Proof of Corollary \ref{pe7_cl_vanish_cor_also}:} This follows immediately from Lemma \ref{cl3} and Proposition \ref{E7_Condition_prp_using_D6}. \qed \\ 

\ni \textbf{Proof of Proposition \ref{X8_sharp_Condition_prp}:} Equation \eqref{x8_sharp_and_flat_equation} is the {\it definition(s)} 
of the spaces $\XC_8^{\#}$ and $\XC_{8}^{\# \flat}$ so there is nothing to prove. We have already proved the transversality of the sections $\us_{\AA_0}$, $\us_{\AA_1}$ and $\us_{\DD_4} $ in Proposition \ref{D4_sharp_Condition_prp}.  We will now show the transversality of the sections $\us_{\hat{\XC}_8} $, $\us_{\PP \E_7}$ and $\us_{\J}$.\\ 
\hf Let us start with $\us_{\hat{\XC}_8}$. Observe that by Lemma \ref{wdc}, statement \ref{wd_nabla_all} this section is well defined. Showing that this section is transverse to the zero set is equivalent to showing that the map 
\bgd
\hatii{\ha{\us}}_{\AA_0} \oplus \hatii{\ha{\us}}_{\AA_1} \oplus \hatii{\ha{\us}}_{\DD_4} \oplus \hatii{\ha{\us}}_{\hat{\XC}_8} : \F^* \times \C^3 \lra\C^{10},
\edd
\bgd
(\qq, x, y, \eta)\mapsto \big( \qq(x,y), ~\qq_x, ~\qq_y, ~\qq_{xx}, ~\qq_{xy}, ~\qq_{yy}, ~\qq_{xxx}, ~\qq_{xxy}, ~\qq_{xyy}, ~\qq_{yyy}\big)
\edd
is transverse to zero at $(\qq, 0,0,0)$. The Jacobian matrix of this map at $(\qq,0,0,0)$ is a $10 \times (\delta_d+4)$ matrix which has full rank if $d\geq 3$; take the first ten columns as the partial derivatives with respect to $\qq_{00}$, $\qq_{10}$,  $\qq_{01}$, $\qq_{20}$, 
$\qq_{11}$, $\qq_{02}$, $\qq_{30}$, $\qq_{21}$, $\qq_{12}$ and $\qq_{03}$. \\
\hf By Lemma \ref{wdc}, statement \ref{wd_nabla_all} the section $\us_{\PP \E_7}$ is well defined. Showing that this section is transverse to the zero set is equivalent to showing that the map 
\bgd
\hatii{\ha{\us}}_{\AA_0} \oplus \hatii{\ha{\us}}_{\AA_1} \oplus \hatii{\ha{\us}}_{\DD_4} \oplus \hatii{\ha{\us}}_{\hat{\XC}_8} \oplus \hatii{\ha{\us}}_{\PP \E_7}: 
\F^* \times \C^3 \lra\C^{11}, 
\edd
\bgd
(\qq, x,y, \eta)\mapsto \big(\qq(x,y), ~\qq_x, ~\qq_y, ~\qq_{xx}, ~\qq_{xy}, ~\qq_{yy}, ~\qq_{x x x}, ~\qq_{x x y}, ~\qq_{x yy}, ~\qq_{yyy}, ~\qq_{\xx \xx \xx \xx} \big)
\edd
is transverse to zero at $(\qq, 0,0,0)$, where $\xx = x + \eta y$. The Jacobian matrix of this map at $(\qq,0,0,0)$ is an $11 \times (\delta_d+4)$ matrix which has full rank if $d\geq 4$; take the eleven columns as the partial derivatives with respect to
$\qq_{00}$, $\qq_{10}$,  $\qq_{01}$, $\qq_{20}$, $\qq_{11}$, $\qq_{02}$, $\qq_{30}$, $\qq_{21}$, $\qq_{12}$, $\qq_{03}$ and $\qq_{40}$.\\
\hf Finally, observe that by Lemma \ref{wdc}, statement \ref{wd_nabla_all} the section 
$\us_{\J}$ is well defined. Showing that this section is transverse to the zero set  
is equivalent to showing that the map 
\bgd
\hatii{\ha{\us}}_{\AA_0} \oplus \hatii{\ha{\us}}_{\AA_1} \oplus \hatii{\ha{\us}}_{\DD_4} \oplus \hatii{\ha{\us}}_{\hat{\XC}_8} \oplus \hatii{\ha{\us}}_{\J}  : 
\F^* \times \C^3 \lra\C^{11}, 
\edd
\vspace*{-0.8cm}
\begin{align*}
(\qq, x,y, \eta) \mapsto & \Big(\qq(x,y), ~\qq_x, ~\qq_y, ~\qq_{xx}, ~\qq_{xy}, ~\qq_{yy}, ~\qq_{xxx}, ~\qq_{xxy}, ~\qq_{xyy}, ~\qq_{yyy},\\ 
&\qquad \qquad -\frac{\qq_{\xx \xx \xx y}^2}{8} + \frac{3 \qq_{\xx \xx yy} \qq_{\xx \xx \xx y} \qq_{\xx \xx \xx \xx}}{16} - \frac{\qq_{ \xx yyy} \qq_{\xx \xx \xx \xx}^2}{16}  \Big)
\end{align*}
is transverse to zero at $(\qq, 0,0,0)$, where $\xx = x + \eta y$. The Jacobian matrix of this map at $(\qq,0,0,0)$ is an $11 \times (\delta_d+4)$ matrix which has full rank if $d\geq 4$ 
and $f_{40} \neq 0$. This follows by taking the first eleven columns to be the partial derivatives with respect to
$\qq_{00}$, $\qq_{10}$,  $\qq_{01}$, $\qq_{20}$, 
$\qq_{11}$, $\qq_{02}$, $\qq_{30}$, $\qq_{21}$, $\qq_{12}$, $\qq_{03}$ and $\qq_{13}$. 
Notice that the condition $\us_{\PP \E_7}(\ff, \lp) \neq 0$ (which is equivalent to $\qq_{40} \neq 0$) is necessary for the matrix to have full rank. \qed \\

\ni \textbf{Proof of Proposition \ref{X8_and_X8_sharp_up_cl_vanish_cor}:} The result for $\ov{\hat{\XC}_8^{\#}}$ follows from Lemma \ref{cl1} and Proposition \ref{X8_sharp_Condition_prp} while that for $\ov{\hat{\XC}_8^{\# \flat}}$ follows from Lemma \ref{cl3} and Proposition \ref{X8_sharp_Condition_prp}. \qed 


\section{Closure} 
\label{closure_of_spaces}

\hf\hf We compute the closure of the various spaces; this is the main result(s) of this section. 
\begin{lmm}
\label{cl}
Let $\X_k$ be a singularity of type $\A_k$, $\D_k$, 
$\E_k$ or $\XC_8$. Then the closures are given by :
\begin{enumerate} 
\item \label{A0cl} $\ov{\A}_0 = \A_0 \cup \ov{\A}_1$ 
\qquad if  $d \geq 2$.
\item \label{A1cl}
$\ov{\hat{\A}_1} = \ov{\hat{\A}^{\#}_1} = \hat{\A}_1^{\#} \cup 
\ov{\PP \A}_2$ \qquad if $d \geq 3$.
\item \label{D4_cl_no_direction}$ \ov{\hat \D^{\#}_4} = \hat \D^{\#}_4 \cup \ov{\PP \D}_4 $ \qquad if $d \geq 3$.
\item \label{D4cl}$ \ov{\PP \D}_4 = \PP \D_4 \cup \ov{\PP \D}_5 \cup \ov{\PP \D_5^{\vee}}$ 
\qquad if $d \geq 4$. 
\item \label{E6cl}$\ov{\PP \E}_6 = 
\PP \E_6 \cup \ov{\PP \E}_7 \cup \ov{\hat{\XC}^{\#}_8}$ 
\qquad if $d \geq 4$. 
\item \label{D5cl}$\ov{\PP \D}_5 = \PP \D_5 \cup \ov{\PP \D}_6 \cup \ov{\PP \E}_6$ \qquad if $d \geq 4$.
\item \label{D6cl}$\ov{\PP \D}_6 = \PP \D_6 \cup \ov{\PP \D}_7 \cup \ov{\PP \E}_7$ \qquad if $d \geq 5$.
\item \label{A2cl}$\ov{\PP \A}_2 = \PP \A_2 \cup 
\ov{\PP \A}_3 \cup \ov{\hat{\D}_4^{\#}} $ 
\qquad if $d \geq 4$.
\item \label{A3cl}$\ov{\PP \A}_3 = \PP \A_3 \cup  \ov{\PP \A}_4 \cup \ov{\PP \D}_4$ \qquad if $d \geq 5$.
\item \label{A4cl}$\ov{\PP \A}_4 = \PP \A_4 \cup  \ov{\PP \A}_5 \cup \ov{\PP \D}_5$ \qquad if $d \geq 6$.
\item \label{A5cl}$\ov{\PP \A}_5 = \PP \A_5 \cup  \ov{\PP \A}_6 \cup \ov{\PP \D}_6 \cup \ov{\PP \E}_6 $ \qquad if $d \geq 7$.
\item \label{A6cl}$\ov{\PP \A}_6 = \PP \A_6 \cup  \ov{\PP \A}_7 \cup 
\ov{\PP \D}_7 \cup \ov{\PP \E}_7 \cup \ov{\hat{\XC}_8^{\# \flat}}$  
 \qquad if $d \geq 8$.\\
\end{enumerate}
\end{lmm}

\textbf{Proof of Lemma \ref{cl} (\ref{A0cl}):} Follows from Corollary \ref{a0_cl_vanish_cor}, Proposition \ref{ift_ml},  Corollary \ref{a1_cl_vanish_cor}.  \qed \\

\textbf{Proof of Lemma \ref{cl} (\ref{A1cl}):} Follows from Corollary  \ref{a1_up_cl_vanish_cor}, \ref{a1_sharp_up_cl_vanish_cor}, Proposition \ref{A1_sharp_Condition_prp}, Corollary \ref{pa2_cl_vanish_cor}. \qed \\

\textbf{Proof of Lemma \ref{cl} (\ref{D4_cl_no_direction}):} 
Follows from Corollary \ref{d4_up_sharp_cl_vanish_cor},  Proposition \ref{D4_sharp_Condition_prp}, Corollary \ref{pd4_cl_vanish_cor}, \ref{pa2_cl_vanish_cor}, \ref{a1_sharp_up_cl_vanish_cor} and \ref{a0_up_cl_vanish_cor}. \qed \\

\textbf{Proof of Lemma \ref{cl} (\ref{D4cl}):} Follows from Corollary \ref{pd4_cl_vanish_cor}, Proposition \ref{D4_Condition_prp}, Corollary \ref{pd5_cl_vanish_cor} 
and \ref{pd5_cl_vanish_cor_dual}.  \qed \\

\textbf{Proof of Lemma \ref{cl} (\ref{E6cl}):} Corollary \ref{X8_and_X8_sharp_up_cl_vanish_cor} implies that 
\bgd
\ov{\XX_8^{\#}} = \{ (\ff,\lp) \in  \D \times \P T\P^2:  
~\us_{\AA_0}(\ff,\lp)=0, ~\us_{\AA_1}(\ff,\lp),   ~\us_{\DD_4}(\ff,\lp) =0, ~\us_{\XX_8}(\ff,\lp) =0 \}.
\edd
Standard linear algebra implies that
\begin{align*}
\ov{\XX_8^{\#}}  &= \{ (\ff,\lp) \in  \D \times \P T\P^2 :  
\us_{\AA_0}(\ff, \lp) =0, \us_{\AA_1}(\ff, \lp)=0, 
\us_{\PP \A_2}(\ff,\lp) =0, \us_{\PP \A_3}(\ff,\lp) =0,\\ 
& \qquad \qquad  \us_{\PP \D_4}(\ff,\lp) =0, \us_{\PP \D_5}^{\UL}(\ff,\lp) =0, 
\us_{\PP \E_6}(\ff,\lp) =0, ~~\us_{\PP \XC_8}(\ff,\lp)=0 \}.
\end{align*}
By Corollary \ref{pe6_cl_vanish_cor}, \ref{pd5_cl_vanish_cor}, \ref{pd4_cl_vanish_cor}, \ref{pa2_cl_vanish_cor}, \ref{a1_sharp_up_cl_vanish_cor}, \ref{a0_up_cl_vanish_cor}, we conclude that 
\begin{align}
\label{x8_in_closure_of_e6_equation}
 \ov{\XX_8^{\#}} & = \{ (\ff, \lp) \in \ov{\PP \E}_6: \us_{\PP \XC_8} (\ff, \lp) =0  \}.
\end{align}
Lemma \ref{cl}, statement \ref{E6cl} now follows from \eqref{x8_in_closure_of_e6_equation}, Corollary \ref{pe7_cl_vanish_cor}, Proposition \ref{E6_Condition_prp} and Corollary \ref{pe6_cl_vanish_cor}. \qed \\  

\textbf{Proof of Lemma \ref{cl} (\ref{D5cl}):} Follows from Corollary \ref{pd5_cl_vanish_cor}, Proposition \ref{D5_Condition_prp}, Corollary \ref{pd6_cl_vanish_cor} and \ref{pe6_cl_vanish_cor}. \qed \\

\textbf{Proof of Lemma \ref{cl} (\ref{D6cl}):} Follows from Proposition \ref{Dk_Condition_prp}. Since the section
\bgd
\us_{\PP \E_6}:\ov{\PP \D}_6 \lra \UL_{\PP \E_6}
\edd
is transverse to the zero set (as proved in Proposition \ref{D6_Condition_prp}), the hypothesis of the second part of Lemma \ref{cl2} is satisfied. By Proposition \ref{D6_Condition_prp} and Corollary \ref{pd6_cl_vanish_cor} we conclude that 
\bgd
\PP \D_6 = \{ (\ff, \lp) \in \ov{\PP \D}_6: \us_{\PP \D_7}(\ff, \lp) \neq 0, ~\us_{\PP \E_6}(\ff, \lp) \neq 0 \}.
\edd
Corollary \ref{pe7_cl_vanish_cor_also} now proves our claim.\qed \\

\textbf{Proof of Lemma \ref{cl} (\ref{A2cl}):} Follows from Corollary \ref{pa2_cl_vanish_cor}, Proposition \ref{A2_Condition_prp},  Corollary \ref{pa3_cl_vanish_cor}, 
\ref{d4_up_sharp_cl_vanish_cor}, \ref{a1_sharp_up_cl_vanish_cor} and \ref{a0_up_cl_vanish_cor}.  
\qed \\

\textbf{Proof of Lemma \ref{cl} (\ref{A3cl}):} Follows from the second part of Lemma \ref{cl2} and Proposition \ref{Ak_Condition_prp}. 
Since the section $\us_{\PP \D_4} : \ov{\PP \A}_3 \lra \UL_{\PP \D_4}$ is transverse to the zero set (cf. Proposition \ref{A3_Condition_prp}), 
the hypothesis of the second part of Lemma \ref{cl2} is satisfied. Corollary \ref{pd4_cl_vanish_cor} and \ref{pa3_cl_vanish_cor} now imply that 
\begin{align}
\label{pa3_cap_pd4_equal_pd4} 
\ov{\PP \D}_4 &= \{ (\ff,\lp) \in \ov{\PP \A}_3: \us_{\PP \D_4}(\ff,\lp) = 0 \}.
\end{align}
By Proposition \ref{A3_Condition_prp} and Corollary \ref{pa3_cl_vanish_cor} we get that 
\begin{align*}
\PP \A_3 &= \{ (\ff,\lp) \in \ov{\PP \A}_3: \us_{\PP \A_4}(\ff,\lp) \neq 0, ~~\us_{\PP \D_4}(\ff,\lp) \neq 0 \}.
\end{align*}
The second part of Lemma \ref{cl2} now proves our claim.\qed \\

\textbf{Proof of Lemma \ref{cl} (\ref{A4cl}):} By Lemma \ref{cl2} applied to
\bgd
M=\ov{\PP \A}_3,\,\,\ts_0=\us_{\PP \A_4},\,\,\ts_1=\us_{\PP \A_5},\,\,\ts_2=\us_{\PP \A_6},\,\,\nts=\us_{\PP\D_4},
\edd
and Proposition \ref{Ak_Condition_prp}, it suffices to show that
\begin{align} 
\{ (\ff,\lp) \in \ov{\PP \A}_4: \us_{\PP \D_4}(\ff,\lp) = 0 \} = \ov{\PP \D}_5. \label{pa4_cap_pd4_equal_pd5}
\end{align}
Let us show that the left hand side of \eqref{pa4_cap_pd4_equal_pd5} is a subset of the right hand side. We claim that
\bge
\ov{\PP \A}_4 \cap \PP \D_4 = \varnothing. \label{pa4closure_cap_pd4_empty}  
\ede
To see this, first we observe that if $(\ff, \lp) \in \PP \D_4 $ then $\us_{\PP \D_4}(\ff, \lp)=0$ and $\us_{\PP \D_5}^{\UL}(\ff, \lp)\neq 0$. Therefore, 
\bgd
\us_{\PP \A_4}(\ff, \lp) = \us_{\PP \D_4}(\ff, \lp) \us_{\PP \D_6}(\ff, \lp) - 3 \us_{\PP \D_5}^{\UL}(\ff, \lp)^2   = -3 \us_{\PP \D_5}^{\UL}(\ff, \lp)^2 \neq 0.
\edd
This implies that if $(\ff(t), \lp(t))$ lies in a small neighborhood of $(\ff, \lp)$ then $\us_{\PP \A_4}(\ff(t), \lp(t)) \neq 0$, proving \eqref{pa4closure_cap_pd4_empty}. By Lemma \ref{cl}, statement \ref{A3cl} we have $\ov{\PP \A}_4 \subset \ov{\PP \A}_3$. Therefore, 
\bgd
\{ (\ff, \lp) \in \ov{\PP \A}_4: \us_{\PP \D_4}(\ff,\lp) = 0 \} \subset \{ (\ff,\lp) \in \ov{\PP \A}_3: \us_{\PP \D_4}(\ff,\lp) = 0 \}.
\edd
The right hand side above equals 
\bgd
\PP \D_4 \cup \ov{\PP \D}_5
\edd
by \eqref{pa3_cap_pd4_equal_pd4} and Lemma \ref{cl}, statement \ref{D4cl}. Hence,
\bgd
\{ (\ff,\lp) \in \ov{\PP \A}_4: \us_{\PP \D_4}(\ff,\lp) = 0 \} \subset \ov{\PP \D}_5
\edd
by \eqref{pa4closure_cap_pd4_empty}. This proves that the left hand side of \eqref{pa4_cap_pd4_equal_pd5} is a subset of the right hand side. For the converse note 
since $\ov{\PP \A}_4$ is a closed set, it suffices to show that
\begin{align} 
\label{pa4_cap_pd4_supset_pd5}
\{ (\ff,l_p) \in \ov{\PP \A}_4: \us_{\PP \D_4}(\ff,l_p) = 0 \}  & \supset  \PP \D_5.  
\end{align}
We will simultaneously prove statement \eqref{pa4_cap_pd4_supset_pd5} 
and also prove the following statement
\begin{align}
\ov{\PP \A}_5 \cap \PP \D_5 &= \varnothing . \label{pa5closure_cap_pd5_empty}  
\end{align}
Since $\PP \D_5$ and $\PP \A_5$ are both subsets of $\ov{\PP \A}_3$, we can consider closures 
inside $\ov{\PP \A}_3$. 
\begin{claim}
\label{claim_a4_closure_simultaneous}
Let $(\ff,\lp) \in \PP \D_5$.
Then there exists a solution 
$(\ff(t), \lp(t) ) \in \ov{\PP \A}_3$ \textit{near} $(\ff,\lp)$ to the set of equations
\begin{align}
\label{pd5_neighborhood_inside_a4}
\us_{\PP \D_4}( \ff(t), \lp(t) ) & \neq 0,  ~~\us_{\PP \A_4}( \ff(t), \lp(t)) = 0.  
\end{align}
Moreover, \textit{whenever} such a solution $(\ff(t), \lp(t))$ is sufficiently close to $(\ff,\lp)$ 
it lies in $\PP \A_4$, i.e.,
$\us_{\PP \A_5}( \ff(t), \lp(t)) \neq 0.$ 
In particular $(\ff(t), \lp(t))$ \textit{does not} lie in $ \PP \A_5$, 
\end{claim}
It is easy to see that claim \ref{claim_a4_closure_simultaneous} proves statements  \eqref{pa4_cap_pd4_supset_pd5} and \eqref{pa5closure_cap_pd5_empty} simultaneously. \\

\pf Let $v \in \G$, $\w \in \pi^*T\P^2/\G$ and $f_{ij}$ be as defined in \eqref{abbreviation}, subsection \ref{summary_sections_of_vector_bundle_definitions}. Equation \eqref{pd5_neighborhood_inside_a4} is a \textit{functional} equation since 
the quantities $\us_{\PP \D_4}( \ff(t), \lp(t) )$ and  $\us_{\PP \A_4}( \ff(t), \lp(t))$ are \textit{functionals}, i.e., they act on vectors $v$ and $\w$ and produce a number. We will first solve the corresponding equation 
\begin{align}
\label{pd5_neighborhood_inside_a4_numbers}
\{ \us_{\PP \D_4}(\ff(t), \lp(t))\}(p \otimes f^{\otimes d} \otimes \w^{\otimes 2})  = f_{02}(t) & \neq 0  \nonumber \\
\{ \us_{\PP \A_4}(\ff(t), \lp(t))\}(p \otimes f^{\otimes d} \otimes v^{\otimes 2} \otimes \w)  = f_{02}(t) \A^{f(t)}_4 & = 0.  
\end{align}
In \eqref{pd5_neighborhood_inside_a4_numbers} equality holds as \textit{numbers}. It is easy to see that the only solutions to \eqref{pd5_neighborhood_inside_a4_numbers} 
are of the form 
\begin{align}
\label{pd5_neighborhood_inside_a4_numbers_solution}
f_{21}(t) & = u, \qquad f_{02}(t) = \frac{3 u^2}{f_{40}(t)}.
\end{align}
Equation \eqref{pd5_neighborhood_inside_a4_numbers_solution} implies that the only solutions to the functional equation \eqref{pd5_neighborhood_inside_a4} 
is of the form
\begin{align}
\label{pd5_neighborhood_inside_a4_functional_solution}
\us_{\PP \D_5}^{\UL}(\ff(t), \lp(t)) & = t, \qquad \us_{\PP \D_4}(\ff(t), \lp(t))  = \frac{3 t^2}{ \us_{\PP \D_6}(\ff(t), \lp(t)) }.
\end{align}
where equality holds as \textit{functionals}.
\begin{rem} 
To avoid confusion, let us explain our notation carefully. The notation $(\ff(t), \lp(t))$ is simply used to indicate that $(\ff(t), \lp(t))$ is some point sufficiently close to $(\ff, \lp)$. We denote the \textit{functional} $\us_{\PP \D_5}^{\UL}(\ff(t), \lp(t))$ by the letter $t$. Hence, $t$ is close to the \textit{zero functional}, since $\us_{\PP \D_5}^{\UL}(\ff, \lp) =0$. Next, we denote the \textit{number} 
$$ \{\us_{\PP \D_5}^{\UL}(\ff(t), \lp(t))\}(v^{\otimes 2} \otimes \w) $$
by the symbol $f_{21}(t)$. This \textit{number} is close to the number zero. We also need to denote this number with some symbol. We decided to use the symbol $u$. Hence, we have this seemingly awkward equation $ f_{21}(t) = u.$
\end{rem}
To summarize, $t$ is \textit{functional}, while $u$ is a \textit{number}. Now comes a crucial observation: the sections 
$$\us_{\PP \D_4}: \ov{\PP \A}_3 \lra \UL_{\PP \D_4} \qquad \textnormal{and} \qquad  \us_{\PP \D_5}^{\UL}: \us_{\PP \D_4}^{-1}(0) \lra \UL_{\PP \D_5}$$ 
are transverse to the zero set (cf .Proposition \ref{pa3_cl_vanish_cor}). Since $\us_{\PP \D_4}^{-1}(0)$ is a smooth manifold, we can extend the section $\us_{\PP \D_5}^{\UL}$ outside 
a small neighborhood of $\us_{\PP \D_4}^{-1}(0)$ using the exponential map
(recall that $\us_{\PP \D_5}^{\UL}$ is well defined only on $\us_{\PP \D_4}^{-1}(0)$).
Therefore, there exists a solution $(\ff(t), \lp(t))$ close to $(\ff, \lp)$ to \eqref{pd5_neighborhood_inside_a4_functional_solution}. 
This proves our first assertion. Now we need to show that any such solution satisfies the condition 
$\us_{\PP \A_5}( \ff(t), \lp(t)) \neq 0$ if $t$ is sufficiently small. Observe that   
\begin{align}
\label{pa5_around_pd5_again}
f_{02}(t)^2 \A^{f(t)}_5 & = 15 f_{12}(t)  u^2 + O(u^3) \qquad \textnormal{using \eqref{pd5_neighborhood_inside_a4_numbers_solution}}. \nonumber \\
\implies \us_{\PP \A_5}( \ff(t), \lp(t)) & = 15 \us_{\PP \E_6}( \ff(t), \lp(t)) t^2 + O(t^3)
\end{align}    
Since $(\ff, \lp) \in \PP \D_5$, we get that $\us_{\PP \E_6}( \ff, \lp)\neq 0$ (cf. Proposition \ref{D5_Condition_prp}). Hence, by \eqref{pa5_around_pd5_again}, if $t$ is sufficiently small then  $\us_{\PP \A_5}( \ff(t), \lp(t)) \neq 0$. This proves claim \ref{claim_a4_closure_simultaneous}.  \qed \\ 

\hf\hf Before proving the next Lemma, we prove a corollary which will be used in the proof of \eqref{algopa5}. Since this corollary follows immediately from the previous discussion, we prove it here.   
\begin{cor}
\label{mult_of_pa5_section_around_pd5}
Let $\W \lra \D \times \P T\P^2$ be a vector bundle such that 
the rank of $\W$ is same as the dimension of $\PP \D_5$ and 
$\Q: \D \times \P T\P^2 \lra \W$  a \textit{generic} 
smooth section. Suppose $(\ff,\lp) \in \PP \D_5 \cap \Q^{-1}(0)$. 
Then the section $$ \us_{\PP \A_5} \oplus \Q: \ov{\PP \A}_4 \lra \UL_{\PP \A_5} \oplus \W$$
vanishes around $(\ff, \lp)$ with a multiplicity of $2$.
\end{cor}
\pf  Since the section $\Q$ is generic, $\Q^{-1}(0)$ intersects $ \PP \D_5$ \textit{transversely}. 
Since the rank of $\W$ is equal to the dimension of $\PP \D_5$ there exists a 
\textit{unique} solution $(\ff(t), \lp(t)) \in \ov{\PP \A}_3$ near $(\ff, \lp)$ to the set of equations 
\begin{align*}
\us_{\PP \D_5}^{\UL}(\ff(t), \lp(t)) & = t, \qquad \us_{\PP \D_4}(\ff(t), \lp(t))  = \frac{3 t^2}{ \us_{\PP \D_6}(\ff(t), \lp(t)) }, \qquad \Q(\ff(t), \lp(t)) =0.
\end{align*}
The claim follows from \eqref{pd5_neighborhood_inside_a4_functional_solution} combined with the added condition $\Q(\ff(t), \lp(t)) =0$. Equation \eqref{pa5_around_pd5_again} now proves our claim, since $\us_{\PP \E_6}(\ff, \lp) \neq 0.$ \qed 
\begin{rem}
This idea is due to Aleksey Zinger - the crucial observation that we can use the transversality of the bundle sections to describe the neighborhood of a point.    
\end{rem}
\textbf{Proof of Lemma \ref{cl} (\ref{A5cl}):}  By Lemma \ref{cl2} applied to
\bgd
M=\ov{\PP \A}_3,\,\,\ts_0=\us_{\PP\A_4},\,\,\ts_1=\us_{\PP\A_5},\,\,\ts_2=\us_{\PP\A_6},\,\,\ts_3=\us_{\PP\A_7},\,\,\nts=\us_{\PP\D_4},
\edd
and Proposition \ref{Ak_Condition_prp}, it suffices to show that 
\begin{align}
\{ (\ff,\lp) \in \ov{\PP \A}_5: \us_{\PP \D_4}(\ff,\lp) = 0 \} = \ov{\PP \D}_6 \cup \ov{\PP \E}_6 . \label{pa5_cap_pd4_equal_pd6_and_pe6}
\end{align}
\ni We will do this in two steps. We will show that 
\begin{align} 
\{ (\ff,\lp) \in \ov{\PP \A}_5: \us_{\PP \D_4}(\ff,\lp) = 0, ~~\us_{\PP \E_6}(\ff, \lp) \neq 0 \} 
& = \{ (\ff,\lp) \in \ov{\PP \D}_6: ~~\us_{\PP \E_6}(\ff,\lp) \neq 0  \}  \label{pa5_cap_pd4_pe6_not_zero_equal_pd6} \\ 
\{ (\ff,\lp) \in \ov{\PP \A}_5: \us_{\PP \D_4}(\ff,\lp) = 0, ~~\us_{\PP \E_6}(\ff,\lp) = 0 \} & = \ov{\PP \E}_6. \label{pa5_cap_pd4_pe6_equal_zero_equal_pe6}
\end{align}
It follows that \eqref{pa5_cap_pd4_pe6_not_zero_equal_pd6} and 
\eqref{pa5_cap_pd4_pe6_equal_zero_equal_pe6} imply \eqref{pa5_cap_pd4_equal_pd6_and_pe6}. To see this, note that 
\begin{align}
\label{pd6_closure_psi_pe6_vanish_subset_of_pe6}
\{ (\ff,\lp) \in \ov{\PP \D}_6: \us_{\PP \E_6}(\ff,\lp) = 0 \} & \subset \ov{\PP \E}_6 \qquad \textnormal{using Corollary \ref{pd6_cl_vanish_cor} and \ref{pe6_cl_vanish_cor}}. 
\end{align}
It is now easy to see that \eqref{pd6_closure_psi_pe6_vanish_subset_of_pe6}, \eqref{pa5_cap_pd4_pe6_not_zero_equal_pd6} and \eqref{pa5_cap_pd4_pe6_equal_zero_equal_pe6} imply \eqref{pa5_cap_pd4_equal_pd6_and_pe6}. \\
\hf We will now start with the proof of \eqref{pa5_cap_pd4_pe6_not_zero_equal_pd6}. We will show that the left hand side of \eqref{pa5_cap_pd4_pe6_not_zero_equal_pd6} is a subset of the right hand side. This follows from \eqref{pa5closure_cap_pd5_empty}. By Lemma \ref{cl}, statement \ref{A4cl} we know that $\ov{\PP \A}_5 \subset \ov{\PP \A}_4$. Therefore, in conjunction with \eqref{pa4_cap_pd4_equal_pd5}, we gather that
\bgd
 \{ (\ff,\lp) \in \ov{\PP \A}_5: \us_{\PP \D_4}(\ff, \lp) = 0 \} \subset \{ (\ff,\lp) \in \ov{\PP \A}_4: \us_{\PP \D_4}(\ff,\lp) = 0 \}=  \ov{\PP \D}_5. 
\edd
The right hand side above equals $\PP \D_5 \cup \ov{\PP \D}_6 \cup \ov{\PP \E}_6$  by Lemma \ref{cl}, statement \ref{D5cl}. Hence, by \eqref{pa5closure_cap_pd5_empty} 
\bgd
\{ (\ff,\lp) \in \ov{\PP \A}_5: \us_{\PP \D_4}(\ff,\lp) = 0, ~\us_{\PP \E_6}(\ff,\lp) \neq 0 \} \subset \{ (\ff,\lp) \in \ov{\PP \D}_6: \us_{\PP \E_6}(\ff, \lp) \neq 0 \}.
\edd
\hf Now we will show the converse. Since $\ov{\PP \A}_5$ is a closed space, it suffices to show that  
\begin{align}
\{ (\ff,l_p) \in \ov{\PP \A}_5: \us_{\PP \D_4}(\ff,l_p) = 0, ~\us_{\PP \E_6}(\ff,l_p) \neq 0 \} 
\supset  \{ (\ff,l_p) \in \PP \D_6: ~\us_{\PP \E_6}(\ff,l_p) \neq 0  \}.  \label{pa5_cap_pd4_pe6_neq_zero_supset_pe6}
\end{align}
As before, we will simultaneously prove \eqref{pa5_cap_pd4_pe6_neq_zero_supset_pe6} and also prove that
\bge
\label{pa6closure_cap_pd6_empty}
\ov{\PP \A}_6 \cap \PP \D_6 = \varnothing.   
\ede

\begin{claim}
\label{claim_a5_closure_simultaneous_pd6}
Let $(\ff,\lp) \in \PP \D_6$.
Then there exists a solution 
$(\ff(t), \lp(t) ) \in \ov{\PP \A}_3$ near $(\ff,\lp)$ to the set of equations
\begin{align}
\us_{\PP \D_4}( \ff(t), \lp(t) ) & \neq 0,  ~~\us_{\PP \A_4}( \ff(t), \lp(t)) = 0, ~~\us_{\PP \A_5}( \ff(t), \lp(t)) = 0. \label{pd6_neighborhood_inside_a5} 
\end{align}
Moreover, \textit{whenever} such a solution $(\ff(t), \lp(t))$ is sufficiently close to $(\ff,\lp)$ 
it lies in $\PP \A_5$, i.e., $\us_{\PP \A_6}( \ff(t), \lp(t)) \neq 0.$ In particular, $(\ff(t), \lp(t))$ \textit{does not} lie in $ \PP \A_6$.
\end{claim}
It is clear that claim \ref{claim_a5_closure_simultaneous_pd6} proves \eqref{pa5_cap_pd4_pe6_neq_zero_supset_pe6} and \eqref{pa6closure_cap_pd6_empty} simultaneously. \\

\pf As before, we will first solve the equation 
\begin{align}
\label{pd6_neighborhood_inside_a5_numbers}
f_{02}(t) \neq 0, ~~f_{02}(t) \A^{f(t)}_4 = 0, ~~f_{02}(t)^2 \A^{f(t)}_5  = 0.
\end{align}
The only solutions to \eqref{pd6_neighborhood_inside_a5_numbers} are of the form 
\begin{align}
\label{pd6_neighborhood_inside_a5_numbers_solution}
f_{02}(t) & = u  \nonumber \\
f_{21}(t) &= \left( \frac{5 f_{31}(t) \pm \sqrt{-15 f_{12}(t) \D^{f(t)}_7 }}{15 f_{12}(t)} \right) u\nonumber\\
f_{40}(t) &=  3 \left( \frac{5 f_{31}(t) \pm \sqrt{ -15 f_{12}(t) \D^{f(t)}_7 }}{15 f_{12}(t)} \right)^2 u .
\end{align}
The second equation comes from solving a quadratic arising from $f_{02}(t)^2 \A^{f(t)}_5  = 0$ while the third is from solving $f_{02}(t) \A^{f(t)}_4  = 0$ and using $f_{21}$ from the second equation. Since $(\ff, \lp) \in \PP \D_6 $ we know that $f_{12} \neq 0$ and $\D^{f}_7 \neq 0$. Hence, there are always \textit{two} solutions for $(f_{21}(t), f_{40}(t))$. Equation \eqref{pd6_neighborhood_inside_a5_numbers_solution} implies that the only solutions to the functional equation \eqref{pd6_neighborhood_inside_a5} is of the form
\begin{align}
\label{pd6_neighborhood_inside_a5_functional_solution}
\us_{\PP \D_4}(\ff(t), \lp(t)) & = t \nonumber \\
\us_{\PP \D_5}^{\UL}(\ff(t), \lp(t)) & = \left( \frac{5 \us_{\PP \E_8}(\ff(t), \lp(t))  \pm \sqrt{-15 \us_{\PP \E_6}(\ff(t), \lp(t)) \us_{\PP \D_7}(\ff(t), \lp(t)) } }{15 \us_{\PP \E_6}(\ff(t), \lp(t))} \right) t 
\nonumber \\
\us_{\PP \D_6}(\ff(t), \lp(t)) & = 3 \left( \frac{5 \us_{\PP \E_8}(\ff(t), \lp(t))  \pm \sqrt{-15 \us_{\PP \E_6}(\ff(t), \lp(t)) \us_{\PP \D_7}(\ff(t), \lp(t)) }  }{15 \us_{\PP \E_6}(\ff(t), \lp(t))} \right)^2 t
\end{align}
where equality holds as \textit{functionals}. Since the sections 
$$\us_{\PP \D_4}: \ov{\PP \A}_3 \lra \UL_{\PP \D_4}, \qquad  \us_{\PP \D_5}^{\UL}: \us_{\PP \D_4}^{-1}(0) \lra \UL_{\PP \D_5} \qquad \textnormal{and} \qquad  
\us_{\PP \D_6}: \us_{\PP \D_5}^{\UL^{-1}}(0) \lra \UL_{\PP \D_6}$$ 
are transverse to the zero set (as proved in Proposition \ref{D4_Condition_prp} and \ref{D5_Condition_prp}), there exists 
a solution $(\ff(t), \lp(t))$ close to $(\ff, \lp)$ to \eqref{pd6_neighborhood_inside_a5_functional_solution}. 
This proves our first assertion. Next we need to show that any such solution satisfies the condition 
$\us_{\PP \A_6}( \ff(t), \lp(t)) \neq 0$ if $t$ is sufficiently small. To prove that we observe  
\begin{align}
\label{pa6_around_pd6_again}
f_{02}(t)^3 \A^{f(t)}_3 & = \frac{\D^{f(t)}_7}{f_{12}(t)} u^2 + O(u^3) \qquad \textnormal{using \eqref{pd6_neighborhood_inside_a5_numbers_solution}, for either choice of $\sqrt{f_{12} \D^f_7}$ }. \nonumber \\
\implies \us_{\PP \A_6}( \ff(t), \lp(t)) & = \frac{\us_{\PP \D_7} (\ff(t), \lp(t))}{\us_{\PP \E_6} (\ff(t), \lp(t))^2} t^2 + O(t^3)
\end{align}    
Since $(\ff, \lp) \in \PP \D_6$, we get that $\us_{\PP \E_6}( \ff, \lp), \us_{\PP \D_7}( \ff, \lp) \neq 0$ (see Proposition \ref{D6_Condition_prp}). 
Hence, \eqref{pa6_around_pd6_again} implies that if $t$ is sufficiently small $\us_{\PP \A_6}( \ff(t), \lp(t)) \neq 0$, which proves claim \ref{claim_a5_closure_simultaneous_pd6}. \\ 
\hf Before proving \eqref{pa5_cap_pd4_pe6_equal_zero_equal_pe6}, we prove a corollary which will be used in the proof of \eqref{algopa6}. 
\begin{cor}
\label{mult_of_pa6_section_around_pd6}
Let $\W \lra \D \times \P T\P^2$ be a vector bundle such that 
the rank of $\W$ is same as the dimension of $\PP \D_6$. 
Let $\Q: \D \times \P T\P^2 \lra \W$ be a \textit{generic} 
smooth section. 
Suppose $(\ff,\lp) \in \PP \D_6 \cap \Q^{-1}(0)$. 
Then the section $$ \us_{\PP \A_6} \oplus \Q: \ov{\PP \A}_5 \lra \UL_{\PP \A_6} \oplus \W$$
vanishes around $(\ff, \lp)$ with a multiplicity of $4$.
\end{cor}
\pf As before, in the proof of claim \ref{mult_of_pa5_section_around_pd5}, this follows from \eqref{pa6_around_pd6_again} and the fact 
that $\Q^{-1}(0)$ intersects $\PP \D_6$ transversely. Each branch of $\sqrt{f_{12} \D^f_7}$ contributes with 
a multiplicity of $2$. Hence, the total multiplicity is $4$.   \qed  \\

\hf Next we will prove \eqref{pa5_cap_pd4_pe6_equal_zero_equal_pe6}. We will show that the left hand side is a 
subset of the right hand side. Note that $\ov{\PP \A}_5 \subset \ov{\PP \A}_4$ is implied by Lemma \ref{cl}, statement \ref{A4cl}, whence 
\bgd
\{ (\ff,\lp) \in \ov{\PP \A}_5: \us_{\PP \D_4}(\ff,\lp) = 0, \us_{\PP \E_6}(\ff,\lp) = 0 \} 
\subset \{ (\ff,\lp) \in \ov{\PP \A}_4: \us_{\PP \D_4}(\ff,\lp) = 0, \us_{\PP \E_6}(\ff,\lp) = 0 \}. 
\edd
The right hand side above can be simplified  \eqref{pa4_cap_pd4_equal_pd5} and Corollary \ref{pe6_cl_vanish_cor} as follows :
\bgd
\{ (\ff,\lp) \in \ov{\PP \A}_4: \us_{\PP \D_4}(\ff,\lp) = 0, \us_{\PP \E_6}(\ff,\lp) = 0 \}= \{ (\ff,\lp) \in \ov{\PP \D}_5: \us_{\PP \E_6}(\ff,\lp) = 0 \}= \ov{\PP \E}_6.
\edd
This implies
\bgd
\{ (\ff,\lp) \in \ov{\PP \A}_5: \us_{\PP \D_4}(\ff,\lp) = 0, ~\us_{\PP \E_6}(\ff,\lp) = 0 \}  \subset \ov{\PP \E}_6.
\edd
Now we prove the converse. Since $\ov{\PP \A}_5$ is a closed set, it suffices to show that      
\begin{align}
\{ (\ff,\lp) \in \ov{\PP \A}_5: \us_{\PP \D_4}(\ff,\lp) = 0, ~~\us_{\PP \E_6}(\ff,\lp) = 0 \} 
\supset  \PP \E_6.  \label{pa5_cap_pd4_pe6_eq_zero_supset_pe6}
\end{align}
As before, we will simultaneously prove \eqref{pa5_cap_pd4_pe6_eq_zero_supset_pe6} and also prove the following:
\bge
\label{pa6closure_cap_pe6_empty}
\ov{\PP \A}_6 \cap \PP \E_6 = \varnothing.   
\ede

\begin{claim}
\label{claim_a5_closure_simultaneous_pe6}
Let $(\ff,\lp) \in \PP \E_6$. Then there exists a solution 
$(\ff(t), \lp(t) ) \in \ov{\PP \A}_3$ near $(\ff,\lp)$ to the set of equations
\begin{align}
\us_{\PP \D_4}( \ff(t), \lp(t) ) & \neq 0,  ~~\us_{\PP \A_4}( \ff(t), \lp(t)) = 0, ~~\us_{\PP \A_5}( \ff(t), \lp(t)) = 0. \label{pe6_neighborhood_inside_a5} 
\end{align}
Moreover, \textit{whenever} such a solution $(\ff(t), \lp(t))$ is sufficiently close to $(\ff ,\lp)$, 
it lies in $\PP \A_5$, i.e., $\us_{\PP \A_6}( \ff(t), \lp(t)) \neq 0.$ In particular $(\ff(t), \lp(t))$ \textit{does not} lie in $ \PP \A_6$.
\end{claim}
Note that claim \ref{claim_a5_closure_simultaneous_pe6} proves \eqref{pa5_cap_pd4_pe6_eq_zero_supset_pe6} and \eqref{pa6closure_cap_pe6_empty} 
simultaneously. \\

\pf As before, we will first solve the equations 
\begin{align}
\label{pe6_neighborhood_inside_a5_numbers}
f_{02}(t) \neq 0, ~~f_{02}(t) \A^{f(t)}_4 = 0, ~~f_{02}(t)^2 \A^{f(t)}_5  = 0.
\end{align}
It is easy to see that the only solutions to \eqref{pe6_neighborhood_inside_a5_numbers} 
are  of the form 
\begin{align}
\label{pe6_neighborhood_inside_a5_numbers_solution}
f_{21}(t) = u, \qquad  f_{02} = \frac{3 u^2}{f_{40}(t)}, \qquad \textnormal{and} \qquad 
f_{12} =  \frac{2 f_{31}(t)}{f_{40}(t)} u - \frac{3 f_{50}(t)}{5 f_{40}(t)^2} u^2. 
 \end{align}
Note that since $(\ff, \lp) \in \PP \E_6 $ we get that $f_{40} \neq 0$. 
Equation \eqref{pe6_neighborhood_inside_a5_numbers_solution} implies that the only solutions to the functional equation \eqref{pe6_neighborhood_inside_a5} 
is of the form
\begin{align}
\label{pe6_neighborhood_inside_a5_functional_solution}
\us_{\PP \D_5}^{\UL}(\ff(t), \lp(t)) & = t \nonumber \\
\us_{\PP \D_4}(\ff(t), \lp(t)) & = \frac{3 t^2}{ \us_{\PP \D_6}(\ff(t), \lp(t))}  \nonumber \\
\us_{\PP \E_6}(\ff(t), \lp(t)) & =  \frac{2 \us_{\PP \E_8}(\ff(t), \lp(t)) }{ \us_{\PP \D_6}(\ff(t), \lp(t)) } t + O(t^2) 
\end{align}
where equality holds as \textit{functionals}. Since the sections 
\bgd
\us_{\PP \D_4}: \ov{\PP \A}_3 \lra \UL_{\PP \D_4}, \qquad  \us_{\PP \D_5}^{\UL}: \us_{\PP \D_4}^{-1}(0) \lra \UL_{\PP \D_5} \qquad \textnormal{and} \qquad  
\us_{\PP \E_6}: \us_{\PP \D_5}^{\UL^{-1}}(0) \lra \UL_{\PP \E_6}
\edd
are transverse to the zero set (as proved in Proposition \ref{D4_Condition_prp} and \ref{D5_Condition_prp}), there exists 
a solution $(\ff(t), \lp(t))$ close to $(\ff, \lp)$ to \eqref{pe6_neighborhood_inside_a5_functional_solution}. 
This proves our first assertion. Next we need to show that any such solution satisfies the condition 
$\us_{\PP \A_6}( \ff(t), \lp(t)) \neq 0$ if $t$ is sufficiently small. To prove that we observe  
\begin{align}
\label{pa6_around_pe6_again}
f_{02}(t)^3 \A^{f(t)}_6 & = -15 f_{03}(t) u^3 + O(u^4) \qquad \textnormal{using \eqref{pe6_neighborhood_inside_a5_numbers_solution}.} \nonumber \\
\implies \us_{\PP \A_6}( \ff(t), \lp(t)) & = -15 \us_{\PP \XC_8}(\ff(t), \lp(t)) t^3 + O(t^3) 
\end{align}    
Since $(\ff, \lp) \in \PP \E_6$, we get that $\us_{\PP \XC_6}( \ff, \lp) \neq 0$ (see Proposition \ref{E6_Condition_prp}). 
Hence, \eqref{pa6_around_pe6_again} implies that if $t$ is sufficiently small $\us_{\PP \A_6}( \ff(t), \lp(t)) \neq 0$ 
which proves claim \ref{claim_a5_closure_simultaneous_pe6} \qed \\  

\hf We now prove a corollary which will be used in proving \eqref{algopa6}. 
\begin{cor}
\label{mult_of_pa6_section_around_pe6}
Let $\W \lra \D \times \P T\P^2$ be a vector bundle such that 
the rank of $\W$ is same as the dimension of $\PP \E_6$. Let $\Q: \D \times \P T\P^2 \lra \W$ be a \textit{generic} smooth section. Suppose $(\ff,\lp) \in \PP \E_6 \cap \Q^{-1}(0)$. Then the section $$ \us_{\PP \A_6} \oplus \Q: \ov{\PP \A}_5 \lra \UL_{\PP \A_6} \oplus \W$$
vanishes around $(\ff, \lp)$ with a multiplicity of $3$.
\end{cor}
\pf Follows from the fact that $\Q^{-1}(0)$ intersects $\PP \E_6$ transversely and \eqref{pa6_around_pe6_again}. \qed  \\

\textbf{Proof of Lemma \ref{cl} (\ref{A6cl}):} By Lemma \ref{cl2} applied to
\bgd
M=\ov{\PP \A}_3,\,\,\ts_0=\us_{\PP\A_4},\,\,\ts_1=\us_{\PP\A_5},\,\,\ts_2=\us_{\PP\A_6},\,\,\ts_3=\us_{\PP\A_7},\,\,\ts_4=\us_{\PP\A_8}, \,\,\nts=\us_{\PP\D_4},
\edd
and Proposition \ref{Ak_Condition_prp}, it suffices to show that 
\begin{align}
\{ (\ff,\lp) \in \ov{\PP \A}_6: \us_{\PP \D_4}(\ff,\lp) = 0 \} = \ov{\PP \D}_7 \cup \ov{\PP \E}_7 \cup \ov{\XX^{\# \flat}_8} . \label{pa6_cap_pd4_equal_pd7_and_pe7_and_x8}
\end{align}
\ni We will do this in three steps. We will show that 
\begin{align} 
\{ (\ff,\lp) \in \ov{\PP \A}_6 &: \us_{\PP \D_4}(\ff,\lp) = 0, ~\us_{\PP \E_6}(\ff,\lp) \neq 0 \} 
 \equiv \{ (\ff,\lp) \in \ov{\PP \D}_7: ~\us_{\PP \E_6}(\ff,\lp) \neq 0  \} \label{pa6_cap_pd4_pe7_not_zero_equal_pd7} \\
\{ (\ff,\lp) \in \ov{\PP \A}_6 &: \us_{\PP \D_4}(\ff,\lp) = 0, ~\us_{\PP \E_6}(\ff,\lp) = 0,  \nonumber \\ 
                               & \qquad \qquad \qquad \qquad ~\us_{\PP \XC_8}(\ff,\lp) \neq  0 \} \equiv
\{ (\ff,\lp) \in \ov{\PP \E}_7: ~\us_{\PP \XC_8}(\ff, \lp) \neq 0  \} \label{pa6_cap_pd4_pe7_equal_zero_px8_not_zero_equal_pe7} \\
\{ (\ff,\lp) \in \ov{\PP \A}_6 &: \us_{\PP \D_4}(\ff,\lp) = 0, ~\us_{\PP \E_6}(\ff,\lp) = 0, ~\us_{\PP \XC_8}(\ff, \lp) = 0 \} \equiv
\ov{\XX^{\# \flat}_8}. \label{pa6_cap_pd4_pe7_equal_zero_px8_equal_zero_px8}
\end{align}
Observe that \eqref{pa6_cap_pd4_pe7_not_zero_equal_pd7}, \eqref{pa6_cap_pd4_pe7_equal_zero_px8_not_zero_equal_pe7} and 
\eqref{pa6_cap_pd4_pe7_equal_zero_px8_equal_zero_px8} prove \eqref{pa6_cap_pd4_equal_pd7_and_pe7_and_x8}. To see this, note that 
\begin{align}
\{ (\ff,\lp) \in \ov{\PP \E}_7: \us_{\PP \XC_8}(\ff, \lp) = 0 \} & \subset \ov{\XX^{\# \flat}_8} \label{pe7_closure_psi_x8_vanish_subset_of_x8} \\ 
\{ (\ff,\lp) \in \ov{\PP \D}_7: \us_{\PP\E_6}(\ff, \lp) = 0 \} & \subset \ov{\PP \E}_7 \cup  \ov{\XX^{\# \flat}_8}. \label{pd7_closure_psi_pe7_vanish_subset_of_pe7}
\end{align}
Equation \eqref{pe7_closure_psi_x8_vanish_subset_of_x8} follows from Corollary \ref{pe7_cl_vanish_cor} and \ref{X8_and_X8_sharp_up_cl_vanish_cor}. 
Equation \eqref{pd7_closure_psi_pe7_vanish_subset_of_pe7} follows from Proposition \ref{Dk_Condition_prp} and Corollary \ref{pe7_cl_vanish_cor} and \ref{X8_and_X8_sharp_up_cl_vanish_cor}. It is now easy to see that \eqref{pa6_cap_pd4_pe7_not_zero_equal_pd7}, \eqref{pa6_cap_pd4_pe7_equal_zero_px8_not_zero_equal_pe7} and 
\eqref{pa6_cap_pd4_pe7_equal_zero_px8_equal_zero_px8} combined with \eqref{pe7_closure_psi_x8_vanish_subset_of_x8} and \eqref{pd7_closure_psi_pe7_vanish_subset_of_pe7} 
prove \eqref{pa6_cap_pd4_equal_pd7_and_pe7_and_x8}. \\
\hf Let us prove \eqref{pa6_cap_pd4_pe7_not_zero_equal_pd7}. To see why the left hand side is a subset of the right hand side, recall \eqref{pa6closure_cap_pd6_empty}. We also recall a result:
\bgd
 \ov{\hat{\XC}_8^{\# \flat}} = \{ (\ff, \lp)\in \D \times \P T\P^2:    
\us_{\AA_0} (\ff, \lp) =0,\us_{\AA_1}(\ff, \lp) =0, \us_{ \DD_4}(\ff,\lp) =0, \us_{\XX_8}(\ff,\lp)=0  \}.
\edd
Now observe that $\ov{\PP \A}_6 \subset \ov{\PP \A}_5$ is implied by Lemma \ref{cl}, statement \ref{A5cl}, whence 
\bgd
\{ (\ff,\lp) \in \ov{\PP \A}_6: \us_{\PP \D_4}(\ff,\lp) = 0 \} \subset \{ (\ff,\lp) \in \ov{\PP \A}_5: \us_{\PP \D_4}(\ff,\lp) = 0 \}.
\edd
By \eqref{pa5_cap_pd4_equal_pd6_and_pe6}, the right hand side above equals $\ov{\PP \D}_6 \cup \ov{\PP \E}_6$. But by Lemma \ref{cl}, statement \ref{D6cl} and statement \ref{E6cl}, we get
\bgd
\ov{\PP \D}_6 \cup \ov{\PP \E}_6= \PP \D_6 \cup \PP \E_6 \cup \ov{\PP \D}_7 \cup \ov{\PP \E}_7 \cup \ov{\XX^{\#}_8}.
\edd
This implies, by \eqref{pa6closure_cap_pd6_empty} and \eqref{pa6closure_cap_pe6_empty}, that 
\bgd
\{ (\ff,\lp) \in \ov{\PP \A}_6: \us_{\PP \D_4}(\ff,\lp) = 0, ~\us_{\PP \E_6}(\ff,\lp) \neq 0 \} \subset \ov{\PP \D}_7. 
\edd
Hence, the left hand side of \eqref{pa6_cap_pd4_pe7_not_zero_equal_pd7} is a subset of its right hand side.\\
\hf Next let us show that the right hand side of \eqref{pa6_cap_pd4_pe7_not_zero_equal_pd7} is a subset of its left hand side. Since $\ov{\PP \A}_6$ is a closed set, it suffices to show that  
\begin{align}
\{ (\ff,\lp) \in \ov{\PP \A}_6: \us_{\PP \D_4}(\ff,\lp) = 0, ~~\us_{\PP \E_6}(\ff,\lp) \neq 0 \} 
\supset  \{ (\ff,\lp) \in \PP \D_7: ~~\us_{\PP \E_6}(\ff,\lp) \neq 0  \}.  \label{pa6_cap_pd4_pe7_neq_zero_supset_pe7}
\end{align}
We will simultaneously prove \eqref{pa6_cap_pd4_pe7_neq_zero_supset_pe7} and also prove the following:
\bge
\label{pa7closure_cap_pd7_empty}
\ov{\PP \A}_7 \cap \PP \D_7 = \varnothing.
\ede

\begin{claim}
\label{claim_a6_closure_simultaneous_pd6}
Let $(\ff,\lp) \in \PP \D_7$. Then there exists a solution $(\ff(t), \lp(t) ) \in \ov{\PP \A}_3$ near $(\ff,\lp)$ to the set of equations
\begin{align}
\us_{\PP \D_4}( \ff(t), \lp(t) ) & \neq 0,  \us_{\PP \A_4}( \ff(t), \lp(t)) = 0, \us_{\PP \A_5}( \ff(t), \lp(t)) = 0, \us_{\PP \A_6}( \ff(t), \lp(t)) = 0. \label{pd7_neighborhood_inside_a6} 
\end{align}
Moreover, \textit{whenever} such a solution $(\ff(t), \lp(t))$ is sufficiently close to $(\ff,\lp)$ 
it lies in $\PP \A_6$, i.e., $\us_{\PP \A_7}( \ff(t), \lp(t)) \neq 0.$ In particular $(\ff(t), \lp(t))$ \textit{does not} lie in $ \PP \A_7$.
\end{claim}
Note that claim \ref{claim_a6_closure_simultaneous_pd6} proves \eqref{pa6_cap_pd4_pe7_neq_zero_supset_pe7} and \eqref{pa7closure_cap_pd7_empty} 
simultaneously. \\

\pf As before, we will first solve the equation 
\begin{align}
\label{pd7_neighborhood_inside_a6_numbers}
f_{02}(t) \neq 0, ~~f_{02}(t) \A^{f(t)}_4 = 0, ~~f_{02}(t)^2 \A^{f(t)}_5  = 0 ~~\textnormal{and} ~~f_{02}(t)^3 \A^{f(t)}_6  = 0.
\end{align}
We claim that the \textit{only} solutions to \eqref{pd7_neighborhood_inside_a6_numbers} that go to zero as $f_{02}(t)$ goes to zero are 
of the form 
\begin{align}
f_{02}(t) & = u^2 +O(u^4) \label{pd7_neighborhood_inside_a6_numbers_solution_cubic_eqn_f02} \\
f_{21}(t) & = \frac{f_{31}(t)}{3 f_{12}(t)} u^2 + \sqrt{\frac{\beta(t)}{f_{12}(t)}} u^3 + O(u^4)   \label{pd7_neighborhood_inside_a6_numbers_solution_cubic_eqn} \\
f_{40}(t) & = \frac{f_{31}(t)^2}{3 f_{12}(t)^2} u^2 + O(u^3) \nonumber \\
\Big( f_{50}(t) - \frac{5 f_{31}(t)^2}{3 f_{12}(t)}  \Big ) & = -15 \beta(t) u^2 + O(u^3) \nonumber \\ 
\textnormal{where} \qquad \beta(t) &= -\frac{f_{03}(t) f_{31}(t)^3}{162 f_{12}(t)^4}+\frac{f_{22}(t) f_{31}(t)^2}{18 f_{12}(t)^3}-\frac{f_{41}(t) f_{31}(t)}{18 f_{12}(t)^2}+\frac{f_{60}(t)}{90 f_{12}(t)}
 \label{pd7_neighborhood_inside_a6_numbers_solution}
\end{align}
for just \textit{one} choice of a branch of $\sqrt{\beta(t)}$ 
$\footnote{In other words, choosing the other branch of the square root does not give us any extra solutions.}$. We will see shortly that $\beta(t) \neq 0$. The value for $f_{40}$ can be calculated using $f_{21}, f_{02}$ and $f_{02}(t) \A^{f(t)}_4 =0$ while the fourth equation follows by using the first three equations and $f_{02}(t)^2 \A^{f(t)}_5 =0$. Let us now explain how we obtain \eqref{pd7_neighborhood_inside_a6_numbers_solution_cubic_eqn_f02} and \eqref{pd7_neighborhood_inside_a6_numbers_solution_cubic_eqn}. The equation $f_{02}(t)^3 \A^{f(t)}_6 =0$ is a cubic equation in $f_{21}(t)$, i.e., it is of the form 
\bgd
\Z_3(f_{02}(t)) f_{21}(t)^3 + \Z_2(f_{02}(t)) f_{21}(t)^2 + \Z_1(f_{02}(t)) f_{21}(t) + \Z_0(f_{02}(t)) =0.
\edd
As $\us_{\PP \E_6}(\ff(t),\lp(t))\neq 0$, this implies that $f_{12}(t) \neq 0$. It follows that as $f_{02}(t)$ goes to zero $\Z_{2}$ remains non zero.  
Hence, there exists a unique 
holomorphic function $\Pz(f_{02}(t))$, of $f_{02}(t)$ (close to the zero function), such that 
if we make a change of variables
\bgd
f_{21}(t) = \mathrm{H} + \Pz(f_{02}(t))
\edd
then our cubic equation becomes 
\bgd
\hat{\Z}_3(f_{02}(t)) \mathrm{H}^3 + \hat{\Z}_2(f_{02}(t)) \mathrm{H}^2 + \hat{\Z}_0(f_{02}(t)) =0.
\edd
The argument is same as in Lemma \ref{fstr_prp}, where we show the existence of $\Y(x)$
(it is simply an application of the Implicit Function Theorem). In fact, we observe that
\bgd
\Pz(f_{02})=\frac{1}{3A_3}\left(-A_2+\sqrt{A_2^2-3A_1A_3}\right).
\edd 
This is defined even when $A_3=0$ as can be seen by a standard binomial expansion, i.e.,
\bgd
\Pz(f_{02}(t)) = \frac{f_{31}(t)}{3 f_{12}(t)} f_{02}(t) + O(f_{02}(t)^2) = \frac{f_{31}(t)}{3 f_{12}(t)} u^2 + O(u^4).  
\edd
The other root has the property that $\Pz(f_{02}(t))$ goes to a non-zero constant as $f_{02}(t)$ goes to zero. \\
\hf Since $\hat{\Z}_2(0) \neq 0$, we can divide out by $\hat{\Z}_2(f_{20}(t))$ and get 
\begin{align}
\hat{\hat{\Z}}_3(f_{02}(t)) \mathrm{H}^3 + \mathrm{H}^2 + \hat{\hat{\Z}}_0(f_{02}(t)) =0. \label{cubic_divided}
\end{align}
By a simple calculation, it is easy to see that 
\begin{align*}
\hat{\hat{\Z}}_0(f_{02}(t))& = -\frac{\beta(t)}{f_{12}(t)} f_{02}(t)^3 + O(f_{02}(t)^4).
\end{align*}
Assuming $\beta(t) \neq 0$ we can make a change of coordinates 
\bgd
 \hat{f}_{02} = f_{02}(t) \bigg( \frac{f_{12}(t)\,\hat{\hat{\Z}}_0(f_{02}(t))}{-\beta(t) f_{02}(t)^3 } \bigg)^{\frac{1}{3}},  \qquad \hat{\mathrm{H}} = \mathrm{H}(1 + \hat{\hat{\Z}}_3(f_{02}(t)) \mathrm{H})^{\frac{1}{2}}.
\edd
Our cubic equation \eqref{cubic_divided} now becomes 
\begin{align}
\hat{\mathrm{H}}^2 - \frac{\beta(t)}{f_{12}(t)} \hat{f}_{02}^3 =0. \label{cubic_simplified}
\end{align}
Now, it is easy to see that the \textit{only} small solutions to \eqref{cubic_simplified} are of the form
\bgd
\hat{\mathrm{H}} = \sqrt{\frac{\beta(t)}{f_{12}(t)}} u^3, \qquad \hat{f}_{02} = u^2
\edd
for just \textit{one} choice of $\sqrt{\beta(t)}$. In other words, by choosing just \textit{one} branch of $\sqrt{\beta(t)}$, we get \textit{all} the possible small solutions of \eqref{cubic_simplified}. 
By inverting the change of coordinates, $(\mathrm{H}, f_{02}) \lra  (\hat{\mathrm{H}}, \hat{f}_{02})$, we conclude that the \textit{only} small solutions to 
\eqref{cubic_divided} are of the form  
\bgd
\mathrm{H} = \sqrt{\frac{\beta(t)}{f_{12}(t)}} u^3 + O(u^4), \qquad f_{02}(t) = u^2 + O(u^4).
\edd
(Note that the transformation $(\mathrm{H}, f_{02}) \lra  (\hat{\mathrm{H}}, \hat{f}_{02})$ is identity to first order, i.e. the Jacobian matrix of this transformation at the 
origin is the identity matrix.) Combining the expression for $\Pz(f_{02})$ and $\mathrm{H}$ gives us \eqref{pd7_neighborhood_inside_a6_numbers_solution_cubic_eqn} and \eqref{pd7_neighborhood_inside_a6_numbers_solution_cubic_eqn_f02}. 
It remains to show that $\beta(t) \neq 0$.  
To see this, note that
\begin{align}
\label{D8_D7_a1}
\beta(t) &= \frac{\D^{f(t)}_8}{90 f_{12}(t)} -  \frac{f_{30}(t) f_{31}(t) \D^{f(t)}_7}{54 f_{12}(t)^3}.
\end{align}
Since $(\ff, \lp) \in \PP \D_7$, $\D^{f}_7 =0$ and $\D^{f}_8 \neq 0$. 
Therefore, by \eqref{D8_D7_a1} $\beta(t) \neq 0$ for small $t$. \\
\hf Equation \eqref{pd7_neighborhood_inside_a6_numbers_solution} combined with \eqref{D8_D7_a1} imply that the only 
solutions to the functional equation \eqref{pd7_neighborhood_inside_a6} are of the form 
\begin{align}
\label{pd7_neighborhood_inside_a6_functional_solution}
\us_{\PP \D_4}( \ff(t), \lp(t)) &= t^2 +O(t^4)\nonumber \\
\us_{\PP \D_5}^{\UL}( \ff(t), \lp(t)) &=  \frac{\us_{\PP \E_8}( \ff(t), \lp(t))}{3 \us_{\PP \E_6}( \ff(t), \lp(t))} t^2 + \sqrt{\Bz(\ff(t), \lp(t))} t^3 + O(t^4) \nonumber   \\
\us_{\PP \D_6}( \ff(t), \lp(t)) &= \frac{\us_{\PP \E_8}( \ff(t), \lp(t))^2}{3 \us_{\PP \E_6}( \ff(t), \lp(t))^2} t^2 + O(t^3)  \nonumber    \\ 
\us_{\PP \D_7}( \ff(t), \lp(t)) &= -15 \us_{\PP \E_6}( \ff(t), \lp(t))^2  \Bz(\ff(t), \lp(t)) t^2 + O(t^3),\\
\textnormal{where} \qquad \Bz(\ff(t), \lp(t) )  & = \frac{\us_{\PP \D_8}(\ff(t), \lp(t)) }{90  \us_{\PP \E_6}(\ff(t), \lp(t))^4} - 
\frac{\us_{\PP \XC_8}(\ff(t), \lp(t))  \us_{\PP \E_8}(\ff(t), \lp(t))  \us_{\PP \D_7}(\ff(t), \lp(t)) }{54 \us_{\PP \E_6}(\ff(t), \lp(t))^4}   \nonumber
\end{align}
and equality holds as \textit{functionals}. 
Note that $\Bz$  is the functional version of $\beta$ using \eqref{D8_D7_a1}.
Since the sections 
\begin{align*}
\us_{\PP \D_4} &: \ov{\PP \A}_3 \lra \UL_{\PP \D_4},  \qquad \us_{\PP \D_5}^{\UL}: \us_{\PP \D_4}^{-1}(0) \lra \UL_{\PP \D_5},\\ 
\us_{\PP \D_6} &: \us_{\PP \D_5}^{\UL^{-1}}(0) \lra \UL_{\PP \D_6}, \qquad \us_{\PP \D_7}: \us_{\PP \D_6}^{-1}(0) - \us_{\PP \E_6}^{-1}(0)  \lra \UL_{\PP \D_7}
\end{align*}
are transverse to the zero set (as proved in Proposition \ref{D4_Condition_prp}, \ref{D5_Condition_prp} and \ref{Dk_Condition_prp} ), there exists 
a solution $(\ff(t), \lp(t)) $  close to $(\ff, \lp)$ to \eqref{pd6_neighborhood_inside_a5_functional_solution}. 
This proves our first assertion. Next we need to show that any such solution satisfies the condition 
$\us_{\PP \A_7}( \ff(t), \lp(t)) \neq 0$ if $t$ is sufficiently small. To prove that we observe  
\begin{align}
\label{pa7_around_pd7_again}
f_{02}(t)^4 \A^{f(t)}_7  &= 630 f_{12}(t)^2 \beta(t) u^6 + O(u^7)  \qquad \textnormal{using \eqref{pd7_neighborhood_inside_a6_numbers_solution}}\nonumber \\
\implies \us_{\PP \A_7}(\ff(t), \lp(t)) &= 630 \us_{\PP \E_6}(\ff(t), \lp(t))^2 \Bz(\ff(t), \lp(t)) t^6 + O(t^7)
\end{align}    
Since $(\ff, \lp) \in \PP \D_7$, we get that $ \Bz(\ff(t), \lp(t))^2 \neq 0$ and 
$\us_{\PP \E_6}(\ff(t), \p(t)) \neq 0$.  
Hence, \eqref{pa7_around_pd7_again} implies that if $t$ is sufficiently small then $\us_{\PP \A_7}( \ff(t), \lp(t)) \neq 0$, since all the sections are continuous. 
This proves claim \ref{claim_a6_closure_simultaneous_pd6}. \\ 
\hf Before proving \eqref{pa6_cap_pd4_pe7_equal_zero_px8_not_zero_equal_pe7}, let us prove a corollary which will be used in the proof of \eqref{algopa7}. 
\begin{cor}
\label{mult_of_pa7_section_around_pd7}
Let $\W \lra \D \times \P T\P^2$ be a vector bundle such that 
the rank of $\W$ is same as the dimension of $\PP \D_7$. 
Let $\Q: \D \times \P T\P^2 \lra \W$ be a \textit{generic} 
smooth section. 
Suppose $(\ff,\lp) \in \PP \D_7 \cap \Q^{-1}(0)$. 
Then the section $$ \us_{\PP \A_7} \oplus \Q: \ov{\PP \A}_6 \lra \UL_{\PP \A_7} \oplus \W$$
vanishes around $(\ff, \lp)$ with a multiplicity of $6$.
\end{cor}
\pf Follows from the fact that $\Q^{-1}(0)$ intersects $\PP \D_7$ transversely and \eqref{pa7_around_pd7_again}. \qed  \\

\hf Next we will prove \eqref{pa6_cap_pd4_pe7_equal_zero_px8_not_zero_equal_pe7}. To show that the left hand side is a subset of the right hand side note that $\ov{\PP \A}_6 \subset \ov{\PP \A}_5$ by Lemma \ref{cl}, statement \ref{A5cl}. Consequently,
\bgd
\{ (\ff,l_p) \in \ov{\PP \A}_6: \us_{\PP \D_4}(\ff,l_p) = 0, \us_{\PP \E_6}(\ff,l_p) = 0, \us_{\PP \XC_8}(\ff,l_p) \neq 0 \}
\edd
is contained in
\bgd
\{ (\ff,l_p) \in \ov{\PP \A}_5: \us_{\PP \D_4}(\ff,l_p) = 0, \us_{\PP \E_6}(\ff,l_p) = 0, \us_{\PP \XC_8}(\ff,l_p) \neq  0 \}.
\edd
But the quantity above, by \eqref{pa5_cap_pd4_pe6_equal_zero_equal_pe6}, equals 
\bgd
 \{ (\ff,l_p) \in \ov{\PP \E}_6: \us_{\PP \XC_8}(\ff,l_p) \neq 0 \} 
\subset \PP \E_6 \cup \ov{\PP \E}_7 \cup \ov{\hat{\XC}_8^{\#}},
\edd
where the last inclusion follows from Lemma \ref{cl}, statement \ref{E6cl}. Therefore, the left hand side of \eqref{pa6_cap_pd4_pe7_equal_zero_px8_not_zero_equal_pe7} is a subset of its right hand side,  using \eqref{pa6closure_cap_pe6_empty}.\\
\hf For the converse, since $\ov{\PP \A}_6$ is a closed set, it suffices to show that
\begin{align} 
\label{pa6_cap_pd4_pe7_neq_zero_supset_pee7}
\{ (\ff,l_p) \in \ov{\PP \A}_6: \us_{\PP \D_4}(\ff,l_p) = 0, \us_{\PP \E_6}(\ff,l_p) = 0, \us_{\PP \XC_8}(\ff,l_p) \neq 0 \} 
\supset  \{ (\ff,l_p) \in \PP \E_7: \us_{\PP \XC_8}(\ff,l_p) \neq 0  \}.
\end{align}
We will simultaneously prove this statement and also prove 
\bge
\label{pa7closure_cap_pe7_empty}
\ov{\PP \A}_7 \cap \PP \E_7 = \varnothing.
\ede

\begin{claim}
\label{claim_a6_closure_simultaneous_pee7}
Let $(\ff,\lp) \in \PP \E_7$.
Then there exists a solution 
$(\ff(t), \lp(t) ) \in \ov{\PP \A}_3$ near $(\ff,\lp)$ to the set of equations
\begin{align}
\us_{\PP \D_4}( \ff(t), \lp(t) ) & \neq 0,  \us_{\PP \A_4}( \ff(t), \lp(t)) = 0, \us_{\PP \A_5}( \ff(t), \lp(t)) = 0, \us_{\PP \A_6}( \ff(t), \lp(t)) = 0. \label{pe7_neighborhood_inside_a6} 
\end{align}
Moreover, \textit{whenever} such a solution $(\ff(t), \lp(t))$ is sufficiently close to $(\ff,\lp)$ 
it lies in $\PP \A_6$, i.e., $\us_{\PP \A_7}( \ff(t), \lp(t)) \neq 0.$ In particular $(\ff(t), \lp(t))$ \textit{does not} lie in $ \PP \A_7$.
\end{claim}
Note that claim \ref{claim_a6_closure_simultaneous_pee7} proves \eqref{pa6_cap_pd4_pe7_neq_zero_supset_pee7} and \eqref{pa7closure_cap_pe7_empty} simultaneously. \\

\pf As before, we will first solve the equation 
\begin{align}
\label{pe7_neighborhood_inside_a6_numbers}
f_{02}(t) \neq 0, ~~f_{02}(t) \A^{f(t)}_4 = 0, ~~f_{02}(t)^2 \A^{f(t)}_5  = 0, ~~f_{02}(t)^3 \A^{f(t)}_6  = 0.
\end{align}
The only solutions to \eqref{pe7_neighborhood_inside_a6_numbers}, that converge to zero as $f_{02}(t)$ and $f_{12}(t)$ go to zero 
are 
\begin{align}
\label{pe7_neighborhood_inside_a6_numbers_solution}
f_{12}(t) & = u \nonumber \\
f_{21}(t) & = -\frac{3 }{2 f_{03}(t)}u^2 + O(u^3)  \nonumber\\
f_{02}(t) & =  - \frac{9 }{4 f_{31}(t) f_{03}(t)} u^3 + O(u^4)  \nonumber \\ 
f_{40}(t) & = - \frac{3 f_{31}(t) }{f_{03}(t)} u + O(u^2)
\end{align}
We set $f_{12}(t)=u$ and use $f_{02}(t)^2\A^{f(t)}_5=0$ to solve for $f_{02}/f_{21}$; we get
\bgd
\frac{f_{02}(t)}{f_{21}(t)}=\frac{3u}{2f_{31}(t)}+ O(u^2).
\edd 
We now use $f_{02}(t)^3 \A^{f(t)}_6  = 0$ ans replace $f_{02}$ with the expression above to solve for $f_{21}$ and obtain the second equation. This also implies the third equation by the expression above. 
The last equation can now be obtained using $f_{02}(t) \A^{f(t)}_4  = 0$ and the first three equations. \\
\hf Let us explain the method in detail. Using $f_{02}(t)^2 \A^{f(t)}_5 =0$ we can solve for 
$\frac{f_{02}}{f_{21}}$ and get 
\begin{align}
\label{pe7_only_one_branch_of_square_root}
\frac{f_{02}(t)}{f_{21}(t)} = 
\frac{10 f_{31}(t) - \sqrt{ 100 f_{31}(t)^2 - 60 f_{50}(t) u }}{2 f_{50}(t)} 
= \frac{3}{2 f_{31}(t)} u + O(u^2)
\end{align}
It is easy to see that we never really used the fact that $f_{50}(t) \neq 0$; the equality of the 
first and last term remains valid even when $f_{50}=0$. 
However, we do need to justify why we did not choose the other branch of the square root. 
We will explain that shortly.  
Plugging in the value of $f_{02}$ from \eqref{pe7_only_one_branch_of_square_root} in equation $f_{02}(t)^3 \A^{f(t)}_6=0$ 
and by using the Implicit Function Theorem, 
we get the expression for $f_{21}(t)$ in \eqref{pe7_neighborhood_inside_a6_numbers_solution}. 
And now using the 
value of $f_{21}(t)$ and \eqref{pe7_only_one_branch_of_square_root} we get the expression for $f_{02}(t)$ 
in \eqref{pe7_neighborhood_inside_a6_numbers_solution}. \\
\hf It remains to show that why we did not chose 
the other branch of the square root. It is easy to see that if we chose the 
other branch, it would imply that as $f_{02}(t)$  and $f_{21}(t)$ go to zero, the 
ratio $L_t := \frac{f_{21}(t)}{f_{02}(t)} $ tends to a finite number 
$L$, since $f_{31} \neq 0$. Using $f_{03}(t)^3 \A^{f(t)}_6 =0$ we can solve 
for $f_{31}(t)$ as a quadratic equation and get that  
\begin{align*}
f_{31}(t)= \frac{30 L_t f_{12}(t) \pm 
\sqrt{10} \sqrt{-15 L_t^3 f_{02}(t) f_{03}(t) + 
45 L_t^2 f_{02}(t) f_{22}(t) -15 L_t f_{02}(t) f_{41}(t) + f_{02}(t) f_{60}(t)}}
{10}.
\end{align*}
It is now easy to see that $f_{31}(t)$ tends to zero as $f_{12}(t)$ and $f_{02}(t)$ tend to zero. 
This gives us a contradiction, since $f_{31} \neq 0$.\\
\hf Since $(\ff, \lp) \in \PP \E_7$ we get that $ f_{03}, ~f_{31} \neq 0 $. 
Equation \eqref{pe7_neighborhood_inside_a6_numbers_solution} implies that the only solutions to the functional equation \eqref{pe7_neighborhood_inside_a6} 
are of the form
\begin{align}
\label{pe7_neighborhood_inside_a6_functional_solution}
\us_{\PP \E_6}( \ff(t), \lp(t)) &= t  \nonumber \\
\us_{\PP \D_5}^{\UL}( \ff(t), \lp(t)) &= -\frac{3 }{2 \us_{\PP \XC_8}( \ff(t), \lp(t)) } t^2 + O(t^3) \nonumber \\
\us_{\PP \D_4}( \ff(t), \lp(t)) &= - \frac{9 }{4 \us_{\PP \E_8}( \ff(t), \lp(t))  \us_{\PP \XC_8}( \ff(t), \lp(t))} t^3 + O(t^4) \nonumber \\ 
\us_{\PP \E_7}( \ff(t), \lp(t)) &=  - \frac{3 \us_{\PP \E_8}( \ff(t), \lp(t))}{\us_{\PP \XC_8}( \ff(t), \lp(t))} t + O(t^2)   
\end{align}
where equality holds as \textit{functionals}. Since the sections 
\begin{align*}
\us_{\PP \D_4} &: \ov{\PP \A}_3 \lra \UL_{\PP \D_4},  \qquad \us_{\PP \D_5}^{\UL}: \us_{\PP \D_4}^{-1}(0) \lra \UL_{\PP \D_5},\\ 
\us_{\PP \E_6} &: \us_{\PP \D_5}^{\UL^{-1}}(0) \lra \UL_{\PP \E_6},  \qquad  \us_{\PP \E_7}: \us_{\PP \E_6}^{-1}(0)  \lra \UL_{\PP \E_7}
\end{align*}
are transverse to the zero set (as proved in Proposition \ref{D4_Condition_prp}, \ref{D5_Condition_prp} and \ref{E6_Condition_prp}), there exists 
a solution $(\ff(t), \lp(t))$ close to $(\ff, \lp)$ for \eqref{pe7_neighborhood_inside_a6_functional_solution}. This proves our first assertion. Next we need to show that any such solution satisfies the condition 
$\us_{\PP \A_7}( \ff(t), \lp(t)) \neq 0$ if $t$ is sufficiently small. To prove that we observe  
\begin{align}
\label{pa7_around_pe7_again}
f_{02}(t)^4 \A^{f(t)}_7  &= -\frac{2835}{16 f_{03}(t)^2} u^7 + O(u^8)  \qquad \textnormal{using \eqref{pe7_neighborhood_inside_a6_numbers_solution}}\nonumber \\
\implies \us_{\PP \A_7}(\ff(t), \lp(t)) &= -\frac{2835}{16 \us_{\PP \XC_8}( \ff(t), \lp(t))^2 } t^7 + O(t^8)
\end{align}    
Since $(\ff, \lp) \in \PP \E_7$, we get that $ \us_{\PP \XC_8}( \ff, \lp) \neq 0$.  
Hence, \eqref{pa7_around_pe7_again} implies that if $t$ is sufficiently small $\us_{\PP \A_7}( \ff(t), \lp(t)) \neq 0$, since all the sections are continuous. This proves claim \ref{claim_a6_closure_simultaneous_pee7}. \\ 
\hf Before proving \eqref{pa6_cap_pd4_pe7_equal_zero_px8_equal_zero_px8}, let us prove a corollary which will be used in the proof of \eqref{algopa7}.
\begin{cor}
\label{mult_of_pa7_section_around_pe7}
Let $\W \lra \D \times \P T\P^2$ be a vector bundle such that 
the rank of $\W$ is same as the dimension of $\PP \E_7$. 
Let $\Q: \D \times \P T\P^2 \lra \W$ be a \textit{generic} 
smooth section. 
Suppose $(\ff,\lp) \in \PP \E_7 \cap \Q^{-1}(0)$. 
Then the section $$ \us_{\PP \A_7} \oplus \Q: \ov{\PP \A}_6 \lra \UL_{\PP \A_7} \oplus \W$$
vanishes around $(\ff, \lp)$ with a multiplicity of $7$.
\end{cor}
\pf Follows from the fact that $\Q^{-1}(0)$ intersects $\PP \E_7$ transversely and \eqref{pa7_around_pe7_again}. \qed  \\

\hf Finally, we will prove \eqref{pa6_cap_pd4_pe7_equal_zero_px8_equal_zero_px8}. Let us show that the left hand side is contained in the right hand side. Note that $\ov{\PP \A}_6  \subset \ov{\PP \A}_5$ by Lemma \ref{cl}, statement \ref{A5cl}. Therefore,
\bgd
\{ (\ff,\lp) \in \ov{\PP \A}_6: \us_{\PP \D_4}(\ff,\lp) = 0, \us_{\PP \E_6}(\ff, \lp) = 0, \us_{\PP \XC_8}(\ff, \lp) = 0 \}
\edd
is contained in
\bgd
\{ (\ff, \lp) \in \ov{\PP \A}_5: \us_{\PP \D_4}(\ff, \lp) = 0,\us_{\PP \E_6}(\ff, \lp) = 0, \us_{\PP \XC_8}(\ff, \lp) = 0\}.
\edd
The last quantity above equals, due to \eqref{pa5_cap_pd4_pe6_equal_zero_equal_pe6}, and \eqref{x8_in_closure_of_e6_equation} and Corollary \ref{X8_and_X8_sharp_up_cl_vanish_cor}
\bgd
\textup{LHS}=\{ (\ff, \lp) \in \ov{\PP \E}_6: \us_{\PP \XC_8}(\ff, \lp) = 0 \}=  \ov{\hat{\XC}_8^{\# \flat}}.
\edd
This implies that the left hand side of \eqref{pa6_cap_pd4_pe7_equal_zero_px8_equal_zero_px8} is a subset of its right hand side.\\
\hf Next let us show that the right hand side of \eqref{pa6_cap_pd4_pe7_equal_zero_px8_equal_zero_px8} is a subset of its left hand side. 
Since $\ov{\PP \A}_6$ is a closed set, it suffices to show that 
\begin{align} 
\{ (\ff,l_p) \in \ov{\PP \A}_6: \us_{\PP \D_4}(\ff,l_p) = 0, ~~\us_{\PP \E_7}(\ff,l_p) = 0, ~~\us_{\PP \XC_8}(\ff,l_p) = 0 \} \supset 
\XX^{\# \flat }_8 \label{pa6_cap_pd4_pe7_equal_zero_px8_equal_zero_supset_px8}
\end{align}
We will simultaneously prove this statement and also prove 
\begin{align}
\label{pa7closure_cap_x8_empty}
\ov{\PP \A}_7 \cap \XX^{\#\flat }_8 &= \varnothing
\end{align}

\begin{claim}
\label{claim_a6_closure_simultaneous_x8}
Let $(\ff,\lp) \in \XX^{\# \flat}_8$.
Then there exists a solution 
$(\ff(t), \lp(t) ) \in \ov{\PP \A}_3$ near $(\ff,\lp)$ to the set of equations
\begin{align}
\us_{\PP \D_4}( \ff(t), \lp(t) ) & \neq 0, \us_{\PP \A_4}( \ff(t), \lp(t)) = 0, \us_{\PP \A_5}( \ff(t), \lp(t)) = 0, \us_{\PP \A_6}( \ff(t), \lp(t)) = 0. \label{x8_neighborhood_inside_a6} 
\end{align}
Moreover, \textit{whenever} such a solution $(\ff(t), \lp(t))$ is sufficiently close to $(\ff ,\lp)$ 
it lies in $\PP \A_6$, i.e., $\us_{\PP \A_7}( \ff(t), \lp(t)) \neq 0$. In particular, $(\ff(t), \lp(t))$ \textit{does not} lie in $ \PP \A_7$.
\end{claim}
Notice that claim \ref{claim_a6_closure_simultaneous_x8} proves \eqref{pa6_cap_pd4_pe7_equal_zero_px8_equal_zero_supset_px8} and \eqref{pa7closure_cap_x8_empty} simultaneously.  \\

\pf As before, we will first solve the equation 
\begin{align}
\label{x8_neighborhood_inside_a6_numbers}
f_{02}(t) \neq 0, ~~f_{02}(t) \A^{f(t)}_4 = 0, ~~f_{02}(t)^2 \A^{f(t)}_5  = 0, ~~f_{02}(t)^3 \A^{f(t)}_6  = 0.
\end{align}
The only solutions to \eqref{x8_neighborhood_inside_a6_numbers} that converge to zero as $f_{02}(t)$, $f_{12}(t)$ and $f_{03}$ go to zero  
are
\begin{align}
\label{x8_neighborhood_inside_a6_numbers_solution}
f_{21}(t)  & = u \nonumber   \\
f_{02}(t)  & =  \frac{3u^2}{f_{40}} \qquad \textnormal{using} \qquad  f_{02}(t) \A^{f(t)}_4 = 0 \nonumber     \\
f_{12}(t)  & = \frac{2 f_{31}}{f_{40}} u - \frac{3 f_{50}}{5 f_{40}^2} u^2 \qquad \textnormal{using} \qquad f_{02}(t)^2 \A^{f(t)}_5  = 0 \nonumber   \\ 
f_{03}(t) & = \Big( - \frac{6 f_{31}^2}{f_{40}^2} + \frac{9 f_{22}}{f_{40}} \Big) u + 
                                                                           \Big( - \frac{9 f_{41}}{f_{40}^2} + \frac{36 f_{31} f_{50}}{5 f_{40}^3} \Big)u^2 \nonumber  \\  
                 & + \Big( - \frac{54 f_{50}^2}{25 f_{40}^4} + \frac{9 f_{60}}{5 f_{40}^3} \Big) u^3 \qquad \textnormal{using} \qquad f_{02}(t)^3 \A^{f(t)}_6  =0.
\end{align}
Note that since $(\ff, \lp) \in \XX^{\# \flat}_8$ we get that $f_{40} \neq 0$. Equation \eqref{x8_neighborhood_inside_a6_numbers_solution} implies that the only solutions to the functional equation \eqref{x8_neighborhood_inside_a6} are of the form 
\begin{align}
\label{x8_neighborhood_inside_a6_functional_solution}
\us_{\PP \D_5}^{\UL}( \ff(t), \lp(t)) &= t  \nonumber \\
\us_{\PP \D_4}( \ff(t), \lp(t)) &=  \frac{3}{ \us_{\PP \D_6}( \ff(t), \lp(t)) } t^2  \nonumber \\
\us_{\PP \E_6}( \ff(t), \lp(t)) &=  \frac{2 \us_{\PP \E_8}( \ff(t), \lp(t))}{\us_{\PP \D_6}( \ff(t), \lp(t))} t  + O(t^2) 
\nonumber \\ 
\us_{\PP \XC_8}( \ff(t), \lp(t)) &= \Big( - \frac{6 \us_{\PP \E_8}( \ff(t), \lp(t))^2}{ \us_{\PP \D_6}( \ff(t), \lp(t))^2} + \frac{9\, \Si(\ff(t), \lp(t))}{ \us_{\PP \E_7}( \ff(t), \lp(t)) } \Big) t + O(t^2).
\end{align}
The functional $\Si$ is given by
\bgd
\{\Si(\ff(t), \lp(t))\}( f \otimes p^{\otimes d} \otimes v^{\otimes 2} \otimes \w^{\otimes 2}) := f_{22},
\edd
where notations are as defined in subsection \ref{summary_sections_of_vector_bundle_definitions}. Equality holds here as \textit{functionals}. Since the sections 
\begin{align*}
\us_{\PP \D_4} &: \ov{\PP \A}_3 \lra \UL_{\PP \D_4},  \qquad \us_{\PP \D_5}^{\UL}: \us_{\PP \D_4}^{-1}(0) \lra \UL_{\PP \D_5},\\ 
\us_{\PP \E_6} &: \us_{\PP \D_5}^{\UL^{-1}}(0) \lra \UL_{\PP \E_6},  \qquad  \us_{\PP \XC_8}: \us_{\PP \E_6}^{-1}(0) \lra \UL_{\PP \XC_8}
\end{align*}
are transverse to the zero set (as proved in Proposition \ref{D4_Condition_prp}, \ref{D5_Condition_prp} and \ref{E6_Condition_prp}), there exists 
a solution $(\ff(t), \lp(t))$ close to $(\ff, \lp)$ to \eqref{x8_neighborhood_inside_a6_functional_solution}. This proves our first assertion. Next we need to show that any such solution satisfies the condition 
$\us_{\PP \A_7}( \ff(t), \lp(t)) \neq 0$ if $t$ is small. Observe  
\begin{align}
\label{x8_around_pd7_again}
f_{02}(t)^4 \A^{f(t)}_7  &= \Big(- \frac{f_{31}(t)^3}{8 f_{40}(t)^3} + \frac{3 f_{22}(t) f_{31}(t)}{16 f_{40}(t)^2} - \frac{f_{13}(t)}{16 f_{40}(t)}  \Big) u^5 + O(u^6) 
\qquad \textnormal{using \eqref{x8_neighborhood_inside_a6_numbers_solution}}\nonumber \\
\implies \us_{\PP \A_7}(\ff(t), \lp(t)) &= \frac{ \us_{\PP \J} (\ff(t), \lp(t))}{ \us_{\PP \E_7} (\ff(t), \lp(t))^3 } t^5 + O(t^6) 
\end{align}    
Since $(\ff, \lp) \in \XC_{8}^{\# \flat}$, we get that $\us_{\PP \J} (\ff, \lp)\neq 0$ and $\us_{\PP \E_7} (\ff, \lp) \neq 0$. Hence, \eqref{x8_around_pd7_again} implies that if $t$ is small $\us_{\PP \A_7}( \ff(t), \lp(t)) \neq 0$, which proves claim \ref{claim_a6_closure_simultaneous_x8}. This finishes the proof of Lemma \ref{cl}, statement \ref{A6cl}. \qed


\section{Euler class}  
\label{Euler_class_computation}

\hf\hf Finally, we are ready to prove the recursive formulas stated in section \ref{algorithm_for_numbers}. 
The notations are as in section \ref{summary_notation_def} and notations \ref{tau_bundle_defn} and \ref{tau_bundle_pv_defn}. \\

\ni \textbf{Proof of Equation \eqref{algoa1}:} Let  $\Q: \D \times \P^2 \lra \WL$ be a generic smooth section to
\bgd
\WL := \bigg(\bigoplus_{i=1}^{\delta_d -(n+1)}\gD^* 
\bigg)\oplus\bigg(\bigoplus_{i=1}^{n} \gP^*\bigg) \lra \D \times \P^2.
\edd
By Lemma \ref{gpl} and Theorem \ref{Main_Theorem_pseudo_cycle} we conclude 
\begin{align*}
\Num(\A_1, n) = \lan e(\WL), ~[\ov{\A}_1] \ran =  \pm |\A_1 \cap  \Q^{-1}(0)|.
\end{align*}
By Lemma \ref{cl}, statement \ref{A0cl},  
$~\ov{\A}_0 = \A_0 \cup \ov{\A}_1$.  
The section $\ds_{\A_1} : \ov{\A}_0 \lra \DV_{\A_1}$ 
vanishes on $\A_1$  transversely and doesn't vanish on $\A_0$ 
(cf. Proposition \ref{ift_ml}). Therefore, the zeros of the section 
$$ \ds_{\A_1} \oplus \Q : \ov{\A}_0 \lra \DV_{\A_1} \oplus \WL $$
counted with a sign is our desired number, whence
\bgd
\Num(\A_1, n) = \lan e(\DV_{\A_1}) e(\WL), [\ov{\A}_0] \ran =  \lan \PD [\ov{\A}_0] e(\DV_{\A_1})  e(\WL), [\D \times \P^2] \ran.
\edd
Now Proposition \ref{a0_cl_vanish_cor}, \ref{ift_ml} and Theorem \ref{Main_Theorem} imply that the Poincar\'{e} dual $\PD[\ov{\A}_0]$ of $\ov{\A}_0$ is the Euler class $e(\DL_{\A_0})$. We may now use the splitting principle and Lemma \ref{total_chern_class_of_tpn} to conclude that
\bgd
\Num(\A_1, n) = \lan (\y + d \a)( (\y + d \a)^2 -3\a(\y+ d \a) + 3 \a^2 ) \y^ {\delta_d -(n+1)} \a^{n},  [\D \times \P^2] \ran.
\edd
Equation \eqref{algoa1} now follows.  $\hfill\square$\\

\ni \textbf{Proof of Equation \eqref{algopa20} and \eqref{algopa21}:} Let $\W_{n,m,2}$ and $\Q$ be as in \eqref{generic_Q} with $k=2$. By definition, $~\Num(\PP\A_2, n,m)$ is the signed 
cardinality of the intersection of $\PP \A_2$ with $\Q^{-1}(0)$. By Lemma \ref{cl}, statement \ref{A1cl} we gather that
\bgd
\ov{\hat{\A}_1} = \ov{\hat{\A}^{\#}_1} = \hat{\A}_1^{\#} \cup \ov{\PP \A}_2.
\edd
By Proposition \ref{A1_sharp_Condition_prp}, the section $\us_{\PP \A_2} : \ov{\hat{\A}_1} \lra \UV_{\PP \A_2}$ vanishes on $\PP \A_2$ transversely and by \textit{definition} it  
doesn't vanish on $\hat{\A}_1^{\#}$. Hence, the zeros of the section 
$$ \us_{\PP \A_2} \oplus \Q : \ov{\hat{\A}_1} \lra \UV_{\PP \A_2} \oplus \W_{n,m,2}, $$
counted with a sign, is our desired number. Via the splitting principle and Lemma \ref{total_chern_class_of_tpn} we have
\BAA
\Num(\PP\A_2, n,m) & = & \lan e(\UV_{\PP \A_2}) e(\W_{n,m,2}), [\ov{\hat{\A}_1}] \ran \\
                   & = & \lan ( (\lm+ \y + d \a)^2 -3\a(\lm+ \y+ d \a) + 3 \a^2) \y^{\delta_d-(n+m+2)}\a^n \lm^m, [\ov{\hat{\A}_1}]  \ran.
\EAA
\ni  Next we use the fact that 
\begin{align*}
\lan \pi^*(\y^{\delta_d-(n_1+1)} \a^{n_1}) \lm, [\ov{\hat{\A}_1}] \ran  = \lan \y^{\delta_d-(n_1+1)} \a^{n_1}, [\ov{\A}_1] \ran  
\qquad\textnormal{and} \qquad \lan \pi^*(\y^{\delta_d-n_1} \a^{n_1}) , [\ov{\hat{\A}_1}] \ran & = 0  
\end{align*}
for all $n_1$. This follows from Lemma \ref{cohomology_ring_of_pv} and the fact that $\ov{\A}_1$ is a smooth manifold (cf .Corollary \ref{a1_cl_vanish_cor}). 
Finally, using the ring structure of $H^*(\D \times \P T\P^2; \mathbb{Z})$ (cf. Lemma \ref{cohomology_ring_of_pv}), we obtain equations \eqref{algopa20} and \eqref{algopa21}. Here $\pi: \D \times \P T\P^2 \lra \D \times \P^2 $ is the projection map. \qed \\

\ni \textbf{Proof of Equation \eqref{algopa3}:} Let $\W_{n,m,3}$ and $\Q$ be as in \eqref{generic_Q} with $k=3$.
By Lemma \ref{cl}, statement \ref{A2cl} we have 
$$\ov{\PP \A}_2 = \PP \A_2 \cup 
\ov{\PP \A}_3 \cup \ov{\hat{\D}_4^{\#}}.$$
The section 
$ \us_{\PP \A_3} : \ov{\PP \A}_2 \lra \UL_{\PP \A_3}$
doesn't vanish on $\PP \A_2$ and 
vanishes transversely on $\ov{\PP \A}_3$. Furthermore, it does not 
vanish on any point of $\hat{\D}_4^{\#}$ (by \textit{definition}). 
Hence, the zeros of the section 
$$ \us_{\PP \A_3} \oplus \Q : \ov{\PP \A}_2 \lra \UL_{\PP \A_3} \oplus \W_{n,m,3} $$
counted with a sign is $\Num(\PP\A_3, n, m)$. A similar computation using the product formula for the 
first Chern class of a product of line bundles, proves the equation. \qed \\ 

\ni \textbf{Proof of Equation \eqref{algopa4}:} Let $\W_{n,m,4}$ and $\Q$ be as in \eqref{generic_Q} with $k=4$. By Lemma \ref{cl}, statement \ref{A3cl} we have that 
$$\ov{\PP \A}_3 = \PP \A_3 \cup \ov{\PP \A}_4 \cup \ov{\PP \D}_4.$$
The section $ \us_{\PP \A_4} : \ov{\PP \A}_3 \lra \UL_{\PP \A_4}$
doesn't vanish on $\PP \A_3$ and vanishes transversely on $\PP \A_4$
(cf. Proposition \ref{Ak_Condition_prp}). Furthermore, it does not 
vanish on any point of $\PP \D_4$. Hence, the zeros of the section 
$$ \us_{\PP \A_4} \oplus \Q : \ov{\PP \A}_3 \lra \UL_{\PP \A_4} \oplus \W_{n,m,4} $$
counted with a sign is $\Num(\PP\A_4, n, m)$, which proves the equation. \qed \\

\ni \textbf{Proof of Equation \eqref{algopa5}:} Let $\W_{n,m,5}$ and $\Q$ be as in \eqref{generic_Q} with $k=5$. By Lemma \ref{cl}, statement \ref{A4cl} we have that 
$$\ov{\PP \A}_4 = \PP \A_4 \cup \ov{\PP \A}_5 \cup \ov{\PP \D}_5.$$
The section $ \us_{\PP \A_5} : \ov{\PP \A}_4 \lra \UL_{\PP \A_5}$
doesn't vanish on $\PP \A_4$ and vanishes transversely on $\PP \A_5$ (see Proposition \ref{Ak_Condition_prp}). Furthermore, the section 
$$\us_{\PP \A_5} \oplus \Q: \ov{\PP \A}_4 \lra \UL_{\PP \A_5} \oplus \W_{n,m,5} $$ 
vanishes on $\PP \D_5$ with a multiplicity of $2$ (cf. Corollary \ref{mult_of_pa5_section_around_pd5}). Hence,
\begin{align*}
\lan e(\UL_{\PP \A_5}) e(\W_{n,m,5}), [\ov{\PP \A}_4]  \ran = \Num(\PP \A_5,n, m)+ 2\Num(\PP \D_5, n, m) 
\end{align*}
completing the proof. \qed\\

\ni \textbf{Proof of Equation \eqref{algopa6}:} Let $\W_{n,m,6}$ and $\Q$ be as in \eqref{generic_Q} with $k=6$. By Lemma \ref{cl}, statement \ref{A5cl} we have that 
$$\ov{\PP \A}_5 = \PP \A_5 \cup \ov{\PP \A}_6 \cup \ov{\PP \D}_6 \cup \ov{\PP \E}_6.$$
The section $ \us_{\PP \A_6} : \ov{\PP \A}_5 \lra \UL_{\PP \A_6}$
doesn't vanish on $\PP \A_5$ and vanishes transversely on $\PP \A_6$.
Furthermore, the section 
$$\us_{\PP \A_6} \oplus \Q: \ov{\PP \A}_5 \lra \UL_{\PP \A_6} \oplus \W_{n,m,6} $$ 
vanishes on $\PP \D_6$ and $\PP \E_6$ with a multiplicity of $4$ and $3$ respectively (cf. Corollary \ref{mult_of_pa6_section_around_pd6} and \ref{mult_of_pa6_section_around_pe6}).  \qed  \\

\ni \textbf{Proof of Equation \eqref{algopa7}:} Let $\W_{n,0,7}$ and $\Q$ be as in \eqref{generic_Q} with $m=0$ and $k=7$. By Lemma \ref{cl}, statement \ref{A5cl} we have that 
$$\ov{\PP \A}_6 = \PP \A_6 \cup \ov{\PP \A}_7 \cup \ov{\PP \D}_7 \cup \ov{\PP \E}_7 
\cup \ov{\hat{\XC}_8^{\# \flat}}.$$
The section $ \us_{\PP \A_7} : \ov{\PP \A}_6 \lra \UL_{\PP \A_7}$
doesn't vanish on $\PP \A_6$ and vanishes transversely on $\PP \A_7$.
Furthermore, the section 
$$\us_{\PP \A_7} \oplus \Q: \ov{\PP \A}_6 \lra \UL_{\PP \A_7} \oplus \W_{n,m,7} $$ 
vanishes on $\PP \D_7$ and $\PP \E_7$ with a multiplicity of $6$ and $7$ respectively (cf. Corollary \ref{mult_of_pa6_section_around_pd6} and \ref{mult_of_pa6_section_around_pe6}).
Let us assume the section vanishes with a multiplicity of $\eta$ on $\hat{\XC}_8^{\# \flat}$.  
Hence, 
\begin{align*}
\lan e(\UL_{\PP \A_7}) e(\W_{n,0,7}), 
[\ov{\PP \A}_6]  \ran &= \Num(\PP \A_7,n, 0)+ 
6 \Num(\PP \D_7, n, 0) + 7 \Num(\PP \E_7, n,0) 
+ \eta \lan e(\W_{n,0,7}), [\ov{\hat{\XC}_8^{\#}}] \ran   
\end{align*}
It is easy to see that $ \lan e(\W_{n,0,7}), [\ov{\hat{\XC}_8^{\# \flat}}] \ran =0$, which proves the equation. \qed \\

\ni \textbf{Proof of Equation \eqref{algopd4}:}  Let $\W_{n,m,4}$ and $\Q$ be as in \eqref{generic_Q} with $k=4$. By Lemma \ref{cl}, statement \ref{A3cl} we have 
$$\ov{\PP \A}_3 = \PP \A_3 \cup \ov{\PP \A}_4 \cup \ov{\PP \D}_4.$$
The section $ \us_{\PP \D_4} : \ov{\PP \A}_3 \lra \UL_{\PP \D_4}$ doesn't vanish on $\PP \A_3$ and vanishes transversely on $\PP \D_4$. Furthermore, this section does not vanish on any point of $\PP \A_4$. Hence, the zeros of the section 
$$ \us_{\PP \D_4} \oplus \Q : \ov{\PP \A}_3 \lra \UL_{\PP \D_4} \oplus \W_{n,m,4} $$
counted with a sign is $\Num(\PP \D_4, n,m)$ which proves the equation. \qed \\

\ni \textbf{Proof of Equation \eqref{algopd5}:} Let $\W_{n,m,5}$ and $\Q$ be as in \eqref{generic_Q} with $k=5$. By Lemma \ref{cl}, statement \ref{D4cl} we have that 
$$\ov{\PP \D}_4 = \PP \D_4 \cup \ov{\PP \D}_5 \cup \ov{\PP \D_5^{\vee}}.$$
The section $ \us_{\PP \D_5} : \ov{\PP \D}_4 \lra \UL_{\PP \D_5}$
doesn't vanish on $\PP \D_4$ and vanishes transversely on $\PP \D_5$. Moreover, the section does not vanish on $\PP \D_5^{\vee}$ by definition.
Hence, the zeros of the section 
$$ \us_{\PP \D_5} \oplus \Q : \ov{\PP \D}_4 \lra \UL_{\PP \D_5} \oplus \W_{n,m,5} $$
counted with a sign is $\Num(\PP \D_5, n,m)$, which proves the equation. \\ 
\hf \hf Here is an alternative way to compute $\Num(\D_5, n)$. First  note that the map $\pi: \PP\D_5^{\vee} \lra \D_5$ is one to one. 
Hence $\lan \W_{n,0,5}, ~ [\ov{\PP \D^{\vee}_5}]\ran = \Num(\D_5, n)$. Now notice that the section 
$\us_{\PP \D_5^{\vee}}: \ov{\PP \D}_4 \lra \UL_{\PP \D_5^{\vee}}$ vanishes transversely on $\PP \D_5^{\vee}$ and 
does not vanish on $\PP \D_5$. Hence 
\begin{align*}
\lan e(\UL_{\PP \D_5^{\vee}} \oplus \W_{n,0,5}), ~[\ov{\PP \D}_4] \ran &= \lan \W_{n,0,5}, ~ [\ov{\PP \D^{\vee}_5}]\ran =  \Num(\D_5, n).  
\end{align*}
It is easy to check directly that these two methods give the same answer for $\Num(\D_5, n)$. \qed \\ 


\ni \textbf{Proof of Equation \eqref{algopd6}:} Let $\W_{n,m,6}$ and $\Q$ be as in \eqref{generic_Q} with $k=6$. By Lemma \ref{cl}, statement \ref{D5cl} we have  
$$\ov{\PP \D}_5 = \PP \D_5 \cup \ov{\PP \D}_6 \cup \ov{\PP \E}_6.$$
The section $ \us_{\PP \D_6} : \ov{\PP \D}_5 \lra \UL_{\PP \D_6}$
doesn't vanish on $\PP \D_5$ and vanishes transversely on $\PP \D_6$. 
Furthermore, it does not vanish on any point of $\PP \E_6$. Hence, the zeros of the section 
$$ \us_{\PP \D_6} \oplus \Q : \ov{\PP \D}_5 \lra \UL_{\PP \D_6} \oplus \W_{n,m,6} $$
counted with a sign is $\Num(\PP \D_6, n,m)$, which proves the equation. \qed \\

\ni \textbf{Proof of Equation \eqref{algopd7}:} Let $\W_{n,m,7}$ and $\Q$ be as in \eqref{generic_Q} with $k=7$. By Lemma \ref{cl}, statement \ref{D6cl} we have that 
$$\ov{\PP \D}_6 = \PP \D_6 \cup \ov{\PP \D}_7 \cup \ov{\PP \E}_7.$$
The section $ \us_{\PP \D_7} : \ov{\PP \D}_6 \lra \UL_{\PP \D_7}$
doesn't vanish on $\PP \D_6$ and vanishes transversely on $\PP \D_7$.
Furthermore, it does not vanish on any point of $\PP \E_7$. Hence, the zeros of the section 
$$ \us_{\PP \D_7} \oplus \Q : \ov{\PP \D}_6 \lra \UL_{\PP \D_7} \oplus \W_{n,m,7} $$
counted with a sign is $\Num(\PP \D_7, n,m)$, which proves the equation. \qed  \\ 

\ni \textbf{Proof of Equation \eqref{algope6}:} Let  $\W_{n,m,6}$ and $\Q$ be as in \eqref{generic_Q} with $k=6$. By Lemma \ref{cl}, statement \ref{D5cl} we have that 
$$\ov{\PP \D}_5 = \PP \D_5 \cup \ov{\PP \D}_6 \cup \ov{\PP \E}_6.$$
The section $ \us_{\PP \E_6} : \ov{\PP \D}_5 \lra \UL_{\PP \E_6}$
doesn't vanish on $\PP \D_5$ and vanishes transversely on $\PP \E_6$.
Hence, the zeros of the section 
$$ \us_{\PP \E_6} \oplus \Q : \ov{\PP \D}_5 \lra \UL_{\PP \E_6} \oplus \W_{n,m,6} $$
counted with a sign is $\Num(\PP \E_6, n,m)$, which proves the equation. \qed \\ 

\ni \textbf{Proof of Equation \eqref{algope7}:} Let $\W_{n,m,7}$ and $\Q$ be as in \eqref{generic_Q} with $k=7$. By Lemma \ref{cl}, statement \ref{D6cl} we have that 
$$ \ov{\PP \D}_6 = \PP \D_6 \cup \ov{\PP \E}_7 \cup \ov{\PP \D}_7.$$
The section $ \us_{\PP \E_6} : \ov{\PP \D}_6 \lra \UL_{\PP \E_6}$
doesn't vanish on $\PP \D_7$ and vanishes transversely on $\PP \E_7$ (cf .Proposition \ref{E7_Condition_prp_using_D6}).  Hence, the zeros of the section 
$$ \us_{\PP \E_6} \oplus \Q : \ov{\PP \D}_6 \lra \UL_{\PP \E_6} \oplus \W_{n,m,7} $$
counted with a sign is $\Num(\PP \E_7, n,m)$,, which proves the equation. \qed 
 
\appendix

\section{Low degree checks} \label{low_degree_check}
\ni \textbf{Verification of the number $\Num(\A_1,0) = 3(d-1)^2$:} \\
\ni $d=1:$ There are no nodal lines. \\
\ni $d=2:$ The number of line pairs that pass through $4$ general points
is $\frac{1}{2}\binom{4}{2}\!=\!3$. \\
\ni $d = 3:$ The number of nodal cubics passing through $8$ general points are 
the rational cubics passing through these points; this number $12$ can also be computed through Kontsevich's recursion formula. \\  

\ni \textbf{Verification of the number $\Num(\A_1,1) = 3(d-1)$:} \\
\ni $d=1:$ There are no nodal lines. \\
\ni $d=2:$ The number of line pairs that pass through $3$ 
points
and meet on a line is $\binom{3}{2} = 3$. \\

\ni \textbf{Verification of the number $\Num(\A_2,0) = 12(d-1)(d-2)$:} \\
\ni $d=1:$ There are no lines with a cusp. \\
\ni $d=2:$ The only way a conic can have a cusp is if its a double line.
There are no double lines through three generic points. \\
\ni $d=3:$ The number of nodal cubics passing through $7$ general points are 
the rational cubics passing through these points; 
according to \cite{Rahul1} this number $24$. 
This number can also be computed by the algorithm described in \cite{g2p2and3}. \\
\ni $d = 4:$ The number of quartics with a cusp is $72$. This is same as the number of genus two curves with a cusp and equals $72$ (cf. \cite{Va}, pp. 19).\\

\ni \textbf{Verification of the number $\Num(\PP \A_3,n,m)$:} \\ 
\ni $d=3:$ The number $\Num(\PP \A_3,n,m)$ can be verified  by direct geometric means for all values 
of $n$ and $m$ in the case of cubics (cf. \cite{BM_Detail}). \\

\ni \textbf{Verification of the number $\Num(\A_4,0) = 60(d-3)(3d-5)$:} \\
\ni $d = 3:$ There are no cubics with an $\A_4$-node. \\

\ni \textbf{Verification of the number $\Num(\D_4,0) = 15(d-2)^2$:} \\ 
\ni $d = 2$: There are no conics with a $\D_4$-node. \\
\ni $d = 3$: The only way a cubic can have a $\D_4$-node is, if it breaks into three distinct lines intersecting at a common point. The number of such configurations passing 
through $5$ points is $\frac{1}{3}\times \binom{5}{2}\times \binom{3}{2} = 15.$ \\

\ni \textbf{Verification of the number $\Num(\D_4,1) = 6(d-2)$:} \\
\ni $d = 2$: There are no conics with a $\D_4$-node on a line. \\
\ni $d = 3$: The number of triple lines, having a common point
at a given line and passing through four points is $\binom{4}{2} = 6.$ \\

\ni \textbf{Verification of the number $\Num(\PP \D_4,n,1)$:} \\ 
\ni $d=3:$ The number $\Num(\PP \D_4,n,1)$ can be verified  by direct geometric means for all values 
of $n$ in the case of cubics (cf. \cite{BM_Detail}). \\

\ni \textbf{Verification of the number $\Num(\PP \D_6,n,m)$:} \\ 
\ni $d=4:$ The number $\Num(\PP \D_6,n,m)$ can be verified  by direct geometric means for all values 
of $n$ and $m$ in the case of quartics (cf. \cite{BM_Detail}). \\

\ni \textbf{Verification of the number $\Num(\E_6,0) = 21(d-3)(4d-9)$:} \\
\ni $d = 3:$ There are no cubics with an $\E_6$-node. \\
\ni $d=4:$ An $\E_6$-node contributes three to the genus of a 
curve. Since a smooth quartic has genus three, the quartics 
with an $\E_6$-node have genus zero. The number of such quartics 
through $8$ points is $147$ (cf. \cite{g3}, pp. 24).

\section{Some details}

\subsection{Proof of some lemmas}
\label{appendix_proof_of_Lemmas}
\ni \textbf{Proof of Lemma \ref{up_to_down}:} Let $\X_k = \A_k$. If $(\ff, \lp) \in \PP \A_k$ then the projection map has to be one to one, because 
otherwise the kernel of the Hessian would have two linearly independent vectors (which 
would imply it is identically zero). Similar argument holds if $\X_k = \D_k$ for $k \geq 5$ or 
$\E_6, \E_7$ or $\E_8$. In all of the above cases, the map $\PP \X_k \longrightarrow \X_k$ is a diffeomorphism (actually, a biholomorphism). Since $\X_k$ and $\pi(\PP \X_k)$ are cobordant, the pseudocycles $[\ov{\X_k}]$ and $\pi_\ast [\ov {\PP \X_k}]$, defined by $\X_k$ and $\pi(\PP \X_k)$ respectively, define the same homology class. This proves 
\eqref{up_down_equation} (using definition \ref{up_number_defn} and Lemma \ref{gpl}) for $\X_k\neq \D_4$.\\
\hf If $(\ff, \p) \in \D_4$ then there exists three distinct directions $\lp$ along which the third derivative vanishes. Hence, the projection map is three to one and orientation preserving. If we count the signed intersection number then we see that each transverse intersection point in $\D_4$ is accounted for thrice (with the same sign) when counted in $\PP \D_4$. This proves \eqref{up_down_equation} (using definition \ref{up_number_defn} in its transverse intersection form) for $\X_k=\D_4$. We note that this method of counting transverse intersection points also work identically for the first part above where the map $\PP \X_k \longrightarrow \X_k$ is one to one.\qed \\

\ni \textbf{Proof of Lemma \ref{cl1}:} It is evident that the left hand side 
of \eqref{cl1_equation1} is a subset of its right hand side, since all the sections are \textit{continuous}. For the converse, we will show that if $p$ belongs to the right hand side, then there exists a sequence $p_n$ in $\SS_{k-1}$ that converges to $p$. Since the section
$$ \ts_{k}: 
\{ p\in M: \ts_0(p) =0, \ldots, \ts_{k-1}(p) =0 \} \lra V_{k} $$
is transverse to the zero set, there exists a solution $p_t$ \textit{near} $p$ to the set of equations 
$$ \ts_0(p_t) =0, \ldots, \ts_{k-1}(p_t) =0, ~\ts_{k}(p_t) = t $$
if $t$ is sufficiently small. By definition, $p_t$ belongs to $\SS_{k-1}$ if $t\neq 0$. This  
gives us a sequence that lies in $\SS_{k-1}$  and converges to $p$. To finish the proof, it suffices to prove \eqref{cl1_equation2} by showing that 
\begin{align*}
 \ov{\SS}_{k} &= \{ p\in M: \ts_0(p)=0, \ldots, \ts_{k-1}(p)=0, ~\ts_{k}(p) =0 \} 
\end{align*}
This follows from an identical argument as before, using transversality of the section 
\bgd
\pushQED{\qed}
\ts_{k+1}: \{ p\in M: \ts_0(p) =0, \ldots, \ts_{k}(p) =0 \} \lra V_{k+1}.\qedhere 
\popQED
\edd

\vspace*{0.1cm}\textbf{Proof of Lemma \ref{cl2}:} We will show that the left hand side of \eqref{cl2_eqn1} is a subset of its right hand side. We may assume that $p \in \ov{\SS}_{k-1} - \SS_{k-1}$ and $\nts(p) \neq 0$. We claim that $p \in \ov{\SS}_{k}$. In other words, we need to show that there exists a sequence in $\SS_{k}$ that converges to $p$. 
Since the section
$$ \ts_{k+1}: 
\{ p\in M: \ts_0(p) =0, \ldots, ~\ts_{k}(p) =0, ~\nts(p) \neq 0  \} \lra V_{k+1} $$
is transverse to the zero set, there exists a solution $p_t$ \textit{near} 
$p$ to the set of equations 
$$\ts_0(p_t) =0, \ldots, \ts_{k}(p_t) =0, \ts_{k+1}(p_t) = t$$
if $t$ is sufficiently small.  
By definition, $p_t$ belongs to $\SS_{k}$ if $t\neq 0$. This gives us a sequence that lies in $\SS_{k}$  and converges to $p$. This proves that the left hand side of \eqref{cl2_eqn1} is a subset of its right hand side; if $p \in \ov{\SS}_{k-1} - \SS_{k-1}$ and $\nts(p) = 0$, then  $p \in \B$. Next, let us show that the right hand side of \eqref{cl2_eqn1} is a subset of its left hand side. Since $ \B \subset \ov{\SS}_{k-1}$, it suffices to show that $ \ov{\SS}_{k} \subset \ov{\SS}_{k-1} $.
We need to show that if $p \in \ov{\SS}_{k}$, then there exists a sequence in $\SS_{k-1}$ that converges to $p$. Since by hypothesis $p\in \ov{\SS}_k$, there exists a sequence $p_n \in \SS_k$ that converges to $p$, i.e., $p_n$ satisfies the equations: 
$$ \ts_0(p_{n}) =0, \ldots,\ts_{k-1}(p_{n}) =0, ~\ts_{k}(p_{n}) =0, ~\ts_{k+1}(p_{n}) \neq 0, ~\nts(p_{n}) \neq 0.$$
However, since the section 
$$ \ts_{k}: 
\{ p\in M: \ts_0(p) =0, \ldots, \ts_{k-1}(p) =0, \nts(p) \neq 0 
\} \lra V_k $$
is transverse to the zero set, we can conclude that there exists a solution $p_{n,t}$ near $p_{n}$ for the 
set of equations:
$$ \ts_0(p_{n, t}) =0, \ldots,\ts_{k-1}(p_{n, t}) =0, \ts_{k}(p_{n,t}) =t.$$
Moreover, since $\nts(p_n) \neq 0$ and $\nts$ is continuous, we get that $\nts(p_{n,t}) \neq 0$ if $t$ is sufficiently 
small. Hence, $p_{n,t}$ lies in $\SS_{k-1}$.
This gives us a sequence in $\SS_{k-1}$ that converges to $p$, which proves \eqref{cl2_eqn1}. 
Finally, we will prove \eqref{cl2_extra}. Let us prove that
\begin{align}
\label{cl2_extra_main_equation_subset} 
\ov{\SS}_0 \supset \{ p\in M: ~\ts_0(p) =0, ~\nts(p) \neq 0\}.
\end{align}
In other words, we need to show that if $\ts_0(p) =0$ and $\nts(p) \neq 0$, 
then there exists a sequence in $\SS$ that converges to $p$.
Since the section 
$$ \ts_{1}: \{ p\in M: \ts_0(p) =0, ~\nts(p) \neq 0 \} \lra V_1 $$
is transverse to the zero set, 
there exists a solution $p_t$ near $p$ to the set of equations 
$$ \ts_0(p_t) =0, ~\ts_{1}(p_t) =t.$$ 
Since $\nts$ is continuous and $\nts(p) \neq 0$, we get that $\nts(p_t) \neq 0$ if $t$ is sufficiently small. This gives us the desired sequence in $\SS_0$. Now using \eqref{cl2_extra_main_equation_subset} and the definition of $\SS_{-1}$ and $\B$, we get that
\bgd
M=\{p \in M: \ts_0(p) \neq 0, \nts(p) \neq 0 \} \cup \{ p\in M: \ts_0(p) =0, \nts(p) \neq 0\} \cup \{ p\in M: \nts(p) =0\}                                   
\edd
is contained in $\SS_{-1} \cup \ov{\SS}_0 \cup \B^{\prime}$. The reverse inclusion is vacuous.  \qed\\

\textbf{Proof of Lemma \ref{cl3}:} Observe that the left hand side of \eqref{cl3_eqn1} is a subset of its right hand side. For the converse, assume that $p$ belongs to the right hand side. We need to show that there exists a sequence in $\SS_0$ that converges to $p$. Let us assume that 
$\nts(p) \neq 0$. 
Since the section  
$$ \ts_1 : \{ p \in M: ~\ts_0(p) =0, ~\nts(p) \neq 0 \} \lra V_1$$
is transverse to the zero set, there exists a solution $p_t$ near $p$ to the set of equations $\ts_0(p_t) =0$ and $\ts_1(p_t) = t$. Moreover, since $\nts(p) \neq 0$, we conclude that $\nts(p_t) \neq 0$ if $t$ is small. Hence, $p_t \in \SS_0$ for all small $t$ which gives us the desired sequence. Next let us assume that $\nts(p) =0$. Since the section 
$$ \nts: \{ p \in M: ~\ts_0(p) =0 \} \lra W $$
is transverse to the zero set,  there exists a solution $p_{t_1}$ near $p$  to the set of equations 
$$ \ts_0(p_{t_1}) =0, ~\nts(p_{t_1}) = t_1.$$
Since $\nts(p_{t_1}) \neq 0$, there exists a solution $p_{t_1, t_2}$ near $p_{t_1}$ to the set of equations $\ts_0(p_{t_1,t_2}) =0$ and $\ts_1(p_{t_1, t_2}) = t_2$ because the section   
$$ \ts_1 : \{ p \in M: ~\ts_0(p) =0, ~\nts(p) \neq 0 \} \lra V_1$$
is transverse to the zero set. Moreover, since $\nts(p_{t_1}) \neq 0$, we infer that 
$\nts(p_{t_1, t_2}) \neq 0$ if $t_2$ is sufficiently small. Hence, $p_{t_1, t_2}$ lies in $\SS_0$. This  gives us the desired sequence which proves Lemma \ref{cl3}. \qed  

\subsection{Chern classes and projectivized bundle} 

\begin{lmm}{\bf (\cite{MiSt}, Theorem 14.10)}
\label{total_chern_class_of_tpn}
Let $\P^n$ be the $n$-dimensional complex projective space and $\gamma \lra \P^n$ the tautological line bundle. 
Then the total Chern class $c(T\P^n)$ is given by $~c(T\P^n) = (1+ c_1(\gamma^*))^{n+1}.$
\end{lmm}

\begin{lmm}{\bf (\cite{BoTu}, pp. 270)}
\label{cohomology_ring_of_pv}
Let $V\lra M$ be a complex vector bundle, of rank $k$, over a smooth manifold $M$ and $\pi: \P V \lra M$ the \textit{projectivization} of 
$V$. Let $\G\lra \P V$ be the tautological line bundle over $\P V$ and $\lm = c_1(\G^*)$. 
There is a linear isomorphism
\begin{align}
H^*(\P V; \mathbb{Z}) & \cong  H^*(M;\mathbb{Z}) \otimes H^*(\P^{k-1};\mathbb{Z})
\end{align}
and an isomorphism of rings
\begin{align}
H^*(\P V; \mathbb{Z})& \cong  \frac{H^*(M;\mathbb{Z})[\lm]}{\lan \lm^k + \lm^{k-1} \pi^*c_1(V) + \lm^{k-2} \pi^*c_2(V) + \ldots +\pi^*c_{k}(V)\ran}.
\end{align}
In particular, if $\omega \in H^{*}(M; \mathbb{Z}) $ is a top cohomology class then 
\bgd
\lan \pi^*(\omega) \lm^{k-1}, [\P V] \ran = \lan \omega, [M] \ran,
\edd
i.e., $\lm^{k-1}$ is a \textit{cohomology extension} of the fibre.
\end{lmm} 

%
%

\bibliographystyle{siam}
\bibliography{Myref_bib.bib}

\vspace*{0.4cm}

\hf {\small D}{\scriptsize EPARTMENT OF }{\small M}{\scriptsize ATHEMATICAL }{\small S}{\scriptsize CIENCES, }{\small B}{\scriptsize INGHAMTON }{\small U}{\scriptsize NIVERSITY, }{\small NY} {\footnotesize 13902-6000, }{\small USA}\\
\hf{\it E-mail address} : \texttt{somnath@math.binghamton.edu}\\

\hf {\small D}{\scriptsize EPARTMENT OF }{\small M}{\scriptsize ATHEMATICS, }{\small I}{\scriptsize NSTITUTE OF }{\small M}{\scriptsize ATHEMATICAL }{\small S}{\scriptsize CIENCES, }{\small C}{\scriptsize HENNAI }{\footnotesize 600113, }{\small INDIA}\\
\hf{\it E-mail address} : \texttt{ritwikm@imsc.res.in}\\[0.2cm]

\end{document}